\documentclass{amsart}

\usepackage{amssymb}
\usepackage{amsmath}
\usepackage{amsthm}
\usepackage{amscd}
\usepackage{bbm}
\usepackage{changepage}
\usepackage{enumerate}
\usepackage{fancybox}
\usepackage{fancyhdr}
\usepackage{float}
\usepackage[bookmarksopen,bookmarksdepth=3]{hyperref}
\usepackage{marvosym}
\usepackage{mathrsfs}
\usepackage[pdftex]{graphicx}
\usepackage{setspace}
\usepackage{subfig}
\usepackage{theoremref}
\usepackage{verbatim}
\usepackage{xy}

\def\Xint#1{\mathchoice
   {\XXint\displaystyle\textstyle{#1}}%
   {\XXint\textstyle\scriptstyle{#1}}%
   {\XXint\scriptstyle\scriptscriptstyle{#1}}%
   {\XXint\scriptscriptstyle\scriptscriptstyle{#1}}%
   \!\int}
\def\XXint#1#2#3{{\setbox0=\hbox{$#1{#2#3}{\int}$}
     \vcenter{\hbox{$#2#3$}}\kern-.5\wd0}}

\def\dashint{\Xint-}

\hypersetup{
    colorlinks,
    citecolor=green,
    linkcolor=red
  }

\makeatletter
\def\blfootnote{\xdef\@thefnmark{}\@footnotetext}
\makeatother

\newcommand{\xx}{\times}
\newcommand{\C}{\mathbb{C}}
\newcommand{\R}{\mathbb{R}}
\newcommand{\N}{\mathbb{N}} 
\newcommand{\br}[1]{\left( #1 \right)}
\newcommand{\brs}[1]{\left[ #1 \right]}
\newcommand{\norm}[1]{\left\Vert #1 \right\Vert}
\newcommand{\abs}[1]{\left\vert #1 \right\vert}
\newcommand{\lb}[0]{\left\lbrace}
\newcommand{\rb}[0]{\right\rbrace}
\newcommand\BoldSquare{%
  \setlength\fboxrule{1.1pt}\setlength\fboxsep{0pt}\fbox{\phantom{\rule{5pt}{5pt}}}}

\newenvironment{prof}{   \textsc{Proof.}}{\hfill \BoldSquare \vspace{5pt} }

\newtheorem{thm}{Theorem}[section]
\newtheorem{prop}{Proposition}[section]
\newtheorem{lem}{Lemma}[section]
\newtheorem{cor}{Corollary}[section]
\newtheorem{conj}{Conjecture}[section]

\theoremstyle{definition}

\newtheorem{deff}{Definition}[section]

\newtheorem{rmk}{Remark}[section]

\title[Exponential Weights for Schr\"{o}dinger Operators]{Weights of
  Exponential Growth and Decay for Schr\"{o}dinger-Type Operators}
\author{Julian Bailey}

\date{}

\begin{document}

\begin{abstract}
Fix $d \geq 3$ and $1 < p < \infty$. Let $V : \mathbb{R}^{d} \rightarrow [0,\infty)$ belong to the reverse
H\"{o}lder class $RH_{d/2}$ and consider the Schr\"{o}dinger
operator $L_{V} := - \Delta + V$.  In this article, we introduce
classes of weights $w$ for which the Riesz transforms $\nabla
L_{V}^{-1/2}$, their adjoints $L_{V}^{-1/2} \nabla$ and the heat maximal operator $\sup_{t >
  0} e^{- t L_{V}} |f|$ are bounded on the weighted Lebesgue space
$L^{p}(w)$. The boundedness
of the $L_{V}$-Riesz potentials $L_{V}^{-\alpha/2}$ from $L^{p}(w)$ to
$L^{\nu}(w^{\nu/p})$ for $0 < \alpha \leq 2$ and $\frac{1}{\nu} =
\frac{1}{p} - \frac{\alpha}{d}$ will also be
proved. These 
weight classes are strictly larger than a class previously introduced by
B. Bongioanni, E. Harboure and O. Salinas in
\cite{bongioanni2011classes} that shares these properties and they contain weights of exponential growth and decay.

\vspace*{0.1in}

The classes will also be considered in relation to different generalised forms of Schr\"{o}dinger
operator. In particular, the Schr\"{o}dinger operator with
measure potential $- \Delta + \mu$, the uniformly elliptic operator with potential
$- \mathrm{div} A \nabla + V$ and the magnetic
Schr\"{o}dinger operator $(\nabla - i a)^{2} + V$ will all be considered. It will
be proved that, under suitable conditions, the standard operators corresponding
to these second-order differential operators
are bounded on $L^{p}(w)$ for weights $w$ in these classes.
  \end{abstract}

  \maketitle

  \setcounter{tocdepth}{1}
  \tableofcontents

\section{Introduction}
\label{sec:Introduction}

\blfootnote{\textit{Key words and phrases.}
  Schr\"{o}dinger operator; weights; Riesz transforms. \\
  Mathematics Subject Classification. 42B20, 42B25, 42B35. \\
This research was supported by the Engineering and Physical
Sciences Research Council (EPSRC) grant EP/P009239/1.}

The Laplacian operator $\Delta$ is inextricably linked to classical weighted
theory. The characterisation of the Muckenhoupt class $A_{p}$ for $1 <
p < \infty$ in terms
of various operators related to the Laplacian, such as the Riesz
transforms $R_{0} := \nabla (-\Delta)^{-\frac{1}{2}}$, the heat
maximal operator $T^{*}_{0}f := \sup_{t > 0} e^{t \Delta} \abs{f}$ and
the Hardy-Littlewood maximal operator $M_{0}$, is widely regarded as the apex
of classical weighted theory. This string of results, the culmination
of work from a number of great mathematicians in the 1970's, states that a weight $w$ will be contained in $A_{p}$ if and
only if any one of these operators is bounded from the weighted
$p$-Lebesgue space $L^{p}(w) :=
L^{p}(\R^{d}; w \, dx)$ to itself (see \cite{muckenhoupt1974equivalence}, \cite{hunt1973weighted} and
\cite{coifman1974weighted} for sufficiency and \cite{hunt1973weighted} and
\cite{stein1993harmonic} for necessity).

A current area of active research is the study of the harmonic analysis of
differential operators other than the Laplacian. A natural question to
ask in such a setting is whether it is possible to construct a
Muckenhoupt-type class adapted to the underlying differential operator
and whether this class can be characterised in terms of the
corresponding operators such as the associated Riesz transforms and heat maximal operator. In this article,  the
differential operators of interest are Schr\"{o}dinger operators.

Fix
dimension $d \geq 3$. The Schr\"{o}dinger operator with non-negative
potential $V \in L^{1}_{loc}\br{\R^{d}}$ is defined to be the operator
$$
L_{V} :=  - \Delta + V.
$$
Let $R_{V}$ and $R^{*}_{V}$ denote the Riesz transforms associated with $L_{V}$ and
their adjoints,
$$
R_{V} := \nabla L_{V}^{-\frac{1}{2}} \quad and \quad R^{*}_{V} := L_{V}^{-\frac{1}{2}} \nabla.
$$
Define the $L_{V}$-Riesz potential for $0 < \alpha \leq 2$ and heat
maximal operator for $L_{V}$ respectively through
$$
I^{\alpha}_{V} := L_{V}^{-\frac{\alpha}{2}} \quad and \quad T^{*}_{V}f(x) :=
\sup_{t > 0} e^{- t L_{V}} \abs{f}(x),
$$
for $f \in L^{1}_{loc}\br{\R^{d}}$ and $x \in \R^{d}$. 
Recall the definition of the reverse H\"{o}lder class of
potentials. A non-negative function $V \in L^{1}_{loc}\br{\R^{d}}$ is said to belong to the class $RH_{q}$ for $1 < q < \infty$ if there
exists a constant $C > 0$ for which
$$
\br{\frac{1}{\abs{B}}\int_{B} V^{q}}^{\frac{1}{q}} \leq C \br{\frac{1}{\abs{B}}
\int_{B} V}
$$
for all open Euclidean balls $B \subset \R^{d}$. The smallest constant $C$ for which
the above estimate is satisfied will be denoted $\brs{V}_{RH_{q}}$.
In the article \cite{shen1995lp},
Z. Shen introduced the notion of the critical radius function
$\rho_{V}$ for a potential $V \in RH_{\frac{d}{2}}$. This is the function $\rho_{V} : \R^{d} \rightarrow
[0,\infty)$ defined through
\begin{equation}
  \label{eqtn:CritRad}
\rho_{V}(x) := \sup \lb r > 0 : \frac{1}{r^{d - 2}}
\int_{B(x,r)} V(x) \, dx \leq 1 \rb
\end{equation}
for $x \in \R^{d}$, where the notation $B(x,r)$ is used to denote the
open Euclidean
ball of radius $r$ centered at $x$. The critical radius function determines a scale
below which the operators associated with the Schr\"{o}dinger
operator behave locally like their classical counterparts. This allowed
for the analysis of operators such as $R_{V}$ to be split up into two
regions; the local region where $R_{V}$ resembles the
classical Riesz transforms $R_{0}$ and a
global region where the singular kernel for $R_{V}$ will have substantially
better decay properties than the kernel of $R_{0}$.

In the article \cite{bongioanni2011classes}, the authors
B. Bongioanni, E. Harboure and O. Salinas initiated a study into the
weighted theory of Schr\"{o}dinger operators for $V \in RH_{\frac{d}{2}}$ by introducing a new class of
weights, denoted by $A_{p}^{V,\infty}$, that
was adapted to $L_{V}$ and whose definition was based on
the Schr\"{o}dinger operator machinery of Shen. 
What made the class $A_{p}^{V, \infty}$ so compelling was that not only were the
operators $R_{V}$, $R^{*}_{V}$ and $T^{*}_{V}$ all bounded
on $L^{p}(w)$ for any $w \in A_{p}^{V,\infty}$, but also that the class was strictly larger
than the classical Muckenhoupt class $A_{p}$. In addition, in
\cite{bongioanni2011classes} it was also
proved that the fractional
integral operator $I^{\alpha}_{V}$ is bounded on appropriate weighted Lebesgue spaces that correspond to
the classical case for weights in $A_{p}^{V,\infty}$.

Unfortunately, an inescapable deficiency with the class
$A_{p}^{V,\infty}$ is that the reverse implication does not hold in
general.  That is, there exist weights not contained in
$A_{p}^{V,\infty}$ for which $R_{V}$ and $T^{*}_{V}$ are bounded on
$L^{p}(w)$. 
 Such an example was found in 
\cite{bailey2018hardy} for the harmonic oscillator potential,
$V(x) = \abs{x}^{2}$. Indeed for polynomial potentials of order zero
or higher, as will be shown in this article, there
exist a wealth of such counterexamples that are non-doubling weights of
exponential growth or decay.
The existence of such weights demonstrates that the class
$A_{p}^{V,\infty}$ is not characterised completely by the boundedness of $R_{V}$
or $T^{*}_{V}$
and thus the harmonic analytic aspects of $L_{V}$ are not fully
captured by this class.
The first aim of this article is thus to improve upon
the class $A_{p}^{V,\infty}$ by introducing a strictly larger class
that accounts for the counterexamples of exponential-type weights.

 Let $d_{V}(x,y)$ denote the Agmon distance for the
potential defined through
$$
d_{V}(x,y) := \inf_{\gamma} \int^{1}_{0} \rho_{V}(\gamma(t))^{-1}
\abs{\gamma^{'}(t)} \, dt,
$$
where the infimum is taken over all curves $\gamma : [0,1] \rightarrow
\R^{d}$
connecting the points $x, \, y \in \R^{d}$. We introduce the notation
$B_{V}(x,r)$ to denote the open ball of
radius $r > 0$ centered at the point $x \in \R^{d}$ in the metric
$d_{V}$,
$$
B_{V}(x,r) := \lb y \in \R^{d} : d_{V}(x,y) < r \rb.
$$
 Our weight class is defined as follows.

\begin{deff} 
  \label{def:ApVc}
   Let $1 < p < \infty$ and $c > 0$. 
 $S_{p,c}^{V}$ is the class
 of all weights for which
$$
\brs{w}_{S_{p,c}^{V}} := \sup_{B_{V}} \br{\frac{1}{\abs{B_{V}}e^{c \cdot r}}
  \int_{B_{V}} w}^{\frac{1}{p}} \br{\frac{1}{\abs{B_{V}}e^{c \cdot r}} \int_{B_{V}}
  w^{-\frac{1}{p-1}}}^{\frac{p - 1}{p}} < \infty,
$$
where the supremum is taken over all balls $B_{V} = B_{V}(x,r) \subset \R^{d}$
in the metric $d_{V}$ with $x \in \R^{d}$ and $r > 0$.
\end{deff}

One of the main results for this article is the below theorem whose
proof will be provided in Section \ref{subsec:Riesz}.

\begin{thm} 
 \label{thm:Riesz} 
 Let $V \in RH_{q}$ for some $q > \frac{d}{2}$. The following
 statements are true.

\vspace*{0.1in}
 
 \begin{enumerate}[(i)]
\item \label{Riesz1} Suppose that $q \geq d$. There must exist some $c_{1} > 0$ for which both $R_{V}$ and $R_{V}^{*}$
  are bounded on $L^{p}(w)$ for $1 < p < \infty$ when $w \in S_{p,c_{1}}^{V}$.

  \vspace*{0.1in}

  \item \label{Riesz2} Suppose instead that $\frac{d}{2} < q < d$ and let $s$ be
    defined through $\frac{1}{s} = \frac{1}{q} - \frac{1}{d}$. There
    exists a constant $c_{2} > 0$ for which the operator
    $R_{V}^{*}$ is bounded on $L^{p}(w)$ for $s' < p < \infty$ when $w
    \in S_{p/s',c_{2}}^{V}$ and the operator $R_{V}$ is bounded on $L^{p}(w)$
    for $1 < p < s$ when $w^{-\frac{1}{p - 1}} \in S_{p'/s',c_{2}}^{V}$.

\vspace*{0.1in}
    
    \item \label{Riesz3} For any $q > \frac{d}{2}$ and $0 < \alpha
      \leq 2$, there exists $c_{3} > 0$ for which the operator
      $I^{\alpha}_{V}$ is bounded from $L^{p}(w)$ to
$L^{\nu}(w^{\nu/p})$ for $w^{\nu/p} \in S_{1 +
  \frac{\nu}{p'},c_{3}}^{V}$ and $1 < p < \frac{d}{\alpha}$,
where $\frac{1}{\nu} = \frac{1}{p} - \frac{\alpha}{d}$.
\end{enumerate}
In each of above statements, the constants $c_{1}, \, c_{2}$ and
$c_{3}$ will depend on $V$ only through $[V]_{RH_{\frac{d}{2}}}$ and they will
be independent of $p$.
\end{thm}

 This theorem
provides an improvement upon all existing weight classes that possess
these properties since, as will be proved in Section \ref{sec:Adapted},
$A_{p}^{V,\infty}$ is strictly contained in $S_{p,c}^{V}$ for any $c >
0$. Indeed, our class will allow for exponential growth and decay in
the weights as opposed to only polynomial growth and decay as in
$A_{p}^{V,\infty}$. Our class also has the advantage that it has a more
natural geometric definition than other classes since it is defined in
terms of the inherent geometry associated with the potential. 

Let $\Gamma_{V}$ denote the fundamental solution of the
Schr\"{o}dinger operator $L_{V}$. This is a function defined on $\lb
(x,y) \in \R^{d} \xx \R^{d} : x \neq y \rb$ with the properties that
$\Gamma_{V} (\cdot, y) \in L^{1}_{loc}\br{\R^{d}}$ and $L_{V}
\Gamma_{V}(\cdot,y) = \delta_{y}$ for each $y \in \R^{d}$, where
$\delta_{y}$ is the Dirac delta distribution with pole at $y$. The proof of Theorem \ref{thm:Riesz} relies heavily on sharp pointwise exponential decay
estimates for $\Gamma_{V}$ that were obtained by Shen in the article
\cite{shen1999fundamental}. These estimates state that there must exist
some $C_{1}, \, C_{2}, \, \varepsilon_{1}, \, \varepsilon_{2} > 0$ for which
\begin{equation}
  \label{eqtn:ShenSharp}
C_{1} \frac{e^{- \varepsilon_{1} d_{V}(x,y)}}{\abs{x - y}^{d - 2}} \leq
\Gamma_{V}(x,y) \leq C_{2}
\frac{e^{-\varepsilon_{2}d_{V}(x,y)}}{\abs{x - y}^{d - 2}} \quad
\forall \ x, \, y \in \R^{d}.
\end{equation}
It should be stressed that only the upper estimate is required for the
proof of Theorem \ref{thm:Riesz}. However, the lower estimate is
useful in that it hints towards the optimality of these weight classes
in a way that will be stated explicitly in Conjecture
\ref{conj:Upper}.

 One should also expect for the heat maximal operator to be
bounded on $L^{p}(w)$ for weights in the class $S_{p,c}^{V}$.
Unfortunately, it is still an open problem as to whether the heat
kernel of $L_{V}$ satisfies sharp estimates that are analogous to
\eqref{eqtn:ShenSharp}. The best known estimates were obtained by K. Kurata
in \cite{kurata2000estimate} and are non-sharp.
This indicates that the proof of the boundedness of
$T_{V}^{*}$ on $L^{p}(w)$ for $w \in S_{p,c}^{V}$ will be just beyond our
reach and will remain so until sharp estimates for the heat kernel are
proved. Therefore, instead of proving the boundedness of $T^{*}_{V}$ on
$L^{p}(w)$ for $w \in S_{p,c}^{V}$, a second smaller class of
exponential-type weights will be introduced and boundedness will be
proved for this class instead.

\begin{deff} 
  \label{def:HeatWeights}
 Let $1 < p < \infty$ and $m, \, c > 0$. Introduce the notation
 $\Phi_{m,c}^{V}$ to denote the function
$$
\Phi_{m,c}^{V}(x,r) := \exp \br{c \br{1 + \frac{r}{\rho_{V}(x)}}^{m}},
$$
for $x \in \R^{d}$ and $r > 0$.
 Let $H_{p,c}^{V,m}$ denote the class
 of all weights for which
 $$
\brs{w}_{H_{p,c}^{V,m}} := \sup_{B} \br{\frac{1}{\abs{B} \Phi_{m,c}^{V}(x,r)}
\int_{B} w}^{\frac{1}{p}} \br{\frac{1}{\abs{B} \Phi_{m,c}^{V}(x,r)} \int_{B}
w^{-\frac{1}{p-1}}}^{\frac{p - 1}{p}} < \infty,
$$
where the supremum is taken over all Euclidean balls $B = B(x,r) \subset \R^{d}$
with radius $r > 0$ and center $x \in \R^{d}$.
\end{deff}

The class $H_{p,c}^{V,m}$ is very similar in nature to
$A_{p}^{V,\infty}$ except that it will allow for
exponential growth and decay (refer to Section \ref{sec:Adapted} for a
rigorous definition of $A_{p}^{V,\infty}$). Due to the similarities between the two classes,
it is not difficult to see that $A_{p}^{V,\infty} \subset H_{p,c}^{V,m}$ for any
$c, \, m > 0$.
An inclusion that is less obvious is that for any $c > 0$ we must have
$H_{p,c}^{V,m} \subset S_{p,c}^{V}$ provided that $m \leq 2 m_{0}$ for
some $m_{0} > 0$ that will be defined later in the article. This
statement will be proved in Section \ref{sec:Adapted}.

One of our main results for this
class of weights is the following theorem that will be proved in
Section \ref{subsec:Heat}.

\begin{thm} 
 \label{thm:Heat} 
 Suppose that $V \in RH_{\frac{d}{2}}$ and let $1 < p < \infty$. There
 must exist $c, \, m_{0} > 0$, independent of $p$,
 such that for any $w \in
 H_{p,c}^{V,m}$ with $m \leq m_{0}$ the operator $T^{*}_{V}$ is
 bounded on $L^{p}(w)$. 
\end{thm}

 The second aim of this article is to extend the results for these
 freshly minted classes of weights to more general forms of
 Schr\"{o}dinger operator. We will consider three different
 types of generalised Schr\"{o}dinger operator and prove
 statements analogous to Theorems \ref{thm:Riesz} and \ref{thm:Heat}
 for each of these forms. The first generalised form of
 Schr\"{o}dinger operator that will be considered is the
 Schr\"{o}dinger operator with measure potential. This is the operator
 defined through
 $$
L_{\mu} := -\Delta + \mu,
$$
where $\mu$ is a non-negative Radon measure on $\R^{d}$ that satisfies the
conditions
$$
\mu(B(x,r)) \leq C_{\mu} \br{\frac{r}{R}}^{d - 2 + \delta_{\mu}} \mu(B(x,R))
$$
and
$$
\mu(B(x,2 r)) \leq D_{\mu} \br{\mu(B(x,r)) + r^{d - 2}}
$$
for all $x \in \R^{d}$ and $0 < r < R$, for some $\delta_{\mu}, \,
C_{\mu}, \, D_{\mu} > 0$. In Section \ref{sec:Schrodinger}, it will be
proved that the $L_{\mu}$-counterpart of Theorem \ref{thm:Riesz} is true.

Next, we
will consider uniformly elliptic operators with potential. Let $A$ be
a $d \xx d$ matrix-valued function with real-valued coefficients in
$L^{\infty}\br{\R^{d}}$. Suppose that $A$ satisfies the ellipticity
condition,
$$
\lambda \abs{\xi}^{2}\leq \langle A(x) \xi, \xi \rangle \leq \Lambda \abs{\xi}^{2}
$$
for some $\lambda, \, \Lambda > 0$, for all $\xi \in \R^{d}$ and almost every $x \in \R^{d}$.
Define for $V \in RH_{\frac{d}{2}}$ the operator
$$
L_{A,V} := - \mathrm{div} A \nabla + V.
$$
In Section \ref{sec:UniformlyElliptic} it will be proved that
if the critical radius function $\rho_{V}$ is bounded from above, as
is the case for any polynomial potential of order zero or higher,
and $A$ is H\"{o}lder continuous then
the first and second parts of Theorem \ref{thm:Riesz} will be true for
$L_{A,V}$.
It will also be proved that the $L_{A,V}$-counterpart of
Theorem \ref{thm:Heat} is true subject
to the constraint $A = A^{*}$. The third part of Theorem
\ref{thm:Riesz} will be proved subject to no additional
constraints. For a rigorous statement of these results refer to
Section \ref{sec:UniformlyElliptic}.

\begin{rmk}
  \label{rmk:PerturbationEffect}
 Our result that $R_{A,V} := \nabla L_{A,V}^{-\frac{1}{2}}$ is bounded on $L^{p}\br{w}$ for all
  weights contained in a class strictly larger than the
  classical Muckenhoupt class $A_{p}$ when $V \in RH_{\frac{d}{2}}$
  and $\rho_{V}$ is bounded from above sharply contrasts with the
  potential free case. Without the presence of the potential, it is well-known
  known that there exists $A$ for which the
  operator $R_{A,0} := \nabla L_{A,0}^{-\frac{1}{2}}$ is unbounded
  on $L^{p}(w)$ for some $w \in A_{p}$. This follows by combining
  Remarks $1.7$ and $1.8$ of \cite{shen2005bounds} for example. Thus,
  the perturbation $A$ has the tendency to decrease the size of the
  associated weight class. The inclusion of a sufficiently large
  potential $V$ will have the opposite effect and increase the
  size of the weight class, effectively cancelling out the presence of $A$.
\end{rmk}

Finally, for $a = (a_{1},\cdots,a_{d})$  a vector of real-valued
functions in $C^{1}\br{\R^{d}}$, we will consider the magnetic Schr\"{o}dinger operator
$$
L_{V}^{a} := \br{\nabla - i a}^{*} \br{\nabla - i a} + V.
$$
In Section \ref{sec:Magnetic} it will
be proved that subject to additional constraints on $a$ and $V$ the third part of
Theorem \ref{thm:Riesz} and Theorem \ref{thm:Heat} will hold. A weaker
form of the first part of Theorem \ref{thm:Riesz} will also be
proved.

In the article \cite{mayboroda2019exponential}, the authors S.
Mayboroda and B. Poggi proved that Shen's sharp exponential decay
estimates on the fundamental solution \eqref{eqtn:ShenSharp} could be generalised to the operators
$L_{A,V}$ and $L^{a}_{V}$. Our weighted results for both the magnetic Schr\"{o}dinger
operator and the uniformly elliptic operator with potential will rely
on this result. In addition, in order to prove the boundedness of the
heat maximal operator for $L_{A,V}$ and $L^{a}_{V}$, we will also
require the exponential decay estimates for the heat kernel provided
by \cite{kurata2000estimate}.

\vspace*{0.1in}

In the last part of this article, we will consider necessary
conditions for a weight to satisfy in order for the operators $R_{V}$
and $T^{*}_{V}$ to be bounded on $L^{p}(w)$. It will first be proved
that for $V \in RH_{d}$, if $R_{V}$ is bounded on $L^{p}(w)$ then the
weight must be contained in the local Muckenhoupt class
$A_{p}^{V,loc}$ (refer to Section \ref{sec:Adapted} for the definition
of $A_{p}^{V,loc}$). Following this, we will discuss
the optimality of the classes $S_{p,c}^{V}$ and $H_{p,c}^{V,m}$.  In
particular, the following conjecture will be considered.

\begin{conj}
  \label{conj:Upper}
  Let $V \in RH_{d}$. There exists $c_{1}, \, c_{2} > 0$
  for which the following chains of inclusions hold,
  $$
S_{p,c_{1}}^{V} \subset \lb w : \norm{R_{V}}_{L^{p}(w)} < \infty \rb
\subset S_{p,c_{2}}^{V}
$$
and
  $$
S_{p,c_{1}}^{V} \subset \lb w : \norm{T^{*}_{V}}_{L^{p}(w)} < \infty \rb
\subset S_{p,c_{2}}^{V}.
$$
\end{conj}

 In Section \ref{sec:Optimality}, it will be proved that the first
 chain of inclusions of the previous conjecture is true
 for constant potentials and the second chain of inclusions is true
 for potentials that are bounded both from above and below.

 It will be proved in Section \ref{sec:Adapted} that for any $c >
 0$ there exists $c_{1}, \, c_{2}, \, m_{1}, \, m_{2} > 0$ for which
 $H_{p,c_{1}}^{V,m_{1}} \subset S_{p,c}^{V} \subset
 H_{p,c_{2}}^{V,m_{2}}$. The following is therefore a weaker form of
 Conjecture \ref{conj:Upper}.

 \begin{conj}
   \label{conj:Weak}
Let $V \in RH_{d}$. There exists $c_{1}, \, c_{2}, \, m_{1},
\, m_{2} > 0$ for which the following chains of inclusions hold,
$$
H_{p,c_{1}}^{V,m_{1}} \subset \lb w : \norm{R_{V}}_{L^{p}(w)} < \infty
\rb \subset H_{p,c_{2}}^{V,m_{2}}
$$
and
$$
H_{p,c_{1}}^{V,m_{1}} \subset \lb w : \norm{T^{*}_{V}}_{L^{p}(w)} <
\infty \rb \subset H_{p,c_{2}}^{V,m_{2}}.
$$
\end{conj}

In Section \ref{sec:Optimality}, the second chain of inclusions in
Conjecture \ref{conj:Weak} will
be proved for the harmonic oscillator potential $V(x) = \abs{x}^{2}$.

\subsection{Acknowledgements} I am indebted to Maria Carmen Reguera
for her sage advice, support and feedback on this project. I am also
very grateful to Pierre Portal for the numerous conversations on the
weighted theory of Schr\"{o}dinger operators that no doubt inspired this article.

\section{Critical Radius Function and Agmon Distance}

Throughout this article, the notation $A \lesssim B$ will be used to denote that there exists a
constant $c > 0$ for which the inequality $A \leq c B$ is
satisfied. Similarly, the notation $A \simeq B$ will denote that there
exists $c > 0$ for which $c^{-1} B \leq A \leq c B$. The dependence of the constant $c$ on the various
parameters should be clear from the context. To emphasise the
dependence of the constant on a particular parameter subscript
notation will be used. For example, $A \lesssim_{b} B$ will indicate that
the implicit constant $c$ depends on $b$. For a weight $w$ on
$\R^{d}$ and measurable set $E \subset \R^{d}$, the below notation
will frequently be used
$$
w(E) := \int_{E} w(x) \, dx.
$$

In this section, the function $\rho_{V}$ from \eqref{eqtn:CritRad} will be generalised and the
inherent geometry attached to such a function will be discussed.

\begin{deff} 
 \label{def:CriticalRadius} 
 A function $\rho : \R^{d} \rightarrow [0,\infty)$ will be called a
 critical radius function if there exist constants $k_{0}, \, B_{0} >
 1$ for which
 \begin{equation}
    \label{eqtn:Shen0}
    B_{0}^{-1} \rho(x) \br{1 + \frac{\abs{x -
          y}}{\rho(x)}}^{- k_{0}} \leq \rho(y) \leq B_{0}
    \rho(x) \br{1 + \frac{\abs{x -
          y}}{\rho(x)}}^{\frac{k_{0}}{k_{0} + 1}}.
  \end{equation}
  for all $x, \, y \in \R^{d}$.
\end{deff}

The following result from \cite{shen1995lp} confirms that the above
definition is indeed a generalisation of the function $\rho_{V}$.

\begin{lem}[{\cite[Lem.~1.4]{shen1995lp}}] 
 \label{lem:Shen0} 
 For $V \in RH_{\frac{d}{2}}$, the function $\rho_{V}$, as defined in \eqref{eqtn:CritRad},  is a
 critical radius function in the sense of Definition
 \ref{def:CriticalRadius}. In particular, \eqref{eqtn:Shen0} is satisfied with constants $k_{0}$ and $B_{0}$ that
 depend on $V$ only through $\brs{V}_{RH_{\frac{d}{2}}}$.
\end{lem}

For any critical radius function $\rho : \R^{d} \rightarrow
[0,\infty)$, it is possible to construct
  a corresponding Agmon distance through
  $$
d_{\rho}(x,y) := \inf_{\gamma} \int^{1}_{0} \rho(\gamma(t))^{-1}
\abs{\gamma'(t)} \, dt,
$$
where the infimum is taken over all possible curves $\gamma : [0,1]
\rightarrow \R^{d}$ connecting the points $x, \, y \in \R^{d}$.
At a local scale, the Agmon distance $d_{\rho}$ will be comparable to the
Euclidean distance. 

\begin{lem}
  \label{lem:LocalAgmon}
  Let $\rho : \R^{d} \rightarrow [0,\infty)$ be a critical radius
  function. There exists $D_{0} > 1$ so that for any $x, \, y \in \R^{d}$ with
  $\abs{x - y} \leq 2 \rho(x)$ we have
  $$
D^{-1}_{0} \abs{x - y} \rho(x)^{-1} \leq d_{\rho}(x,y) \leq D_{0} \abs{x - y} \rho(x)^{-1}.
  $$
\end{lem}

\begin{prof}
  Fix $x, \, y \in \R^{d}$ with $\abs{x - y} \leq 2 \rho(x)$.
Let's first prove the upper estimate. Let $\gamma$ denote the straight
line starting at $x$ and ending at $y$. Then, on applying \eqref{eqtn:Shen0},
\begin{align*}\begin{split}  
 \int^{1}_{0} \rho(\gamma(t))^{-1} \abs{\gamma'(t)} \, dt &\leq B_{0}
 \rho(x)^{-1} \int^{1}_{0} \br{1 + \frac{\abs{x -
       \gamma(t)}}{\rho(x)}}^{k_{0}} \abs{\gamma'(t)} \, dt \\
 &\leq B_{0} \rho(x)^{-1} 3^{k_{0}} \int^{1}_{0} \abs{\gamma'(t)}
 \, dt \\
 &= B_{0} \rho(x)^{-1} 3^{k_{0}}\abs{x - y},
\end{split}\end{align*}
where the second to last line follows from the fact that $\abs{x -
  \gamma(t)} \leq 2 \rho(x)$ for all $t \in [0,1]$.

\vspace*{0.1in}

Let's now prove the lower bound. Let $\gamma$ be a curve for which
$$
d_{\rho}(x,y) \geq \frac{1}{2} \int^{1}_{0} \rho(\gamma(t))^{-1}
\abs{\gamma'(t)} \, dt.
$$
The definition of a critical radius function leads to
$$
 d_{\rho}(x,y) \geq \frac{1}{2} B_{0}^{-1} \rho(x)^{-1} \int^{1}_{0}
 \br{1 + \frac{\abs{x - \gamma(t)}}{\rho(x)}}^{-\frac{k_{0}}{k_{0} +
     1}} \abs{\gamma'(t)} \, dt.
 $$
First suppose that the curve is entirely contained within $B(x,2
\rho(x))$. Then
\begin{align*}\begin{split}  
 d_{\rho}(x,y) &\geq \frac{1}{2} B_{0}^{-1} \rho(x)^{-1}
 3^{-\frac{k_{0}}{k_{0} + 1}} \int^{1}_{0} \abs{\gamma'(t)} \, dt \\
 &\geq \frac{1}{2} B_{0}^{-1} \rho(x)^{-1}
 3^{-\frac{k_{0}}{k_{0} + 1}} \abs{x - y}.
 \end{split}\end{align*}
 Next, suppose that the curve $\gamma$ leaves the ball $B(x,2\rho(x))$.
Due to the continuity of the curve $\gamma$, there must then exist $a \in
(0,1)$ for which $\abs{x - \gamma(a)} = 2 \rho(x)$ but $\abs{x -
  \gamma(t)} < 2 \rho(x)$ for all $t \in [0,a)$. Then
\begin{align*}\begin{split}  
 d_{\rho}(x,y) &\geq \frac{1}{2} B_{0}^{-1} \rho(x)^{-1} \int^{a}_{0}
 \br{1 + \frac{\abs{x - \gamma(t)}}{\rho(x)}}^{-\frac{k_{0}}{k_{0} +
     1}} \abs{\gamma'(t)} \, dt \\
 &\geq \frac{1}{2} B_{0}^{-1} 3^{-\frac{k_{0}}{k_{0} + 1}}
 \rho(x)^{-1} \int^{a}_{0} \abs{\gamma'(t)} \, dt \\
 &\geq \frac{1}{2} B_{0}^{-1} 3^{-\frac{k_{0}}{k_{0} + 1}}
 \rho(x)^{-1} \abs{x - y}.
 \end{split}\end{align*}
 \end{prof}
  
  The following lemma will allow
  us to compare the Agmon distance $d_{\rho}(x,y)$ with the quantity
  $\br{1 + \frac{\abs{x - y}}{\rho(x)}}$ at a global scale.

\begin{lem}
  \label{lem:Shen}
Let $\rho : \R^{d} \rightarrow [0,\infty)$ be a critical radius
function.  There exists $D_{1} > 1$, dependent on $\rho$ only through
$B_{0}$ and $k_{0}$, such that
  \begin{equation}
    \label{eqtn:lem:Shen1}
    d_{\rho}(x,y) \leq D_{1} \br{1 + \frac{\abs{x - y}}{\rho(x)}}^{k_{0} + 1}
  \end{equation}
  for all $x, \, y \in \R^{d}$ and
  \begin{equation}
    \label{eqtn:lem:Shen2}
    d_{\rho}(x,y) \geq D_{1}^{-1} \br{1 + \frac{\abs{x -
          y}}{\rho(x)}}^{\frac{1}{k_{0} + 1}}
  \end{equation}
  for all $x, \, y \in \R^{d}$ satisfying $\abs{x - y} \geq \rho(x)$.
\end{lem}

\begin{prof}  
 This was proved for the case $\rho = \rho_{V}$ in
 \cite{shen1999fundamental}. For a general critical radius function,
 the estimates follow using an identical proof.
 \end{prof}

Define the constant
\begin{equation}
  \label{eqtn:Beta}
\beta := \max \br{B_{0},D_{0},D_{1},2}. 
\end{equation}
It is obvious that \eqref{eqtn:Shen0} and Lemmas  \ref{lem:LocalAgmon} and
\ref{lem:Shen} will all hold with the constant $\beta$ replacing
$B_{0}$, $D_{0}$ and $D_{1}$ respectively. The constant $\beta$ will
thus be used as a way of simplifying notation by consolidating the three
constants $B_{0}$, $D_{0}$ and $D_{1}$ into the single constant $\beta$.

\vspace*{0.1in}

We introduce the notation $B_{\rho}(x,r)$ to denote the open ball in the
metric $d_{\rho}$ centered at the point $x \in \R^{d}$ and of radius $r
> 0$. That is,
$$
B_{\rho}(x,r) := \lb y \in \R^{d} : d_{\rho}(x,y) < r \rb
$$
 As a consequence of Lemmas \ref{lem:LocalAgmon} and \ref{lem:Shen}, we are able to compare
 balls in the Euclidean metric with balls in the metric $d_{\rho}$.

\begin{lem}
  \label{lem:Balls}
Let $\rho : \R^{d} \rightarrow [0,\infty)$ be a critical radius function.  Let $r > 0$ and $x \in \R^{d}$. Suppose that $r \leq 2$. Then
  $$
B(x,r \rho(x)) \subset B_{\rho}(x,\beta r).
$$
Suppose instead that $r > 2$. Then
$$
B(x,r \rho(x)) \subset B_{\rho}(x,\beta \br{1 + r}^{k_{0} + 1}).
$$
\end{lem}

\begin{prof}  
 First suppose that $r \leq 2$ and let $y \in B(x, r \rho(x))$. Then from
 Lemma \ref{lem:LocalAgmon},
 \begin{align*}\begin{split}  
     d_{\rho}(x,y) &\leq \beta \abs{x - y} \rho(x)^{-1} \\
     &< \beta r,
   \end{split}\end{align*}
 which implies that $y \in B_{\rho}(x,\beta r)$.

 Next, suppose that $r > 2$ and let $y \in B(x,r \rho(x))$. Lemma \ref{lem:Shen} implies
 \begin{align*}\begin{split}  
 d_{\rho}(x,y) &\leq \beta \br{1 + \frac{\abs{x - y}}{\rho(x)}}^{k_{0} +
   1} \\
 &< \beta \br{1 + r}^{k_{0} + 1}.
 \end{split}\end{align*}
 \end{prof}

 \begin{lem}
   \label{lem:Balls2}
   Let $\rho : \R^{d} \rightarrow [0,\infty)$ be a critical radius
   function. There exists a constant $A_{0} > 1$, dependent on $\rho$ only
   through $B_{0}$ and $k_{0}$, such that for all $r > 0$ with $r
\leq \beta$ and $x \in \R^{d}$,
   $$
B_{\rho}(x,r) \subset B(x,A_{0} r \rho(x)).
$$
Also, for $x \in \R^{d}$ and $r > \beta$,
$$
  B_{\rho}(x,r) \subset B\br{x, ((r \beta)^{k_{0} + 1} - 1) \rho(x)}.
  $$   

\end{lem}

\begin{prof}
  Fix $A_{0} > 1$ to be a constant large enough so that
  $$
 \frac{A_{0}}{2 \beta \br{1 + A_{0}
    \beta}^{\frac{k_{0}}{k_{0} + 1}}} \geq 1.
  $$
Suppose first that $r \leq \beta$ and fix $y \in B(x, A_{0}r \rho(x))^{c}$. Let $\gamma$ be a curve
 connecting the points $x$ and $y$ such that
 $$
d_{\rho}(x,y) \geq \frac{1}{2} \int^{1}_{0} \rho(\gamma(t))^{-1}
\abs{\gamma'(t)} \, dt.
$$
On applying \eqref{eqtn:Shen0},
$$
d_{\rho}(x,y) \geq \frac{1}{2} \beta^{-1} \rho(x)^{-1} \int^{1}_{0} \br{1
+ \frac{\abs{x - \gamma(t)}}{\rho(x)}}^{-\frac{k_{0}}{k_{0} + 1}}
\abs{\gamma'(t)} \, dt.
$$
Let $a \in (0,1)$ be the unique point in the interval that satisfies
$\abs{x - \gamma(a)} = A_{0} r \rho(x)$ and $\abs{x - \gamma(t)} <
A_{0} r
\rho(x)$ for all $t \in [0,a)$. Then
\begin{align*}\begin{split}  
 d_{\rho}(x,y) &\geq \frac{1}{2} \beta^{-1} \rho(x)^{-1} \int^{a}_{0}
 \br{1 + \frac{\abs{x - \gamma(t)}}{\rho(x)}}^{-\frac{k_{0}}{k_{0} +
     1}} \abs{\gamma ' (t)} \, dt \\
 &\geq \frac{1}{2} \beta^{-1} \br{1 + A_{0}
   \beta}^{-\frac{k_{0}}{k_{0} + 1}} \rho(x)^{-1} \int^{a}_{0}
 \abs{\gamma'(t)} \, dt \\
 &\geq \frac{1}{2} \beta^{-1} \br{1 + A_{0}
   \beta}^{-\frac{k_{0}}{k_{0} + 1}} A_{0} r  \\
 &\geq r ,
\end{split}\end{align*}
which completes the proof of the first inclusion.

\vspace*{0.1in}

For the second inclusion, suppose that $r > \beta$ and fix $y \in
B(x,((r \beta)^{k_{0} + 1} - 1) \rho(x))^{c}$. Since $r > \beta$, it
follows that 
$$
\abs{x - y} \geq \br{\br{r \beta}^{k_{0} + 1} - 1} \rho(x) \geq \br{4^{k_{0} + 1}
- 1} \rho(x) > \rho(x).
$$
This allows us to apply the second part of Lemma
\ref{lem:Shen} to obtain
\begin{align*}\begin{split}  
 d_{\rho}(x,y) &\geq  \beta^{-1} \br{1 + \frac{\abs{x -
       y}}{\rho(x)}}^{\frac{1}{k_{0} + 1}} \\
 &\geq r,
\end{split}\end{align*}
which tells us that $y \in B_{\rho}(x,r)^{c}$.
\end{prof}

\begin{cor}
  \label{cor:NonDoubling}
Let $\rho : \R^{d} \rightarrow [0,\infty)$ be a critical radius function.  For all $x \in \R^{d}$ and $r > 0$,
  $$
\abs{B_{\rho}(x,2r)} \lesssim (1 + r)^{(k_{0} + 1)d} \abs{B_{\rho}(x,r)}.
  $$
\end{cor}

\begin{prof}  
 Suppose first that $r \leq \beta / 2$. Successively applying Lemmas
 \ref{lem:Balls2} and \ref{lem:Balls},
 \begin{align*}\begin{split}  
     \abs{B_{\rho}(x,2r)} &\leq \abs{B(x,2 A_{0}r \rho(x))} \\
     &\simeq \abs{B\br{x,\frac{r}{\beta} \rho(x)}} \\
     &\leq \abs{B_{\rho}(x,r)} \\
     &\leq (1 + r)^{(k_{0} + 1)d}\abs{B_{\rho}(x,r)}.
   \end{split}\end{align*}
 Next suppose that $\beta / 2 < r \leq 2^{(k_{0} + 1)}\beta$. Then Lemmas
 \ref{lem:Balls2} and \ref{lem:Balls} lead to
 \begin{align*}\begin{split}  
     \abs{B_{\rho}(x,2r)} &\leq \abs{B(x,(2 r \beta)^{k_{0} + 1} \rho(x))} \\
     &\lesssim  r^{(k_{0} + 1)d}\abs{B\br{x, \frac{r
           \rho(x)}{2^{(k_{0} + 1)} \beta}}} \\
     &\leq r^{(k_{0} + 1)d} \abs{B_{\rho}(x,r)}.
   \end{split}\end{align*}
 Finally, suppose that $r > 2^{(k_{0} + 1)} \beta$. We would then
 have
 \begin{align*}\begin{split}  
     \abs{B_{\rho}(x,2r)} &\leq \abs{B(x,(2 r \beta)^{k_{0} + 1} \rho(x))} \\
     &\lesssim r^{(k_{0} + 1)d} \abs{B\br{x, \frac{1}{2}
         \br{\frac{r}{\beta}}^{\frac{1}{k_{0} + 1}} \rho(x)}} \\
     &\lesssim r^{(k_{0} + 1)d} \abs{B_{\rho}(x,r)},
   \end{split}\end{align*}
 where the last line follows from Lemma \ref{lem:Balls}.
\end{prof}

For an operator $S$ acting on functions in $L^{1}_{loc}\br{\R^{d}}$,
we define the local and global components through
$$
S_{loc}f(x) := S \br{f \cdot \mathbbm{1}_{B(x,\rho (x))}}(x), \qquad
S_{glob}f(x) := S \br{f \cdot \mathbbm{1}_{B(x,  \rho (x))^{c}}}(x)
$$
for $x \in \R^{d}$ and $f \in L^{1}_{loc}\br{\R^{d}}$. In order to
prove that an operator $S$ is
bounded on some weighted space $L^{p}(w)$, it is sufficient to prove
that both $S_{loc}$ and $S_{glob}$ are bounded on $L^{p}(w)$.
The below proposition will be of vital importance for proving the
$L^{p}(w)$-boundedness of our operators.

\begin{prop}[\cite{bongioanni2011classes}]
 \label{prop:Cover} 
 There exists a sequence of points $\lb x_{j} \rb_{j \in \N} \subset
 \R^{d}$ that satisfies the following two properties,
 \begin{enumerate}[(i)]
 \item $\R^{d} = \bigcup_{j \in \N} B(x_{j}, \rho(x_{j}))$,
 \item There exists $C, \, N_{1} > 0$ such that for every $\sigma \geq
   1$
   $$
\sum_{j \in \N} \mathbbm{1}_{B(x_{j}, \sigma \rho(x_{j}))} \leq C \sigma^{N_{1}}.
   $$
   \end{enumerate}
 \end{prop}

 \section{The Adapted Weight Classes}
 \label{sec:Adapted}

Throughout this section, let $\rho : \R^{d} \rightarrow
[0,\infty)$ be a critical radius function satisfying
\eqref{eqtn:Shen0} with constants $B_{0}, \, k_{0} > 1$. The Agmon
distance corresponding to this critical radius function allows one to define classes $S_{p,c}^{\rho}$
and $H_{p,c}^{\rho,m}$ in an identical manner to the classes
$S_{p,c}^{V}$ and $H_{p,c}^{V,m}$.

\begin{deff} 
  \label{def:ApVc}
   Let $1 < p < \infty$ and $c > 0$. 
 $S_{p,c}^{\rho}$ is the class of all weights $w$ on $\R^{d}$ for which
$$
\brs{w}_{S_{p,c}^{\rho}} := \sup_{B_{\rho}} \br{\frac{1}{\abs{B_{\rho}}e^{c \cdot r}}
  \int_{B_{\rho}} w}^{\frac{1}{p}} \br{\frac{1}{\abs{B_{\rho}}e^{c \cdot r}} \int_{B_{\rho}}
  w^{-\frac{1}{p-1}}}^{\frac{p - 1}{p}} < \infty,
$$
where the supremum is taken over all balls $B_{\rho} = B_{\rho}(x,r) \subset \R^{d}$
in the metric $d_{\rho}$ with $x \in \R^{d}$ and $r > 0$.
\end{deff}

\begin{deff} 
  \label{def:HeatWeights}
 Let $1 < p < \infty$ and $m, \, c > 0$. Let
 $\Phi_{m,c}^{\rho} : \R^{d} \xx (0,\infty) \rightarrow \R$ be the function defined by
$$
\Phi_{m,c}^{\rho}(x,r) := \exp \br{c \br{1 + \frac{r}{\rho(x)}}^{m}},
$$
for $x \in \R^{d}$ and $r > 0$.
 Let $H_{p,c}^{\rho,m}$ denote the class
 of all weights $w$ on $\R^{d}$ for which
 $$
\brs{w}_{H_{p,c}^{\rho,m}} := \sup_{B} \br{\frac{1}{\abs{B} \Phi_{m,c}^{\rho}(x,r)}
\int_{B} w}^{\frac{1}{p}} \br{\frac{1}{\abs{B} \Phi_{m,c}^{\rho}(x,r)} \int_{B}
w^{-\frac{1}{p-1}}}^{\frac{p - 1}{p}} < \infty,
$$
where the supremum is taken over all Euclidean balls $B = B(x,r) \subset \R^{d}$
with radius $r > 0$ and center $x \in \R^{d}$.
\end{deff}

We clearly have $S_{p,c}^{\rho_{V}} = S_{p,c}^{V}$ and
$H_{p,c}^{\rho_{V},m} = H_{p,c}^{V,m}$. Let $A_{p}^{\rho,loc}$ be as
defined in \cite{bongioanni2011classes}. That is, $A_{p}^{\rho,loc}$
is the collection of all weights $w$ for which
$$
w(B)^{\frac{1}{p}} w^{-\frac{1}{p - 1}} (B)^{\frac{p - 1}{p}} \lesssim \abs{B}
$$
for all balls $B = B(x,r)$ with $r \leq \rho(x)$.

 \begin{prop}
   \label{prop:LocalApVc}
   For any $1 < p < \infty$ and $c > 0$,
   $$
S_{p,c}^{\rho} \subset A_{p}^{\rho,loc}.
   $$
 \end{prop}

 \begin{prof}  
 Fix $w \in S_{p,c}^{\rho}$ for some $1 < p < \infty$ and $c > 0$. Let
 $B := B(x, r \rho(x))$ be a ball in the Euclidean metric with $r \leq
 1$. On successively applying Lemma \ref{lem:Balls}, the condition
 that $w \in S_{p,c}^{\rho}$ and finally the bound $r \leq 1$,
 \begin{align*}\begin{split}  
 w \br{B}^{\frac{1}{p}} w^{-\frac{1}{p - 1}} \br{B}^{\frac{p - 1}{p}}
 &\leq w \br{B_{\rho}(x,\beta r)}^{\frac{1}{p}} w^{-\frac{1}{p - 1}}
 \br{B_{\rho}(x,\beta r)}^{\frac{p - 1}{p}} \\
 &\lesssim e^{c \beta r} \abs{B_{\rho}(x,\beta r)} \\
 &\lesssim \abs{B_{\rho}(x,\beta r)}.
\end{split}\end{align*}
Since $\beta r \leq \beta$, Lemma \ref{lem:Balls2} can be
applied to give
\begin{align*}\begin{split}  
 w \br{B}^{\frac{1}{p}} w^{-\frac{1}{p - 1}} \br{B}^{\frac{p - 1}{p}}
 &\lesssim \abs{B(x, A_{0} \beta r \rho(x))} \\
 &\simeq \abs{B(x,r \rho(x))}.
\end{split}\end{align*}
This proves that $w \in A_{p}^{\rho,loc}$.
\end{prof}

Recall the definition of the class $A_{p}^{\rho,\infty}$ introduced in \cite{bongioanni2011classes}.

\begin{deff}[\cite{bongioanni2011classes}]
  \label{def:ApRhoInf}
For $1 < p < \infty$ and $\theta \geq 0$, a weight $w$ on $\R^{d}$ is said to
belong to the class $A_{p}^{\rho,\theta}$ if there exists a constant
$C > 0$ for which
$$
w\br{B}^{\frac{1}{p}} w^{-\frac{1}{p - 1}} \br{B}^{\frac{p - 1}{p}} \leq C \abs{B} \br{1 + \frac{r}{\rho(x)}}^{\theta}
$$
for all balls $B = B(x,r)$ with center $x \in \R^{d}$ and radius $r >
0$. Define
$$
A_{p}^{\rho,\infty} := \bigcup_{\theta \geq 0} A_{p}^{\rho,\theta}.
$$
\end{deff}

For $V \in RH_{\frac{d}{2}}$, the class $A_{p}^{V,\infty}$ discussed
in the introductory section is then defined by $A_{p}^{V,\infty} = A_{p}^{\rho_{V},\infty}$.

\begin{prop} 
 \label{prop:Inclusion} 
Let $1 < p < \infty$. For any $c_{1}, \, c_{2}, \, c_{3} > 0$, $m_{1} < 
\br{k_{0} + 1}^{-1}$ and $m_{2} > (k_{0} + 1)$ we have
 $$
A_{p}^{\rho,\infty} \subset H_{p,c_{1}}^{\rho,m_{1}} \subset
S_{p,c_{2}}^{\rho} \subset H_{p,c_{3}}^{\rho,m_{2}}.
 $$
\end{prop}

\begin{prof}  
 The first inclusion is almost trivial. Fix $w \in
 A_{p}^{\rho,\infty}$. Then there must exist some $\theta \geq 0$ for
 which $w \in A_{p}^{\rho,\theta}$. Since $x^{\theta} \lesssim e^{c_{1} x^{m_{1}}}$ we have
 \begin{align*}\begin{split}   
 w \br{B}^{\frac{1}{p}} w^{-\frac{1}{p - 1}} \br{B}^{\frac{p - 1}{p}}
 &\lesssim \br{1 + \frac{r}{\rho(x)}}^{\theta} \abs{B} \\
 &\lesssim \exp \br{c_{1} \br{1 + \frac{r}{\rho(x)}}^{m_{1}}} \abs{B} \\
 &= \Phi_{m_{1}, c_{1}}(x,r) \abs{B}
\end{split}\end{align*}
for all balls $B := B(x,r) \subset \R^{d}$.

\vspace*{0.1in}

Let's prove the second inclusion. Fix $w \in H_{p,c_{1}}^{\rho,m_{1}}$ for
some $c_{1} > 0$ and $m_{1} < (k_{0} + 1)^{-1}$. Also fix $c_{2} > 0$.  Let $B_{\rho} := B_{\rho}(x,r) \subset \R^{d}$ be a ball in the metric
$d_{\rho}$ for some $x \in \R^{d}$ and $r > 0$. First suppose that
$r > \beta$. Then
$$
w\br{B_{\rho}}^{\frac{1}{p}} w^{-\frac{1}{p - 1}} \br{B_{\rho}}^{\frac{p -
    1}{p}} \leq w \br{B(x,r' \rho(x))}^{\frac{1}{p}} w^{-\frac{1}{p - 1}}
\br{B(x, r' \rho(x))}^{\frac{p - 1}{p}}
$$
by Lemma \ref{lem:Balls2}, where $r' := \br{\br{r \beta}^{k_{0} + 1} - 1}$. Applying
the condition $w \in H_{p,c_{1}}^{\rho,m_{1}}$ gives
\begin{align*}\begin{split}  
 w \br{B_{\rho}}^{\frac{1}{p}} w^{-\frac{1}{p - 1}} \br{B_{\rho}}^{\frac{p -
   1}{p}} &\lesssim \abs{B(x,r' \rho(x))} \exp \br{c_{1} \br{1 +
   \frac{r' \rho(x)}{\rho(x)}}^{m_{1}}} \\
&= \abs{B(x,r'\rho(x))} \exp \br{c_{1} \br{r \beta}^{m_{1} (k_{0} + 1)}} \\
&\lesssim \rho(x)^{d} \br{r \beta}^{(k_{0} + 1)d} \exp
\br{c_{1} \br{r \beta}^{m_{1} (k_{0} + 1)}} \\
&\lesssim \rho(x)^{d} \exp \br{c' r^{m_{1}(k_{0} + 1)}},
\end{split}\end{align*}
for any $c' > c_{1}$. Since $m_{1} < \br{k_{0} + 1}^{-1}$,
$$
 w \br{B_{\rho}}^{\frac{1}{p}} w^{-\frac{1}{p - 1}} \br{B_{\rho}}^{\frac{p -
   1}{p}} \lesssim \rho(x)^{d} e^{c_{2} r},
$$
for any $c_{2} > 0$. 
Lemma \ref{lem:Balls} tells us that we must have the inclusion
$B(x,\rho(x)) \subset B_{\rho}(x,r)$ and therefore
$$
\abs{B_{\rho}(x,r)} \geq \abs{B(x,\rho(x))} = \rho(x)^{d}.
$$
This gives
$$
w(B_{\rho})^{\frac{1}{p}} w^{-\frac{1}{p - 1}}\br{B_{\rho}}^{\frac{p -
    1}{p}} \lesssim \abs{B_{\rho}} e^{c_{2} r}.
$$
We must now consider the case $r \leq \beta$. Lemma
\ref{lem:Balls2} implies that
\begin{align*}\begin{split}  
 w \br{B_{\rho}}^{\frac{1}{p}} w^{-\frac{1}{p - 1}} \br{B_{\rho}}^{\frac{p -
   1}{p}} &\leq w \br{B(x,A_{0} r \rho(x))}^{\frac{1}{p}}
w^{-\frac{1}{p - 1}} \br{B(x,A_{0} r \rho(x))}^{\frac{p - 1}{p}} \\
&\lesssim \abs{B(x,A_{0} r \rho(x))} \exp \br{c_{1} \br{1 + A_{0}
    r}^{m_{1}}}.
\end{split}\end{align*}
Since $r \leq \beta$,
\begin{align*}\begin{split}  
 w \br{B_{\rho}}^{\frac{1}{p}} w^{-\frac{1}{p - 1}} \br{B_{\rho}}^{\frac{p -
     1}{p}} &\lesssim \abs{B(x, A_{0} r \rho(x))} \\
 &\leq \abs{B(x, A_{0} r \rho(x))} e^{c_{2} r} \\
 &\simeq \abs{B(x, \beta^{-1}  r \rho(x))} e^{c_{2} r}.
\end{split}\end{align*}
Lemma \ref{lem:Balls} then gives
$$
w \br{B_{\rho}}^{\frac{1}{p}} w^{-\frac{1}{p - 1}} \br{B_{\rho}}^{\frac{p -
     1}{p}} \lesssim \abs{B_{\rho}(x, r)} e^{c_{2} r}.
$$
This completes the proof of $w \in S_{p,c_{2}}^{\rho}$.

\vspace*{0.1in}

Finally, let's prove the last inclusion. Let $w \in
S_{p,c_{2}}^{\rho}$ with $c_{2} > 0$. Fix $c_{3} > 0$ and $m_{2} >
(k_{0} + 1)$. Let $r > 0$ and $x \in
\R^{d}$. First consider the case $r \leq 1$. On consecutively applying Lemma
\ref{lem:Balls}, the hypothesis $w \in S_{p,c_{2}}^{\rho}$ and Lemma
\ref{lem:Balls2},
\begin{align*}\begin{split}  
 w \br{B(x,r\rho(x))}^{\frac{1}{p}} w^{-\frac{1}{p - 1}} \br{B(x,r
   \rho(x))}^{\frac{p - 1}{p}} &\leq w \br{B_{\rho}(x,\beta
   r)}^{\frac{1}{p}} w^{-\frac{1}{p - 1}} \br{B_{\rho}(x,\beta
   r)}^{\frac{p - 1}{p}} \\
&\lesssim \abs{B_{\rho}(x,\beta r)} e^{c_{2} \beta r} \\
&\lesssim \abs{B(x, A_{0} \beta r \rho(x))} e^{c_{2} \beta
r} \\
&\lesssim r^{d} \rho(x)^{d} e^{c_{3} (1 + r)^{m_{2}}}.
 \end{split}\end{align*}
Next consider $r \geq 1$. On successively applying Lemma
\ref{lem:Balls}, the hypothesis $w \in S_{p,c_{2}}^{\rho}$ and Lemma \ref{lem:Balls2},
\begin{align*}\begin{split}  
& w(B(x,r \rho(x)))^{\frac{1}{p}} w^{-\frac{1}{p - 1}} (B(x,r
 \rho(x)))^{\frac{p - 1}{p}}  \\ & \qquad \qquad \leq w (B_{\rho}(x,\beta(1 + r)^{k_{0} +
   1}))^{\frac{1}{p}} w^{-\frac{1}{p - 1}} (B_{\rho}(x,\beta(1 +
 r)^{k_{0} + 1}))^{\frac{p - 1}{p}} \\
& \qquad \qquad \lesssim \abs{B_{\rho}(x,\beta(1 + r)^{k_{0} + 1})} e^{c_{2} \beta (1 +
  r)^{k_{0} + 1}} \\
& \qquad \qquad \lesssim \abs{B(x, \beta^{2(k_{0} + 1)} (1 + r)^{(1 + k_{0})^{2}} \rho(x))}
e^{c_{2} \beta (1 + r)^{k_{0} + 1}} \\
& \qquad \qquad \lesssim  r^{d} \rho(x)^{d} e^{c_{3}(1 + r)^{m_{2}}}.
 \end{split}\end{align*}
\end{prof}

\begin{deff} 
 \label{def:AdaptedAveraging} 
Let $c > 0$. For each $t > 0$, define the averaging operator
 $$
A_{t,c}^{\rho}f(x) := \frac{1}{\abs{B_{\rho}(x,t)} e^{c t}} \int_{B_{\rho}(x,t)} f(y) \, dy
$$
for $f \in L^{1}_{loc}(\R^{d})$ and $x \in \R^{d}$. Define the centered Hardy-Littlewood operator associated with
$\rho$ through
$$
M_{\rho,c}f(x) := \sup_{t > 0} A_{t,c}^{\rho}\abs{f}(x).
$$
Similarly, define the uncentered Hardy-Littlewood operator by
$$
\mathcal{M}_{\rho,c}f(x) := \sup_{B_{\rho} \ni x} \frac{1}{e^{c r} \abs{B_{\rho}}}
\int_{B_{\rho}} \abs{f(y)} \, dy,
$$
where the supremum is taken over all balls $B_{\rho} = B_{\rho}(x',r) \subset \R^{d}$
in the metric $d_{\rho}$ that contain $x$.
\end{deff}

For any $x \in \R^{d}$ and $f \in L^{1}_{loc}(\R^{d})$, the inequality
$M_{\rho,c}f(x) \leq \mathcal{M}_{\rho,c}f(x)$ is trivial. The below
proposition states that a weak converse will hold.

\begin{prop} 
 \label{prop:ConverseCenter} 
 For any $c_{1}, \, c_{2} > 0$ with $c_{1} > 2 c_{2}$ we have
 $$
\mathcal{M}_{\rho,c_{1}}f(x) \lesssim M_{\rho,c_{2}} f(x)
$$
for all $f \in L^{1}_{loc}(\R^{d})$ and $x \in \R^{d}$.
\end{prop}

\begin{prof}  
Fix $c_{1}, \,  c_{2} > 0$ with $c_{1} > 2 c_{2}$  and let $x \in \R^{d}$ and $f \in
L^{1}_{loc}(\R^{d})$. Let $B_{\rho} = B_{\rho}(x',r) \subset \R^{d}$
for $x' \in \R^{d}$ and $r > 0$ be a
ball in the metric $d_{\rho}$ that contains the point $x$. We have by
Corollary \ref{cor:NonDoubling},
\begin{align*}\begin{split}  
 \abs{B_{\rho}(x,2 r)} &\lesssim e^{\br{\frac{c_{1}}{2} - c_{2}} r}
 \abs{B_{\rho}(x,r)} \\
 &\leq e^{\br{\frac{c_{1}}{2} - c_{2}}r} \abs{B_{\rho}(x',2r)} \\
 &\lesssim e^{(c_{1} - 2 c_{2}) r} \abs{B_{\rho}(x',r)}.
\end{split}\end{align*}
Therefore,
\begin{align*}\begin{split}  
 \frac{1}{e^{c_{1} r} \abs{B_{\rho}(x',r)}} \int_{B_{\rho}(x',r)} \abs{f(y)}
 \, dy &\lesssim \frac{e^{(c_{1} - 2 c_{2})r}}{ e^{c_{1}r}
   \abs{B_{\rho}(x,2r)}} \int_{B_{\rho}(x',r)} \abs{f(y)} \, dy \\
 &\leq \frac{1}{e^{2 c_{2}r} \abs{B_{\rho}(x,2r)}} \int_{B_{\rho}(x,2r)}
 \abs{f(y)} \, dy \\
 &\leq M_{c_{2}}f(x).
\end{split}\end{align*}
Taking the supremum over all $B_{\rho} = B_{\rho}(x',r)$ then proves the proposition.
 \end{prof}

\begin{prop} 
 \label{prop:HardyLittlewood} 
 Fix $1 < p < \infty$. For any $c_{1}, \, c_{2} > 0$ with $c_{1} > 2 c_{2}$,
 $$
\lb w : \norm{M_{\rho,c_{2}}}_{L^{p}(w)} < \infty \rb \subset S_{p,c_{1}}^{\rho}.
 $$
\end{prop}

\begin{prof}
  Fix $1 < p < \infty$ and $c_{1}, \, c_{2} > 0$ with $c_{1} > 2 c_{2}$.
 Let $w$ be a weight for which $\norm{M_{\rho,c_{2}}}_{L^{p}(w)} <
 \infty$. The proof that $w \in S_{p,c_{1}}^{V}$ is similar to the standard
 classical proof that can be found in \cite{grafakos2009modern} for example. Fix
 $B_{\rho} := B_{\rho}(x,r)$ for some $x \in \R^{d}$ and $r > 0$. The
 boundedness of $M_{\rho,c_{2}}$ on $L^{p}(w)$ together with Proposition \ref{prop:ConverseCenter} implies that
 \begin{align*}\begin{split}  
 w(B_{\rho}) \br{\frac{1}{e^{c_{1} r}\abs{B_{\rho}}} \int_{B_{\rho}} \abs{f}}^{p}
 &\lesssim \int_{B_{\rho}} \mathcal{M}_{\rho,c_{1}} (f \mathbbm{1}_{B_{\rho}})(y)^{p} w(y) \,
 dy \\
 &\lesssim \int_{B_{\rho}} M_{\rho,c_{2}}(f \mathbbm{1}_{B_{\rho}})(y)^{p}
 w(y) \, dy \\
 &\lesssim \int_{B_{\rho}} \abs{f}^{p} w(y) \, dy.
\end{split}\end{align*}
For $\epsilon > 0$, take $f := (w + \epsilon)^{-\frac{1}{p - 1}}$ in
the above inequality to obtain
$$
\frac{1}{e^{c_{1} p r} \abs{B_{\rho}}^{p}} w(B_{\rho}) \br{\int_{B_{\rho}} (w +
  \epsilon)^{-\frac{1}{p - 1}}}^{p} \lesssim \br{\int_{B_{\rho}} (w +
  \epsilon)^{-\frac{p}{p - 1}}w } \leq \br{\int_{B_{\rho}} (w +
  \epsilon)^{-\frac{1}{p - 1}}}.
$$
Which leads to
$$
\br{\frac{1}{e^{c_{1} r}\abs{B_{\rho}}} \int_{B_{\rho}}w} \br{\frac{1}{e^{c_{1} r}
    \abs{B_{\rho}}} \int_{B_{\rho}} (w + \epsilon)^{-\frac{1}{p - 1}}}^{p -
  1} \leq C
$$
for some constant $C > 0$. The monotone convergence theorem then
allows us to conclude that $w \in S_{p,c_{1}}^{\rho}$.
 \end{prof}

 \begin{deff} 
 \label{def:AdaptedAveragingH} 
Let $c, \, m > 0$. For each $t > 0$, define the averaging operator
 $$
\tilde{A}_{t,c}^{\rho,m}f(x) := \frac{1}{ \Phi_{m,c}^{\rho}(x,t) \abs{B(x,t)}} \int_{B(x,t)} f(y) \, dy
$$
for $f \in L^{1}_{loc}(\R^{d})$ and $x \in \R^{d}$. Define the corresponding centered Hardy-Littlewood operator through
$$
\tilde{M}_{\rho,c}^{m}f(x) := \sup_{t > 0} \tilde{A}_{t,c}^{\rho,m}\abs{f}(x).
$$
Similarly, define the uncentered Hardy-Littlewood operator by
$$
\tilde{\mathcal{M}}_{\rho,c}^{m}f(x) := \sup_{B \ni x} \frac{1}{\Phi_{m,c}^{\rho}(x',r) \abs{B}}
\int_{B} \abs{f(y)} \, dy,
$$
where the supremum is taken over all Euclidean balls $B = B(x',r)
\subset \R^{d}$ that contain the point $x$.
\end{deff}

For these operators, an analogue of the pointwise bound from Proposition
\ref{prop:ConverseCenter} will hold. 

\begin{prop}
  \label{prop:TempPointwise}
  For any $c_{1}, \, c_{2}, \, m_{1}, \, m_{2} > 0$ with $m_{1} >
  (k_{0} + 1) m_{2}$ we have
  $$
\tilde{\mathcal{M}}_{\rho,c_{1}}^{m_{1}} f(x) \lesssim \tilde{M}_{\rho,c_{2}}^{m_{2}}f(x)
$$
for all $f \in L_{loc}^{1}\br{\R^{d}}$ and $x \in \R^{d}$.
\end{prop}

\begin{prof}  
 Fix $c_{1}, \, c_{2}, \, m_{1}, \, m_{2} > 0$ with $m_{1} > (k_{0} +
 1) m_{2}$ and
 let $x \in \R^{d}$ and $f \in L^{1}_{loc}\br{\R^{d}}$. Let $B(x',r) \subset \R^{d}$ be a ball in the Euclidean metric that
 contains the point $x$. We have by Lemma \ref{lem:Shen0}
 \begin{align*}\begin{split}  
 1 + \frac{2 r}{\rho(x)} &\leq 1 + 2 B_{0} \frac{r}{\rho(x')}
 \br{1 + \frac{r}{\rho(x')}}^{k_{0}} \\
 &\leq 2 B_{0} \br{1 + \frac{r}{\rho(x')}}^{k_{0} + 1}.
\end{split}\end{align*}
Therefore,
\begin{align*}\begin{split}  
 \Phi^{\rho}_{m_{2},c_{2}}(x,2 r) &= \exp \br{c_{2} \br{1 +
     \frac{2 r}{\rho(x)}}^{m_{2}}} \\
 &\leq \exp \br{2 c_{2} B_{0} \br{1 + \frac{r}{\rho(x')}}^{(k_{0} +
     1)m_{2}}} \\
 &\lesssim \exp \br{c_{1} \br{1 + \frac{r}{\rho(x')}}^{m_{1}}} = \Phi^{\rho}_{m_{1},c_{1}}(x',r).
\end{split}\end{align*}
Which gives
\begin{align*}\begin{split}  
 \frac{1}{\Phi^{\rho}_{m_{1},c_{1}}(x',r) \abs{B(x',r)}}
 \int_{B(x',r)} \abs{f(y)} \, dy &\lesssim
 \frac{1}{\Phi^{\rho}_{m_{2},c_{2}}(x,2 r) \abs{B(x,2r)}} \int_{B(x,2r)}
 \abs{f(y)}  \, dy \\
 &\leq \tilde{M}_{\rho,c_{2}}^{m_{2}}f(x).
\end{split}\end{align*}
Taking the supremum over all $B(x',r)$ that contains $x$ then proves
the proposition.

 \end{prof}

The pointwise estimate from the previous proposition then allows us to
deduce the following inclusion. The proof is identical to that of
Proposition \ref{prop:HardyLittlewood}.

\begin{prop} 
 \label{prop:HardyLittlewood2} 
Fix $1 < p < \infty$. For any $c_{1}, \, c_{2}, \, m_{1}, \, m_{2} >
0$ with $m_{1} > (k_{0} + 1) m_{2}$,
 $$
\lb w : \norm{\tilde{M}_{\rho,c_{2}}^{m_{2}}}_{L^{p}(w)} < \infty \rb \subset
H_{p,c_{1}}^{\rho, m_{1}}.
 $$
\end{prop}

\section{Schr\"{o}dinger Operators}
\label{sec:Schrodinger}

In this section, a proof of Theorems \ref{thm:Riesz} and
\ref{thm:Heat} will be provided. Theorem \ref{thm:Riesz} will be
proved by demonstrating that the statements actually hold for a more
general form of operator, the Schr\"{o}dinger operator with measure
potential $-\Delta + \mu$. We now provide a brief
description of these generalised Schr\"{o}dinger operators.

\subsection{Schr\"{o}dinger Operators with Measure Potential}
\label{subsec:SchrodingerMeasure}

Schr\"{o}dinger operators with measure potential, or generalised
Schr\"{o}dinger operators, were considered by Z. Shen in the article
\cite{shen1999fundamental}. For this form of Schr\"{o}dinger operator,
the scalar potential $V$ is replaced by a non-negative Radon measure
$\mu$ on $\R^{d}$. It is assumed that the measure $\mu$
 satisfies the property that there exists $\delta_{\mu}, \, C_{\mu},
 \, D_{\mu} > 0$ such that
\begin{equation}
  \label{eqtn:Measure1}
  \mu(B(x,r)) \leq C_{\mu} \br{\frac{r}{R}}^{d - 2 + \delta_{\mu}} \mu(B(x,R))
\end{equation}
and
\begin{equation}
  \label{eqtn:Measure2}
  \mu(B(x,2r)) \leq D_{\mu} \br{\mu(B(x,r)) + r^{d - 2}}
  \end{equation}
  for all $x \in \R^{d}$ and $0 < r < R$.

 \begin{rmk}
   \label{rmk:Measure}
For $V \in RH_{q}$ with $q > \frac{d}{2}$, the measure $d \mu(x) = V(x) \, dx$ will
satisfy both properties \eqref{eqtn:Measure1} and
\eqref{eqtn:Measure2} with constants $C_{\mu}, \, D_{\mu} > 0$, dependent on $V$ only through $\brs{V}_{RH_{q}}$, and
$\delta_{\mu} = 2 - \frac{d}{q}$ (c.f. {\cite[Lem.~1.2]{shen1995lp}}).  This
tells us that standard Schr\"{o}dinger
operators are instances of this generalised form. If $V \in
RH_{d}$ then, due to the self-improvement property of the reverse
H\"{o}lder classes, $V \in RH_{q'}$ for some $q' > d$ and therefore
\eqref{eqtn:Measure1} will be satisfied with $\delta_{\mu} = 2 -
\frac{d}{q'} > 1$.
\end{rmk}

  A proof of the following
  lemma can be found in {\cite[Lem.~2.6]{bailey2019unbounded}} for the case $d
  \mu(x) = V(x) \, dx$ with $V \in RH_{\frac{d}{2}}$. The proof for general $\mu$ is identical.

 \begin{lem}
  \label{lem:RHn2}
Let $\mu$ be a non-negative Radon measure on $\R^{d}$ that satisfies
\eqref{eqtn:Measure1} and \eqref{eqtn:Measure2} with constants
$C_{\mu}, \, D_{\mu}, \, \delta_{\mu} > 0$. For all $x \in \R^{d}$ and $R > 0$,
   \begin{equation}
    \label{eqtn:0.4}
\int_{B(x,R)} \frac{d \mu(y)}{\abs{y - x}^{d - 2}} \lesssim
\frac{\mu(B(x,R))}{R^{d - 2}}.
\end{equation}
If $\delta_{\mu} > 1$, then we will also have
  \begin{equation}
    \label{eqtn:RHn2}
\int_{B(x,R)} \frac{d \mu(y)}{\abs{y - x}^{d - 1}} \lesssim
\frac{\mu(B(x,R))}{R^{d - 1}}
\end{equation}
for all $x \in \R^{d}$ and $R > 0$.
\end{lem}

  We consider the operator
  $$
L_{\mu} := - \Delta + \mu.
$$
This can be defined rigorously through its corresponding sesquilinear
form as a non-negative unbounded operator on $L^{2}\br{\R^{d}}$ with maximal domain.
Define the operators
$$
R_{\mu} := \nabla L_{\mu}^{-\frac{1}{2}}, \qquad R_{\mu}^{*} :=
L_{\mu}^{-\frac{1}{2}} \nabla, \qquad I_{\mu}^{\alpha} := L_{\mu}^{-\frac{\alpha}{2}}
$$
for $0 < \alpha \leq 2$ and
$$
T^{*}_{\mu}f(x) := \sup_{t > 0} e^{- t L_{\mu}}\abs{f}(x)
$$
for $f \in L_{loc}^{1}\br{\R^{d}}$ and $x \in \R^{d}$.
For this measure form of the electric potential, the corresponding
critical radius function is defined through
$$
\rho_{\mu}(x) := \sup \lb r > 0 : \frac{\mu(B(x,r))}{r^{d - 2}} \leq 1
\rb, \qquad x \in \R^{d}.
$$
\begin{rmk}
  \label{rmk:Critical}
It follows directly from the definition of $\rho_{\mu}$ that for any $x \in \R^{d}$
$$
\frac{\mu(B(x,\rho_{\mu}(x)))}{\rho_{\mu}(x)^{d - 2}} \simeq 1.
$$
\end{rmk}

In \cite{shen1999fundamental} it was proved that $\rho_{\mu}$ is
indeed a critical radius function in the sense of Definition
\ref{def:CriticalRadius}. 

\begin{lem}[{\cite[Prop.~1.8, Rmk.~1.9]{shen1999fundamental}}]
  \label{lem:CriticalMeasure}
The function $\rho_{\mu}$ is a critical radius function in the sense
of Definition \ref{def:CriticalRadius}. In particular,
\eqref{eqtn:Shen0} is satisfied with constants $B_{0}$ and $k_{0}$
depending on $\mu$ only through $C_{\mu}, \, D_{\mu}$ and $\delta_{\mu}$.
\end{lem}

Through the function $\rho_{\mu}$, we can define a corresponding Agmon
distance $d_{\mu} := d_{\rho_{\mu}}$  and balls $B_{\mu}(x,r) := B_{\rho_{\mu}}(x,r)$ for $x \in \R^{d}$ and
$r > 0$. This then allows us to define
appropriate analogues of our weight classes $S_{p,c}^{\mu} := S_{p,c}^{\rho_{\mu}}$ and
$H_{p,c}^{\mu,m} := H_{p,c}^{\rho_{\mu},m}$. For these classes of weights, the following theorem
will be proved.

\vspace*{0.1in}

 \begin{thm} 
 \label{thm:MeasurePotential} 
 Let $\mu$ be a non-negative Radon measure on $\R^{d}$
 that satisfies \eqref{eqtn:Measure1} and
 \eqref{eqtn:Measure2} with constants $C_{\mu}, \, D_{\mu}, \,  \delta_{\mu}
 > 0$.
 
 \begin{enumerate}[(i)]
\item \label{Measure1} Suppose that $\delta_{\mu} > 1$. There exists $c_{1} > 0$ for which both $R_{\mu}$ and $R_{\mu}^{*}$
  are bounded on $L^{p}(w)$ for all $w \in S_{p,c_{1}}^{\mu}$ and $1 <
  p < \infty$.

  \vspace*{0.1in}

  \item \label{Measure2} Suppose instead that $0 < \delta_{\mu} < 1$
    and let $\eta \in (2,(2 - \delta_{\mu})/(1 - \delta_{\mu}))$. There exists
    $c_{2} > 0$ for which the operator $R^{*}_{\mu}$ is bounded
    on $L^{p}(w)$ for $\eta' < p < \infty$ when $w \in
    S_{p/\eta',c_{2}}^{\mu}$ and the operator $R_{\mu}$ is bounded on
    $L^{p}(w)$ for $1 < p < \eta$ when $w^{-\frac{1}{p - 1}} \in S_{p'/\eta',c_{2}}^{\mu}$.
    
    \item \label{Measure3} If $\delta_{\mu} > 0$ and  $0 < \alpha \leq 2$, there
      must exist a constant $c_{3} > 0$ for which the operator
      $I^{\alpha}_{\mu}$ is bounded from $L^{p}(w)$ to
$L^{\nu}(w^{\nu/p})$ for $w^{\nu/p} \in S_{1 +
  \frac{\nu}{p'},c_{3}}^{\mu}$ and $1 < p < \frac{d}{\alpha}$,
where $\frac{1}{\nu} = \frac{1}{p} - \frac{\alpha}{d}$.
\end{enumerate}
The constants $c_{1}$, $c_{2}$ and $c_{3}$ are independent of $p$ and 
depend on $\mu$ only through $C_{\mu}, \, D_{\mu}$ and $\delta_{\mu}$.
 \end{thm}

Let $\Gamma_{\mu}$ denote the fundamental solution of the operator
$L_{\mu}$. Refer to \cite{shen1999fundamental} for further information
and properties.
The proof of Theorem \ref{thm:MeasurePotential} will rely heavily on
the following exponential decay estimates that were proved by
Shen in \cite{shen1999fundamental}.

\begin{thm}[{\cite[Thm.~0.8, Thm.~0.17]{shen1999fundamental}}]
 \label{thm:MeasureExp} 
 Let $\mu$ be a non-negative Radon measure on $\R^{d}$ that satisfies
 \eqref{eqtn:Measure1} and \eqref{eqtn:Measure2} with constants $C_{\mu}, \,
D_{\mu}, \, \delta_{\mu} > 0$. There exist
 constants $C_{1}, \, C_{2}, \, \varepsilon_{1}, \,
 \varepsilon_{2} > 0$ for which
 $$
C_{1} \frac{e^{- \varepsilon_{1} d_{\mu}(x,y)}}{\abs{x - y}^{d - 2}} \leq
\Gamma_{\mu}(x,y) \leq C_{2} \frac{e^{- \varepsilon_{2} d_{\mu}(x,y)}}{\abs{x - y}^{d - 2}}
$$
for all $x, \, y \in \R^{d}$. There will also exist $C_{3}, \, \varepsilon_{3} > 0$
for which
\begin{equation}
  \label{eqtn:DerivativeEstimate}
\abs{\nabla \Gamma_{\mu}(x,y)} \leq C_{3} \frac{e^{-\varepsilon_{3}
    d_{\mu}(x,y)}}{\abs{x - y}^{d - 2}} \br{\int_{B(x,\abs{x - y} / 2)} \frac{d
    \mu(z)}{\abs{z - x}^{d - 1}} + \frac{1}{\abs{x - y}}}
\end{equation}
for all $x, \, y \in \R^{d}$, where the gradient is taken with respect to the first variable.
 If $\delta_{\mu} > 1$ then $C_{4}, \,
\varepsilon_{4} > 0$ can be chosen so that
$$
\abs{\nabla \Gamma_{\mu}(x,y)} \leq C_{4}  \frac{e^{- \varepsilon_{4}
    d_{\mu}(x,y)}}{\abs{x - y}^{d - 1}} \qquad \forall \ x, \, y \in \R^{d}.
$$
The constants $\varepsilon_{i}, \, C_{i}$ for $i = 1, \cdots, 4$ will
depend on $\mu$ only through $C_{\mu}, \, D_{\mu}$ and $\delta_{\mu}$.
\end{thm}

It should be noted that although the estimate
\eqref{eqtn:DerivativeEstimate} is not explicitly stated in
\cite{shen1999fundamental}, its proof is essentially contained within
the proof of {\cite[Thm.~0.17]{shen1999fundamental}}.

We will not attempt to prove the measure potential version of Theorem
\ref{thm:Heat} since, to the best of the author's knowledge, heat
kernel estimates for the general operator $L_{\mu}$ have not yet been
proved. The best known result for heat kernel estimates can be found
in \cite{kurata2000estimate} and requires 
the assumption that $d \mu(x) = V(x) \, dx$ with $V \in
RH_{\frac{d}{2}}$. Therefore, the boundedness of $T^{*}_{\mu}$ will
only be considered for this case.

 \subsection{The Riesz Transforms}
 \label{subsec:Riesz}

Let's consider the boundedness of the operators $R_{\mu}$ and
$R_{\mu}^{*}$ on the weighted Lebesgue space
$L^{p}(w)$. 
 The Riesz transforms $R_{\mu}$ can be expressed as
  \begin{align*}\begin{split}  
      R_{\mu}f(x) &= \nabla \br{-\Delta + \mu}^{-\frac{1}{2}}f(x) \\
      &= \frac{1}{\pi} \int^{\infty}_{0} \lambda^{-\frac{1}{2}} \nabla
      \br{- \Delta + \mu + \lambda}^{-1}f(x) \, d \lambda \\
      &= \frac{1}{\pi} \int^{\infty}_{0} \lambda^{-\frac{1}{2}}
      \int_{\R^{d}} \nabla \Gamma_{\mu + \lambda}(x,y) f(y) \, dy \, d \lambda,
    \end{split}\end{align*}
  where $\Gamma_{\mu + \lambda}$ is the fundamental solution of $\br{-
    \Delta + \mu + \lambda}$ (c.f. {\cite[pg.~282]{kato1980perturbation}}). Fubini's Theorem then gives
  \begin{align*}\begin{split}  
 R_{\mu}f(x) &= \int_{\R^{d}} \frac{1}{\pi} \int^{\infty}_{0}
 \lambda^{-\frac{1}{2}} \nabla \Gamma_{\mu + \lambda}(x,y) \, d\lambda \,
 f(y) \, dy \\
 &= \int_{\R^{d}} K_{\mu}(x,y) f(y) \, dy,
\end{split}\end{align*}
where $K_{\mu}$ is the singular kernel of $R_{\mu}$ given by
\begin{equation}
  \label{eqtn:KernelExpression}
K_{\mu}(x,y) = \frac{1}{\pi} \int^{\infty}_{0} \lambda^{-\frac{1}{2}}
\nabla \Gamma_{\mu + \lambda}(x,y) \, d \lambda.
\end{equation}
The adjoints $R_{\mu}^{*}$ will then be given by
\begin{equation}
  \label{eqtn:RieszAdjoints}
R_{\mu}^{*}f(x) = \int_{\R^{d}} K^{*}_{\mu}(x,y) f(y) \, dy = \int_{\R^{d}} K_{\mu}(y,x) f(y) \, dy.
\end{equation}
In particular, the singular kernel of $R_{\mu}^{*}$, denoted
$K_{\mu}^{*}$, satisfies $K_{\mu}^{*}(x,y) = K_{\mu}(y,x)$ for all $x, \,
y \in \R^{d}$.

\begin{lem}
  \label{lem:RieszKernelEst}
There exists $\varepsilon > 0$, independent of $p$ and depending on
$\mu$ only through $C_{\mu}, \, D_{\mu}$ and $\delta_{\mu}$, for which
\begin{equation}
  \label{eqtn:PointwiseRieszKernel}
\abs{K_{\mu}^{*}(x,y)} \lesssim \frac{e^{-\varepsilon \cdot
    d_{\mu}(x,y)}}{\abs{x - y}^{d-1}} \br{\int_{B(y,\abs{x - y}/2)}
  \frac{d \mu(z)}{\abs{z - y}^{d - 1}} + \frac{1}{\abs{x - y}}}
\end{equation}
for all $x, \, y \in \R^{d}$.
  \end{lem}

\begin{prof}
 First note that from the definition of
the Agmon distance
\begin{equation}
  \label{eqtn:SumAgmon}
d_{\mu + \lambda}(x,y) \geq \frac{1}{2} \br{d_{\mu}(x,y) + d_{\lambda}(x,y)}
\end{equation}
for all $x, \, y \in \R^{d}$. Combining this with
\eqref{eqtn:DerivativeEstimate} from Theorem
\ref{thm:MeasureExp} gives
\begin{equation}
  \label{eqtn:Riesz1}
  \abs{\nabla \Gamma_{\mu + \lambda}(y,x)} \lesssim e^{- \varepsilon
    \cdot d_{\mu}(x,y)} e^{- \varepsilon \cdot d_{\lambda}(x,y)}
  \frac{1}{\abs{x - y}^{d - 2}} \br{\int_{B(y,\abs{x - y}/2)}
    \frac{d\mu(z)}{\abs{z - y}^{d - 1}} + \frac{1}{\abs{x - y}}}.
\end{equation}
for some $\varepsilon > 0$ independent of $p$ and depending on $\mu$
only through $C_{\mu}, \, D_{\mu}$ and $\delta_{\mu}$. Since $d_{\lambda}(x,y) =
\lambda^{\frac{1}{2}} \abs{x - y}$, we then have by
\eqref{eqtn:KernelExpression} and \eqref{eqtn:RieszAdjoints}
\begin{align*}\begin{split}  
 \abs{K_{\mu}^{*}(x,y)} &\lesssim \frac{e^{-\varepsilon \cdot
     d_{\mu}(x,y)}}{\abs{x - y}^{d - 2}} \br{\int_{B(y,\abs{x - y}/2)}
 \frac{d\mu(z)}{\abs{z - y}^{d - 1}} + \frac{1}{\abs{x - y}}} \int_{0}^{\infty}
 \lambda^{-\frac{1}{2}} e^{- \varepsilon \cdot \lambda^{\frac{1}{2}}
   \abs{x - y}} \, d \lambda \\
 &\simeq \frac{ e^{- \varepsilon \cdot d_{\mu}(x,y)}}{\abs{x - y}^{d-1}} \br{\int_{B(y,\abs{x - y}/2)}
 \frac{d\mu(z)}{\abs{z - y}^{d - 1}} + \frac{1}{\abs{x - y}}}
 \int^{\infty}_{0} \lambda^{-\frac{1}{2}} e^{- \varepsilon \lambda^{\frac{1}{2}}}
 \, d \lambda \\
 &\lesssim \frac{e^{- \varepsilon \cdot d_{\mu}(x,y)}}{\abs{x - y}^{d
     - 1}}\br{\int_{B(y,\abs{x - y}/2)}
 \frac{d\mu(z)}{\abs{z - y}^{d - 1}} + \frac{1}{\abs{x - y}}}.
\end{split}\end{align*}
\end{prof}

 As stated previously, in order to prove the $L^{p}(w)$-boundedness of an operator it
is sufficient to prove the boundedness of the global and local
components separately. The following proposition is a measure
generalisation of the local boundedness result proved in 
{\cite[Thm.~1]{bongioanni2011classes}}.

  \begin{prop} 
  \label{prop:MeasureLocalized}
Let $\mu$ be a non-negative Radon measure on $\R^{d}$
that satisfies \eqref{eqtn:Measure1} and
  \eqref{eqtn:Measure2} with constants $C_{\mu}, \, D_{\mu}, \,
  \delta_{\mu} > 0$.
  \begin{enumerate}[(i)]
\item If $\delta_{\mu} > 1$ then the operators $R_{\mu,loc}$ and
  $R_{\mu,loc}^{*}$ are bounded on $L^{p}(w)$ for $1 < p < \infty$ and
  $w \in A_{p}^{\mu,loc}$.
  \item Suppose instead that $0 < \delta_{\mu} < 1$ and let $\eta \in
    (2, (2 - \delta_{\mu})/(1 - \delta_{\mu}))$. Then
    $R^{*}_{\mu,loc}$ will be bounded on $L^{p}(w)$ for $\eta' < p <
    \infty$ when $w \in A_{p/\eta'}^{\mu,loc}$ and $R_{\mu,loc}$ will be
    bounded on $L^{p}(w)$ for $1 < p < \eta$ when $w^{-\frac{1}{p - 1}}
    \in A_{p'/\eta'}^{\mu,loc}$.
  \end{enumerate}
\end{prop}

\begin{prof}
 For any weight $w$ on $\R^{d}$, $1 < p < \infty$ and $f \in L^{p}(w)$ we have
  $$
\norm{R_{\mu,loc}^{*}f}_{L^{p}(w)} \lesssim \norm{(R_{\mu,loc}^{*} -
  R_{0,loc}^{*})f}_{L^{p}(w)} + \norm{R_{0,loc}^{*}f}_{L^{p}(w)}.
$$
As the operator $R_{0,loc}^{*}$ is bounded on $L^{p}(w)$ for any $w \in
A_{p}^{\mu,loc}$ by {\cite[Thm.~1]{bongioanni2011classes}}, it
suffices to prove that the difference term
$R_{\mu,loc}^{*} - R_{0,loc}^{*}$ is bounded on $L^{p}(w)$.
  In the proof of {\cite[Lem.~7.13]{shen1999fundamental}}, it was proved
that for $x, \, y \in \R^{d}$ satisfying $\abs{x - y} \leq \rho_{\mu}(x)$,
$$
\abs{K_{\mu}^{*}(x,y) - K_{0}^{*}(x,y)} \lesssim \frac{1}{\abs{x -
    y}^{d - 1}}
\int_{B(y,\abs{x - y}/2)} \frac{d \mu(z)}{\abs{z - y}^{d - 1}} +
\frac{1}{\abs{x - y}^{d}} \br{\frac{\abs{x - y}}{\rho_{\mu}(x)}}^{\delta_{\mu}}
$$
for all $x, \, y \in \R^{d}$. This gives
\begin{align*}\begin{split}  
 \norm{\br{R_{\mu,loc}^{*} - R_{0,loc}^{*}}f}_{L^{p}(w)}^{p} &\lesssim
 \int_{\R^{d}} \br{\int_{B(x,\rho_{\mu}(x))} \abs{K^{*}_{\mu}(x,y) -
     K_{0}^{*}(x,y)} \abs{f(y)} \, dy}^{p} w(x) \, dx \\
 &\lesssim \int_{\R^{d}} \br{\int_{B(x,\rho_{\mu}(x))} \frac{1}{\abs{x
    - y}^{d}} \br{\frac{\abs{x - y}}{\rho_{\mu}(x)}}^{\delta_{\mu}}
\abs{f(y)} \, dy}^{p} w(x) \, dx  \\
& \qquad + \int_{\R^{d}} \br{\int_{B(x,\rho_{\mu}(x))} \frac{1}{\abs{x
     - y}^{d - 1}} \int_{B(y,\abs{x - y}/2)} \frac{d \mu(z)}{\abs{z -
     y}^{d - 1}} \abs{f(y)} \, dy}^{p} w(x) \, dx \\
&=: J_{1} + J_{2}.
\end{split}\end{align*}
Using the same argument that is used to bound the function $h_{1}(x)$
in {\cite[Thm.~3]{bongioanni2011classes}}, it is clearly true that
$$
J_{1} \lesssim \norm{f}_{L^{p}(w)}
$$
for any $w \in A_{p}^{\mu,loc}$. It therefore suffices to estimate the
term $J_{2}$.

\vspace*{0.1in}

\underline{\textit{Proof of Part (i).}}
Suppose that $\delta_{\mu} > 1$. For this case, Lemma \ref{lem:RHn2},
\eqref{eqtn:Measure1} and Remark \ref{rmk:Critical} imply that for $y \in B(x,\rho_{\mu}(x))$,
\begin{align}\begin{split}
    \label{eqtn:J1J2}
 \int_{B(y, \abs{x - y} / 2)} \frac{d \mu(z)}{\abs{z - y}^{d - 1}}
 &\lesssim \frac{\mu(B(y,\abs{x - y}/2))}{\abs{x - y}^{d - 1}} \\
 &\lesssim \frac{\mu(B(x,2 \abs{x - y}))}{\abs{x - y}^{d - 1}} \\
 &\lesssim \br{\frac{\abs{x - y}}{\rho_{\mu}(x)}}^{d - 2 +
   \delta_{\mu}} \frac{\mu(B(x,2 \rho_{\mu}(x)))}{\abs{x - y}^{d - 1}}
 \\
 &\lesssim \frac{1}{\abs{x - y}} \br{\frac{\abs{x - y}}{\rho_{\mu}(x)}}^{\delta_{\mu}}
 \end{split}\end{align}
and therefore $J_{2} \lesssim J_{1}$.
Since we already know that $J_{1}$ is bounded, this implies that
$R_{\mu,loc}^{*}$ is bounded on $L^{p}(w)$ for $w \in
A_{p}^{\mu,loc}$. The boundedness of $R_{\mu,loc}$ on $L^{p}(w)$
follows from duality. This proves the first
part of our proposition.

\vspace*{0.1in}

\underline{\textit{Proof of Part (ii).}} Suppose that $0 <
\delta_{\mu} < 1$. For this case, Lemma \ref{lem:RHn2} can no longer
be applied to bound the term $J_{2}$. Instead, $J_{2}$ will be handled
by adapting the argument from
{\cite[Thm.~3]{bongioanni2011classes}}. Let $\eta \in (0,(2 -
\delta_{\mu})/(1 - \delta_{\mu}))$, $\eta' < p < \infty$ and assume that $w \in A_{p/\eta'}^{\mu,loc}$. Let $\lb B_{j} \rb_{j \in \N}
= \lb B(x_{j},\rho_{\mu}(x_{j})) \rb_{j \in \N}$
be a covering of $\R^{d}$ by balls as given in Proposition
\ref{prop:Cover}. For each $j, \, k \in \N$, there exists $2^{d k}$
balls $B^{j,k}_{l} = B(x^{j,k}_{l}, 2^{-k} \rho_{\mu}(x_{j}))$, $l = 1,
\cdots, 2^{d k}$, with the properties that $B_{j} \subset \cup^{2^{d
    k}}_{l = 1} B_{l}^{j,k} \subset 2 B_{j}$ and $\sum_{l = 1}^{2^{d
    k}} \chi_{B^{j,k}_{l}} \leq 2^{d}$. Set $\tilde{B}_{l}^{j,k} = 10
B_{0} \cdot B_{l}^{j,k}$ with $B_{0}$ as given in Lemma
\ref{lem:CriticalMeasure}. The construction can be done in such a way so
that
\begin{equation}
  \label{eqtn:RieszOverlap}
\sum_{j} \sum_{l = 1}^{2^{d k}} \mathbbm{1}_{\tilde{B}_{l}^{j,k}} \leq C
\end{equation}
for some $C > 0$ independent of $k$. Define the function
$$
h_{2}(x) := \int_{B(x,\rho_{\mu}(x))} \frac{\abs{f(y)}}{\abs{x - y}^{d
  - 1}} \br{\int_{B(y,\abs{x - y}/2)} \frac{d \mu(z)}{\abs{z - y}^{d -
  1}}} \, dy.
$$
Then we have
\begin{equation}
  \label{eqtn:LocalAnn}
h_{2}(x) \lesssim \sum_{k = 0}^{\infty} 2^{k(d - 1)} h_{2,k}(x),
\end{equation}
where
$$
h_{2,k}(x) := \rho_{\mu}(x)^{-d + 1} \int_{B(x,2^{-k}\rho_{\mu}(x))}
\abs{f(y)} \br{\int_{B(y,\abs{x - y} / 2)} \frac{d \mu(z)}{\abs{z -
      y}^{d - 1}}} \, dy.
$$
If $x \in B_{l}^{j,k}$ then
\begin{align*}\begin{split}  
 h_{2,k}(x) &\lesssim \rho_{\mu}(x_{j})^{-d + 1}
\int_{\tilde{B}_{l}^{j,k}} \abs{f(y)} \br{\int_{\tilde{B}_{l}^{j,k}}
  \frac{d \mu(z)}{\abs{z - y}^{d - 1}}} \, dy  \\
&\lesssim \rho_{\mu}(x_{j})^{-d + 1} \norm{\int_{\tilde{B}_{l}^{j,k}}
  \frac{d\mu(z)}{\abs{z - \cdot}}}_{L^{\eta}\br{\tilde{B}_{l}^{j,k}}}
\norm{f}_{L^{p}(\tilde{B}_{l}^{j,k},w)} \br{\int_{\tilde{B}_{l}^{j,k}}
w^{-\gamma / p}}^{\frac{1}{\gamma}},
 \end{split}\end{align*}
where the second line follows from H\"{o}lder's inequality and
$\frac{1}{\gamma} := 1 - \frac{1}{p} - \frac{1}{\eta}$. Lemma 7.9 of
\cite{shen1999fundamental}, \eqref{eqtn:Measure1} and Remark \ref{rmk:Critical} imply that
\begin{align*}\begin{split}  
 \norm{\int_{\tilde{B}_{l}^{j,k}} \frac{d\mu(z)}{\abs{z - \cdot}^{d -
      1}}}_{L^{\eta}(\tilde{B}_{l}^{j,k})} &\lesssim \frac{\mu(3
  \tilde{B}_{l}^{j,k})}{(2^{-k} \rho_{\mu}(x_{j}))^{\frac{d}{\eta'} - 1}}
\\
&\lesssim (2^{-k})^{d - 1 + \delta_{\mu} - \frac{d}{\eta'}} \rho_{\mu}(x_{j})^{d - 1 - \frac{d}{\eta'}}.
\end{split}\end{align*}
Since $w \in A_{p/\eta'}^{\mu,loc}$ we have
\begin{align*}\begin{split}  
 w(\tilde{B}_{l}^{j,k})^{\frac{1}{p}} \br{\int_{\tilde{B}_{l}^{j,k}}
   w^{-\frac{\gamma}{p}}}^{\frac{1}{\gamma}} &\lesssim
 \abs{\tilde{B}_{l}^{j,k}}^{\frac{1}{\eta'}} \\
 &\lesssim (2^{-k} \rho_{\mu}(x_{j}))^{\frac{d}{\eta'}}.
\end{split}\end{align*}
Therefore,
\begin{align*}\begin{split}  
 \norm{h_{2,k}}_{L^{p}(w)}^{p} &\lesssim  \sum_{j,l} \int_{B_{l}^{j,k}}
 \br{\rho_{\mu}(x_{j})^{-d + 1} \norm{\int_{\tilde{B}_{l}^{j,k}} \frac{d\mu(z)}{\abs{z - \cdot}^{d -
      1}}}_{L^{\eta}(\tilde{B}_{l}^{j,k})}
\norm{f}_{L^{p}(\tilde{B}_{l}^{j,k},w)}w^{-\frac{\gamma}{p}}(\tilde{B}_{l}^{j,k})^{\frac{1}{\gamma}}}^{p}
w(x) \, dx \\
&\lesssim \sum_{j,l} \br{\rho_{\mu}(x_{j})^{-d + 1} \norm{\int_{\tilde{B}_{l}^{j,k}} \frac{d\mu(z)}{\abs{z - \cdot}^{d -
      1}}}_{L^{\eta}(\tilde{B}_{l}^{j,k})}
\norm{f}_{L^{p}(\tilde{B}_{l}^{j,k},w)}
w(\tilde{B}_{l}^{j,k})^{\frac{1}{p}}
w^{-\frac{\gamma}{p}}(\tilde{B}_{l}^{j,k})^{\frac{1}{\gamma}}}^{p} \\
&\lesssim \sum_{j,l} \br{\rho_{\mu}(x_{j})^{-d + 1}
  \rho_{\mu}(x_{j})^{d - 1 - \frac{d}{\eta'}} (2^{-k})^{d - 1 +
    \delta_{\mu} - \frac{d}{\eta'}} (2^{-k}
  \rho_{\mu}(x_{j}))^{\frac{d}{\eta'}}
  \norm{f}_{L^{p}(\tilde{B}_{l}^{j,k},w)}}^{p} \\
&\lesssim 2^{-k p (d - 1 + \delta_{\mu})} \norm{f}_{L^{p}(w)}^{p},
\end{split}\end{align*}
where \eqref{eqtn:RieszOverlap} was used to obtain the final line.
Referring back to \eqref{eqtn:LocalAnn},
\begin{align*}\begin{split}  
 J_{2} = \norm{h_{2}}_{L^{p}(w)} &\lesssim \sum_{k = 0}^{\infty}
 2^{k(d - 1)} \norm{h_{2,k}}_{L^{p}(w)} \\
 &\lesssim \norm{f}_{L^{p}(w)} \br{\sum_{k = 0}^{\infty} 2^{k(d - 1)}
   2^{-k (d - 1 + \delta_{\mu})}} \\
 &= \norm{f}_{L^{p}(w)} \br{\sum_{k = 0}^{\infty} 2^{-k \delta_{\mu}}}
 \\
 &\lesssim \norm{f}_{L^{p}(w)}.
\end{split}\end{align*}
This proves that $R_{\mu,loc}^{*}$ is bounded on $L^{p}(w)$ for $\eta'
< p < \infty$ and $w \in A_{p/\eta'}^{\mu,loc}$. The
$L^{p}(w)$-boundedness of $R_{\mu,loc}$ for $1 < p < \eta$ and
$w^{-\frac{1}{p - 1}} \in A_{p'/\eta'}^{\mu,loc}$ follows by using duality.
\end{prof}

With the boundedness of the local component of our operators
established, it now suffices to consider the boundedness of the global
components.

\begin{thm} 
 \label{thm:RieszGlobal} 
 Let $\mu$ be a non-negative Radon measure on $\R^{d}$
 that satisfies \eqref{eqtn:Measure1} and \eqref{eqtn:Measure2} with
 constants $C_{\mu}, \, D_{\mu}, \, \delta_{\mu} > 0$.
 \begin{enumerate}[(i)]
\item Suppose that $\delta_{\mu} > 1$. There exists $c_{1} > 0$ such that $R_{\mu,glob}$ and $R^{*}_{\mu,glob}$ are both bounded on
 $L^{p}(w)$ for any $w \in S_{p,c_{1}}^{\mu}$ and $1 < p < \infty$. \\
 \item Suppose instead that $0 < \delta_{\mu} < 1$ and let $\eta \in
   (2, (2 - \delta_{\mu}) / (1 - \delta_{\mu}))$. There exists
    $c_{2} > 0$ for which the operator $R^{*}_{\mu,glob}$ is bounded
    on $L^{p}(w)$ for $\eta' < p < \infty$ when $w \in
    S_{p/\eta',c_{2}}^{\mu}$ and the operator $R_{\mu}$ is bounded on
    $L^{p}(w)$ for $1 < p < \eta$ when $w^{-\frac{1}{p - 1}} \in S_{p'/\eta',c_{2}}^{\mu}$.
   \end{enumerate}
The constants $c_{1}$ and $c_{2}$ are independent of $p$ and depend on
$\mu$ only through $C_{\mu}, \, D_{\mu}$ and $\delta_{\mu}$.
 \end{thm}

\begin{prof}  
 For $1 < p < \infty$, weight $w$ on $\R^{d}$ and $f \in L^{p}(w)$,
\begin{align*}\begin{split}
 \norm{R_{\mu,glob}^{*} f}_{L^{p}(w)}^{p} &= \int_{\R^{d}}
 \abs{R_{\mu,glob}^{*}f(x)}^{p} w(x) \, dx \\
 &\leq \int_{\R^{d}} \br{\int_{B(x,\rho_{\mu}(x))^{c}} \abs{K_{\mu}^{*}(x,y)}
   \abs{f(y)} \, dy}^{p}
 w(x) \, dx.
\end{split}\end{align*}
Let $\lb x_{j} \rb_{j \in \N}$ be a collection of points in $\R^{d}$
as given in Proposition \ref{prop:Cover}. Introduce the notation
$B_{j} := B(x_{j}, \rho_{\mu}(x_{j}))$ for $j \in \N$.  
Since the collection of balls $\lb B_{j} \rb_{j \in \N}$ forms
a cover for $\R^{d}$,
\begin{equation} 
      \label{eqtn:Riesz101}
  \norm{R_{glob}^{*} f}_{L^{p}(w)}^{p}  \leq \sum_{j} \int_{B_{j}} \br{\int_{B(x,\rho_{\mu}(x))^{c}} \abs{K_{\mu}^{*}(x,y)} \abs{f(y)} \,
   dy}^{p} w(x) \, dx.
 \end{equation}
 The
estimate \eqref{eqtn:PointwiseRieszKernel} from Lemma \ref{lem:RieszKernelEst} then leads to
\begin{align*}\begin{split}  
    & \norm{R^{*}_{\mu,glob}f}^{p}_{L^{p}(w)} \lesssim J_{1} + J_{2} := \sum_{j}
 \int_{B_{j}} \br{\int_{B(x,\rho_{\mu}(x))^{c}} \frac{e^{-\varepsilon
       d_{\mu}(x,y)}}{\abs{x - y}^{d}} \abs{f(y)} \, dy}^{p} w(x) \,
 dx \\
 & \qquad \qquad + \sum_{j} \int_{B_{j}} \br{\int_{B(x,\rho_{\mu}(x))^{c}} \frac{e^{-
       \varepsilon d_{\mu}(x,y)}}{\abs{x - y}^{d - 1}}
   \int_{B(y,\abs{x - y} / 2)} \frac{d\mu(z)}{\abs{z - y}^{d - 1}}
   \abs{f(y)} \, dy}^{p} w(x) \, dx,
\end{split}\end{align*}
for some $\varepsilon > 0$ independent of $p$ and depending on $\mu$
only through $C_{\mu}$, $D_{\mu}$ and $\delta_{\mu}$.

\vspace*{0.1in}

\textit{\underline{Proof of Part (i).}} Assume that
$\delta_{\mu} > 1$.  Fix
$1 < p < \infty$ and $w
\in S_{p,c_{1}}^{V}$ for some $c_{1} > 0$. Lemma \ref{lem:RHn2} followed by \eqref{eqtn:Measure1} implies that
\begin{align*}\begin{split}  
 \frac{e^{- \varepsilon d_{\mu}(x,y)}}{\abs{x - y}^{d - 1}}
 \int_{B(y,\abs{x - y} / 2)} \frac{d \mu(z)}{\abs{z - y}^{d - 1}}
 &\lesssim \frac{e^{- \varepsilon d_{\mu}(x,y)}}{\abs{x - y}^{d - 1}}
 \frac{\mu(B(y,\abs{x - y} / 2))}{\abs{x - y}^{d - 1}} \\
 &\lesssim \frac{e^{- \varepsilon d_{\mu}(x,y)}}{\abs{x - y}^{d - 1}}
 \br{\frac{\abs{x - y}}{\rho_{\mu}(y)}}^{d - 2 + \delta_{\mu}}
 \frac{\mu(B(y,\rho(y)))}{\abs{x - y}^{d - 1}}.
\end{split}\end{align*}
Remark \ref{rmk:Critical} together
with Lemma \ref{lem:Shen} gives
\begin{align*}\begin{split}  
  \frac{e^{- \varepsilon d_{\mu}(x,y)}}{\abs{x - y}^{d - 1}}
 \int_{B(y,\abs{x - y} / 2)} \frac{d \mu(z)}{\abs{z - y}^{d - 1}}
 &\lesssim \frac{e^{-\varepsilon d_{\mu}(x,y)}}{\abs{x - y}^{d}}
 \br{\frac{\abs{x - y}}{\rho_{\mu}(y)}}^{\delta_{\mu}} \\
 &\lesssim \frac{e^{-\frac{\varepsilon}{2} d_{\mu}(x,y)}}{\abs{x - y}^{d}}.
 \end{split}\end{align*}
This proves that $J_{2}$ is bounded from above by a term that is
identical to $J_{1}$ except that the exponent of the exponential is $\varepsilon /
2$ instead of $\varepsilon$. Due to this similarity, in order to prove
our claim it will then suffice
to show that $c_{1}$ can be set small enough so that
$J_{1} \lesssim \norm{f}_{L^{p}(w)}$.

Notice that Lemma \ref{lem:Balls2} tells us that $B(x, \rho_{\mu}(x))^{c} \subseteq
B_{\mu}(x,A_{0}^{-1})^{c}$. Let's use the shorthand notation $B_{\mu,x}$ to denote
the ball $B_{\mu}\br{x,A_{0}^{-1}}$. Then
\begin{align*}\begin{split}  
 J_{1} &\leq \sum_{j} \int_{B_{j}} \br{\int_{B_{\mu,x}^{c}}
   \frac{e^{-\varepsilon d_{\mu}(x,y)}}{\abs{x - y}^{d}} \abs{f(y)} \,
   dy}^{p} w(x) \, dx \\
 &= \sum_{j} \int_{B_{j}} \br{\sum_{k = 1}^{\infty} \int_{(k + 1)
     B_{\mu,x} \setminus k B_{\mu,x}} \frac{e^{- \varepsilon
       d_{\mu}(x,y)}}{\abs{x - y}^{d}} \abs{f(y)} \, dy}^{p} w(x) \, dx,
 \end{split}\end{align*}
where $k B_{\mu,x} := B_{\mu}(x,
k A_{0}^{-1})$ for $k \geq 1$. 
Fix $x \in B_{j}$ for some $j \in \N$ and $y \in (k + 1) B_{\mu,x}
 \setminus k B_{\mu,x}$ for some $k \geq 1$. Suppose first that $k  \leq 2
A_{0} \beta$. Lemma \ref{lem:Balls} will then imply that 
$$
\br{k B_{\mu,x}}^{c} = B_{\mu}(x, k A_{0}^{-1})^{c} \subset B \br{x,
  \frac{k}{\beta A_{0}} \rho_{\mu}(x)}^{c}.
$$
Therefore,
$$
\abs{x - y}^{-d} \lesssim k^{-d} \rho_{\mu}(x)^{-d} \leq \rho_{\mu}(x)^{-d}.
$$
Next, suppose that $k  > 2 \beta A_{0}$. For this case, Lemma
\ref{lem:Balls} implies that
$$
(k B_{\mu,x})^{c} \subset B(x, 2 \rho_{\mu}(x))^{c}
$$
and therefore
$$
\abs{x - y}^{-d} \lesssim \rho_{\mu}(x)^{-d}.
$$
This produces the estimate,
$$
 \frac{e^{-\varepsilon d_{\mu}(x,y)}}{\abs{x - y}^{d}}
\lesssim  e^{- \delta k}  \rho_{\mu}(x)^{-d}
$$
for any $k \geq 1$, where $\delta := \varepsilon A_{0}^{-1}$.
Since $x \in B_{j}$, we have by Lemma \ref{lem:CriticalMeasure}
$$
\rho_{\mu}(x)^{-d} \lesssim \rho_{\mu}(x_{j})^{-d} \br{1 + \frac{\abs{x -
      x_{j}}}{\rho_{\mu}(x_{j})}}^{d k_{0}} \lesssim \rho_{\mu}(x_{j})^{-d},
$$
which implies that
$$
\frac{e^{-\varepsilon d_{\mu}(x,y)}}{\abs{x - y}^{d}} \lesssim e^{- \delta k} 
\rho_{\mu}(x_{j})^{-d}.
$$
Applying this estimate to $J_{1}$ gives
\begin{equation}
  \label{eqtn:Riesz11}
 J_{1} \lesssim \sum_{j} \int_{B_{j}}
 \br{\sum_{k = 1}^{\infty} e^{- \delta k} 
\rho_{\mu}(x_{j})^{-d} \int_{(k+1) B_{\mu,x}
  \setminus k B_{\mu,x}} \abs{f(y)} \, dy}^{p} w(x) \, dx 
\end{equation}
Define $B_{\mu,j} := B_{\mu}(x_{j}, \beta + A_{0}^{-1})$ for $j \in
\N$. Notice that for $x \in B_{j}$, Lemma \ref{lem:LocalAgmon} implies
that $d_{\mu}(x,x_{j}) \leq \beta$. Therefore, for
$k \geq 1$ and $y \in (k + 1)B_{\mu,x}$,
\begin{align*}\begin{split}  
    d_{\mu}(x_{j},y) &\leq d_{\mu}(x,x_{j}) + d_{\mu}(x,y) \\
    &\leq \beta + A_{0}^{-1} (k + 1) \\
    &\leq \br{\beta + A_{0}^{-1}} (k + 1),
  \end{split}\end{align*}
which implies that $(k + 1)B_{\mu,x} \subset (k + 1)
B_{\mu,j}$. Therefore,
\begin{align*}\begin{split}  
 J_{1} &\lesssim \sum_{j} \int_{B_{j}}
 \br{\sum_{k = 1}^{\infty} e^{- \delta k} 
\rho_{\mu}(x_{j})^{-d} \int_{(k+1) B_{\mu,j}} \abs{f(y)} \, dy}^{p} w(x) \, dx   \\
&\lesssim \sum_{j} \int_{B_{\mu,j}}
 \br{\sum_{k = 1}^{\infty} e^{- \delta k} 
\rho_{\mu}(x_{j})^{-d} \int_{(k+1) B_{\mu,j}} \abs{f(y)} \, dy}^{p} w(x) \, dx,
\end{split}\end{align*}
where the last line follows from the inclusion $B_{j} \subset B_{\mu,j}$
by Lemma \ref{lem:Balls}. H\"{o}lder's inequality then leads to,
\begin{align}\begin{split}
    \label{eqtn:Riesz2}
 J_{1} 
&\lesssim \sum_{j} \int_{B_{\mu,j}}
 \br{\sum_{k = 1}^{\infty} e^{- \delta k} 
\rho_{\mu}(x_{j})^{-d} 
\norm{f}_{L^{p}((k+1) B_{\mu,j},w)} w^{-\frac{1}{p-1}}\br{(k+1)
  B_{\mu,j}}^{\frac{p - 1}{p}}}^{p} w(x) \, dx \\
&\lesssim \sum_{j} 
 \br{\sum_{k = 1}^{\infty} e^{- \delta k} 
\rho_{\mu}(x_{j})^{-d} 
\norm{f}_{L^{p}((k+1) B_{\mu,j},w)} w^{-\frac{1}{p-1}}\br{(k+1)
  B_{\mu,j}}^{\frac{p-1}{p}} w(B_{\mu,j})^{\frac{1}{p}}}^{p}.
\end{split}\end{align}
Since $w \in S_{p,c_{1}}^{\mu}$, we have the estimate
\begin{align}\begin{split}
    \label{eqtn:Riesz21}
 w^{-\frac{1}{p-1}}\br{(k + 1)B_{\mu,j}}^{\frac{p-1}{p}}
 w(B_{\mu,j})^{\frac{1}{p}} &\leq w^{-\frac{1}{p-1}}\br{(k +
   1)B_{\mu,j}}^{\frac{p-1}{p}}
 w((k+1)B_{\mu,j})^{\frac{1}{p}} \\
 &\lesssim \abs{(k + 1)B_{\mu,j}} e^{c'(k + 1)},
\end{split}\end{align}
where $c' := c_{1} (\beta + A_{0}^{-1})$.
For $j, \, k \in \N$, let $B_{j,k}$ denote ball
$$
B_{j,k} := B\br{x_{j}, \br{\beta (k + 1)(\beta + A_{0}^{-1})}^{k_{0} + 1} \rho_{\mu}(x_{j})}.
$$
Lemma \ref{lem:Balls2} can then be applied to obtain
$$
(k + 1)B_{\mu,j} \subset  B_{j,k}
$$
and therefore
\begin{equation}
  \label{eqtn:BallMeasure}
\abs{(k + 1) B_{\mu,j}} \lesssim  (k + 1)^{d(k_{0} + 1)} \rho_{\mu}(x_{j})^{d}.
\end{equation}
Combining this with \eqref{eqtn:Riesz21} and  \eqref{eqtn:Riesz2}
leads to
\begin{align*}\begin{split}  
 J_{1} &\lesssim \sum_{j} 
 \br{\sum_{k = 1}^{\infty} e^{- \delta k} 
\rho_{\mu}(x_{j})^{-d} 
\norm{f}_{L^{p}((k+1) B_{\mu,j},w)} (k + 1)^{d(k_{0} + 1)}
\rho_{\mu}(x_{j})^{d} e^{c' k}}^{p} \\
&\lesssim \sum_{j} \br{\sum_{k = 1}^{\infty} k^{d (k_{0} + 1)} e^{(c' - \delta) k}
  \norm{f}_{L^{p}(B_{j,k},w)}}^{p}.
\end{split}\end{align*}
The bounded overlap property of the balls $B_{j}$, as given in
Proposition \ref{prop:Cover}, then gives
\begin{align*}\begin{split}  
J_{1}^{\frac{1}{p}} &\lesssim \sum_{k = 1}^{\infty}
 k^{d(k_{0} + 1)} e^{(c' - \delta)k} \br{\sum_{j} \norm{f}_{L^{p}(B_{j,k},w)}^{p}}^{\frac{1}{p}} \\
 &\lesssim \br{\sum_{k = 1}^{\infty} e^{(c' - \delta) k}
   k^{\br{\frac{N_{1}}{p} + d}(k_{0} + 1)}}
 \norm{f}_{L^{p}(\R^{d},w)} \\
 &\lesssim \norm{f}_{L^{p}(w)},
\end{split}\end{align*}
so long as we choose $c' = c_{1} (\beta + A_{0}^{-1}) < \varepsilon
A_{0}^{-1} =  \delta$. It remains to prove the boundedness of the global part of the Riesz
transforms $R_{\mu}$. However, this follows from the boundedness of
$R^{*}_{\mu,glob}$ using duality.

\vspace*{0.1in}

\textit{\underline{Proof of Part (ii).}} Assume that
$0 < \delta_{\mu} < 1$. The term $J_{1}$ can be handled in an
identical manner to the case $\delta_{\mu} > 1$. It therefore suffices
to demonstrate that there exists $c_{2} > 0$ for which the term
$J_{2}$ is bounded on $L^{p}(w)$ for $\eta' < p < \infty$ when
$w \in S_{p/\eta',c_{2}}^{\mu}$. From reasoning identical to the case of
$J_{1}$,
\begin{align*}\begin{split}  
J_{2} &\leq \sum_{j} \int_{B_{j}} \br{\sum_{k = 1}^{\infty} \int_{(k + 1)
  B_{\mu,x} \setminus k B_{\mu,x}} \frac{e^{- \varepsilon
    d_{\mu}(x,y)}}{\abs{x - y}^{d - 1}} \int_{B(y,\abs{x - y}/2)}
\frac{d \mu(z)}{\abs{z - y}^{d - 1}} \abs{f(y)} \, dy}^{p} w(x) \, dx
\\
&\lesssim \sum_{j} \int_{B_{j}} \br{\sum_{k = 1}^{\infty} e^{- \delta
    k} \rho_{\mu}(x_{j})^{- (d - 1)} \int_{(k + 1) B_{\mu,x}} \br{\int_{B(y,\abs{x - y}/2)} \frac{d \mu(z)}{\abs{z -
    y}^{d -1}}} \abs{f(y)} \, dy}^{p} w(x) \, dx.
\end{split}\end{align*}
For any $x \in B_{j}$, H\"{o}lder's inequality implies that for
$\frac{1}{\gamma} := 1 - \frac{1}{p} - \frac{1}{\eta}$,
\begin{align}\begin{split}
    \label{eqtn:RestRange1}
& \int_{(k + 1) B_{\mu,x}} \br{\int_{B(y,\abs{x - y}/2)} \frac{d
   \mu(z)}{\abs{z -  y}^{d - 1}}} \abs{f(y)} \, dy \\ & \qquad \qquad
\lesssim \norm{\int_{B(\cdot, \abs{x - \cdot}/2)}
 \frac{d \mu(z)}{\abs{z - \cdot}^{d - 1}}}_{L^{\eta}(B_{j,k})}
\norm{f}_{L^{p}(B_{j,k},w)} \br{\int_{(k + 1) B_{\mu,x}}
  w^{-\gamma/p}}^{\frac{1}{\gamma}} \\
& \qquad \qquad \lesssim  \norm{\int_{2 B_{j,k}} \frac{d \mu(z)}{\abs{z
      - \cdot}^{d - 1}}}_{L^{\eta}(2 B_{j,k})} \norm{f}_{L^{p}(B_{j,k},w)}
\br{\int_{(k + 1) B_{\mu,j}} w^{-\gamma/p}}^{\frac{1}{\gamma}},
\end{split}\end{align}
where we have used the inclusions $(k + 1)B_{\mu,x} \subset (k +
1)B_{\mu,j} \subset B_{j,k}$ and $B(y,\abs{x - y}/2) \subset 2
B_{j,k}$ for $y \in B_{j,k}$. For the first term in
\eqref{eqtn:RestRange1}, successively apply
{\cite[Lem.~7.9]{shen1999fundamental}}, the property
\eqref{eqtn:Measure2} and Remark \ref{rmk:Critical} to obtain
\begin{align}\begin{split}
    \label{eqtn:ShenLem7.9}
 \norm{\int_{2 B_{j,k}} \frac{d\mu(z)}{\abs{z - \cdot }^{d -
       1}}}_{L^{\eta}(2 B_{j,k},dx)} &\lesssim \frac{\mu(6 B_{j,k})}{\rho_{\mu}(x_{j})^{\frac{d}{\eta'} - 1}} \\
 &\lesssim k^{M'} \rho_{\mu}(x_{j})^{d - 1 - \frac{d}{\eta'}},
\end{split}\end{align}
for some $M' > 0$. This implies that
\begin{equation}
  \label{eqtn:RestRange2}
J_{2} \lesssim \sum_{j}  \br{\sum_{k = 1}^{\infty} e^{-
    \delta k} \rho_{\mu}(x_{j})^{-(d - 1)} k^{M'}
  \rho_{\mu}(x_{j})^{d - 1 - \frac{d}{\eta'}}
\norm{f}_{L^{p}(B_{j,k},w)} w^{-\frac{\gamma}{p}} \br{(k + 1)
  B_{\mu,j}}^{\frac{1}{\gamma}} w(B_{j})^{\frac{1}{p}}}^{p}.
\end{equation}
Our assumption $w \in S_{p/\eta',c_{2}}^{\mu}$ then allows us to bound
the term involving the weights by
\begin{align*}\begin{split}  
 w^{-\frac{\gamma}{p}} \br{(k + 1) B_{\mu,j}}^{\frac{1}{\gamma}}
 w(B_{j})^{\frac{1}{p}} &\leq  w^{-\frac{\gamma}{p}} \br{(k + 1) B_{\mu,j}}^{\frac{1}{\gamma}}
 w((k + 1)B_{\mu,j})^{\frac{1}{p}} \\
 &\lesssim \abs{(k + 1)B_{\mu,j}}^{\frac{1}{\eta'}} e^{\frac{c_{2}}{\eta'}(k + 1)(\beta
   + A_{0}^{-1})} \\
 &\lesssim k^{M''} \rho_{\mu}(x_{j})^{\frac{d}{\eta'}}
 e^{\frac{c_{2}}{\eta'}(k + 1)(\beta + A_{0}^{-1})} \\
 &\lesssim k^{M''} \rho_{\mu}(x_{j})^{\frac{d}{\eta'}} e^{c_{2}'k}.
\end{split}\end{align*}
for some $M'' > 0$, where $c_{2}' = c_{2}(\beta + A_{0}^{-1})$. On applying this to \eqref{eqtn:RestRange2},
\begin{align*}\begin{split}  
 J_{2} &\lesssim \sum_{j}  \br{\sum_{k = 1}^{\infty} e^{-
    \delta k}  k^{M'}
  \rho_{\mu}(x_{j})^{-\frac{d}{\eta'}}
\norm{f}_{L^{p}(B_{j,k},w)} k^{M''} \rho(x_{j})^{\frac{d}{\eta'}}
e^{c_{2}' k}}^{p} \\
&\lesssim \sum_{j}  \br{\sum_{k = 1}^{\infty} e^{- (\delta
    - c_{2}')k} k^{M' + M''} \norm{f}_{L^{p}(B_{j,k},w)}}^{p}.
\end{split}\end{align*}
The bounded overlap property of the balls $\lb B_{j} \rb_{j \in \N}$
will then imply that $J_{2} \lesssim \norm{f}_{L^{p}(w)}$ provided
that we choose $c_{2}$ small enough so that $c_{2}' = c_{2}(\beta +
A_{0}^{-1}) < \varepsilon A_{0}^{-1} = \delta$.

\end{prof}

Notice that in Theorem \ref{thm:MeasurePotential}, when $\delta_{\mu} \in (0,1)$, the
range of $p$ for which $R_{\mu}$ is bounded must be restricted to the
interval $(1,(2 - \delta_{\mu})/(1 - \delta_{\mu}))$. When $d\mu(x) = V(x) \, dx$ for some $V \in
RH_{q}$ with $\frac{d}{2} < q < d$, $\delta_{\mu} = 2 - \frac{d}{q}$
by Remark \ref{rmk:Measure}. Therefore, Theorem
\ref{thm:MeasurePotential} restricts the range of $p$ for which $R_{V}$ is bounded to 
$$
1 < p < \frac{d}{d - q}.
$$
Since $\frac{d}{d - q} < \br{\frac{1}{q} - \frac{1}{d}}^{-1}$, this range of $p$
is smaller than the range of $p$ given in
{\cite[Thm.~3]{bongioanni2011classes}}. The following proposition improves
the range of $p$ for the case $d \mu(x) = V(x) \, dx$. This
proposition,
when taken with the boundedness of $R_{V,loc}$ by {\cite[Thm.~3]{bongioanni2011classes}}, completes the
proof of the second part of Theorem \ref{thm:Riesz}.

\begin{prop} 
 \label{prop:ImprovedRange} 
 Suppose that $V \in RH_{q}$ for some $\frac{d}{2} < q < d$ and define
 $s$ through $\frac{1}{s} = \frac{1}{q} - \frac{1}{d}$. There exists
 $c > 0$ for which the operator $R_{V,glob}^{*}$ is bounded on
 $L^{p}(w)$ for $s' < p < \infty$ when $w \in S_{p/s',c}^{V}$ and the
 operator $R_{V}$ is bounded on $L^{p}(w)$ for $1 < p < s$ when
 $w^{-\frac{1}{p - 1}} \in S_{p'/s',c}^{V}$. The constant $c$ will
 depend on $V$ only through $\brs{V}_{RH_{\frac{d}{2}}}$ and will be independent
 of $p$.
\end{prop}

\begin{prof}  
 The proof is essentially identical to that of Theorem
 \ref{thm:RieszGlobal}.(ii). The only difference is that in
 \eqref{eqtn:ShenLem7.9}, Lemma 7.9 of \cite{shen1999fundamental}
 should no longer be used since this leads to a restricted range of
 $p$. Instead, the boundedness of the classical fractional integral operator of
 order one from $L^{q}$ into $L^{s}$ should be exploited as in the
 proof of {\cite[Thm.~3]{bongioanni2011classes}}. Specifically replace
 \eqref{eqtn:ShenLem7.9} with
 \begin{align*}\begin{split}  
 \norm{\int_{2 B_{j,k}} \frac{V(z)}{\abs{z - \cdot}^{d - 1}} \, dz}_{L^{s}(2
   B_{j,k},dx)} &= \norm{I_{0}^{1}(\mathbbm{1}_{2 B_{j,k}}
   V)}_{L^{s}(2 B_{j,k})} \\
 &\lesssim \norm{\mathbbm{1}_{2 B_{j,k}} V}_{q},
\end{split}\end{align*}
where $I_{0}^{1}$ is the classical fraction integral operator of order
$1$ and the well-known property that $I_{0}^{1}$ is bounded from
$L^{q}$ into $L^{s}$ is exploited. On applying the property $V \in
RH_{q}$, followed by \eqref{eqtn:Measure2} and Remark
\ref{rmk:Critical},
\begin{align*}\begin{split}  
 \norm{\int_{2 B_{j,k}} \frac{V(z)}{\abs{z - \cdot}^{d - 1}} \, dz}_{L^{s}(2
   B_{j,k},dx)} &\lesssim \abs{2 B_{j,k}}^{-\frac{1}{q'}} \int_{2
   B_{j,k}} V \\
 &\lesssim k^{M'} \abs{2 B_{j,k}}^{-\frac{1}{q'}} \int_{B_{j}} V \\
 &\lesssim k^{M'} \rho_{\mu}(x_{j})^{-\frac{d}{q'}}
 \rho_{\mu}(x_{j})^{d - 2} \\
 &= k^{M'} \rho_{\mu}(x_{j})^{d - 1 - \frac{d}{s'}},
 \end{split}\end{align*}
for some $M' > 0$.
The rest of the proof proceeds
identically to Theorem \ref{thm:RieszGlobal}.
 \end{prof}

 \subsection{The Riesz Potentials}
 \label{subsec:RieszPotentials}

Next, let's consider the Riesz potentials for the operator
$L_{\mu}$ and prove the third part of Theorem \ref{thm:MeasurePotential}. Notice that the pointwise estimate
$$
\abs{I^{\alpha}_{\mu} f(x)} \leq I^{\alpha}_{0} \abs{f}(x)
$$
holds for all $f \in L^{1}_{loc}\br{\R^{d}}$ and $x \in
\R^{d}$. Therefore, in order to prove the boundedness of the operator
$I^{\alpha,loc}_{\mu}$ from $L^{p}(w)$ to $L^{\nu}(w^{\nu/p})$ for
$w^{\nu/p} \in S_{1 + \frac{\nu}{p'},c}^{\mu}$ with $c > 0$ and $1 < p < \frac{d}{\alpha}$ it is
sufficient to prove the boundedness of $I^{\alpha,loc}_{0}$. This
follows on noting that $S_{1 + \frac{\nu}{p'},c}^{\mu} \subset A_{1 +
  \frac{\nu}{p'}}^{\mu,loc}$ for any $c > 0$ by Proposition \ref{prop:LocalApVc} and $I^{\alpha,loc}_{0}$ is bounded from
$L^{p}(w)$ to $L^{p}(w^{\nu/p})$ for $w^{\nu/p} \in A_{1 +
  \frac{\nu}{p'}}^{\mu,loc}$ and $1 < p < \frac{d}{\alpha}$ by
{\cite[Thm.~1]{bongioanni2011classes}}. It remains to prove the
boundedness of the global operators $I^{\alpha,glob}_{\mu}$.

\begin{thm} 
 \label{thm:RieszPotentials} 
 Fix $0 < \alpha \leq 2$. There exists a constant $c > 0$ for which
 $I^{\alpha}_{\mu,glob}$ is bounded from $L^{p}(w)$ to
 $L^{\nu}\br{w^{\nu/p}}$ for all weights $w$ with  $w^{\nu/p} \in S_{1 +
   \frac{\nu}{p'},c}^{V}$ and $1 < p < \frac{d}{\alpha}$, where
 $\frac{1}{\nu} = \frac{1}{p} - \frac{\alpha}{d}$. Moreover, the
 constant $c$ will be independent of $p$ and will depend on $\mu$ only
 through $C_{\mu}, \, D_{\mu}$ and $\delta_{\mu}$.
\end{thm}

\begin{prof}  
 Let $K_{\mu,\alpha}(x,y)$ denote the singular integral kernel of the operator
 $I^{\alpha}_{\mu}$. It will first be proved that $K_{\mu,\alpha}$
 satisfies the following  pointwise bound
 \begin{equation}
   \label{eqtn:FractionalPointwise}
K_{\mu,\alpha}(x,y) \lesssim \frac{e^{- \varepsilon d_{\mu}(x,y)}}{\abs{x -
    y}^{d - \alpha}} \quad \forall \ x, \, y \in \R^{d},
\end{equation}
for some $\varepsilon > 0$ that depends on $\mu$ only through
$C_{\mu}, \, D_{\mu}$ and $\delta_{\mu}$. The estimate
is obviously satisfied for $\alpha = 2$ owing to the pointwise bound
on the fundamental solution given by Theorem
\ref{thm:MeasureExp}. Consider the case $\alpha \in (0,2)$. The functional
 calculus (c.f. {\cite[pg.~286]{kato1980perturbation}}) implies that
 \begin{align*}\begin{split}  
 I^{\alpha}_{\mu}f(x) &= \frac{1}{\pi} \sin \br{\frac{\pi \alpha}{2}} \int^{\infty}_{0} \lambda^{-\alpha /2} \br{-
   \Delta + \mu + \lambda}^{-1} f(x) \, d \lambda \\
 &= \frac{1}{\pi} \sin \br{\frac{\pi \alpha}{2}} \int^{\infty}_{0} \lambda^{-\alpha / 2} \int_{\R^{d}} \Gamma_{\mu +
   \lambda}(x,y) f(y) \, dy \, d \lambda \\
 &= \int_{\R^{d}} \br{\frac{1}{\pi} \sin \br{\frac{\pi \alpha}{2}} \int^{\infty}_{0} \lambda^{-\alpha / 2} \Gamma_{\mu +
   \lambda}(x,y) \, d \lambda } f(y) \, dy \\
 &=: \int_{\R^{d}} K_{\mu,\alpha}(x,y) f(y) \, dy.
\end{split}\end{align*}
The inequality \eqref{eqtn:SumAgmon} together with the
pointwise bound on our fundamental solution implies 
\begin{align*}\begin{split}
 K_{\mu,\alpha}(x,y) &\lesssim e^{- \frac{\varepsilon'}{2} d_{\mu}(x,y)} \int^{\infty}_{0}
 \lambda^{- \alpha / 2} \frac{e^{- \frac{\varepsilon'}{2} \lambda^{\frac{1}{2}}
     \abs{x - y}}}{ \abs{x - y}^{d - 2}} \, d \lambda \\
 &\lesssim \frac{e^{- \frac{\varepsilon'}{2} d_{\mu}(x,y)}}{\abs{x - y}^{d - \alpha}}
\end{split}\end{align*}
for some $\varepsilon' > 0$, for all $x, \, y \in \R^{d}$. This
completes the proof of the  pointwise bound \eqref{eqtn:FractionalPointwise} for any $0 <
\alpha \leq 2$.

\vspace*{0.1in}

For $j \in
\N$, $k \geq 1$ and $x \in \R^{d}$, let the balls $B_{j}$, $B_{j,k}$,
$B_{\mu,x}$ and $B_{\mu,j}$
be as defined in the proof of the Theorem \ref{thm:RieszGlobal}. On expanding
out $\norm{I^{\alpha}_{\mu,glob}f}_{L^{\nu}(w^{\nu/p})}^{\nu}$,
\begin{align*}\begin{split}  
 \norm{I^{\alpha}_{\mu,glob}f}_{L^{\nu}(w^{\nu/p})}^{\nu} &\leq \sum_{j}
 \int_{B_{j}} \br{\int_{B_{\mu,x}^{c}} K_{\mu,\alpha}(x,y) \abs{f(y)}
   \, dy}^{\nu}
 w^{\nu / p}(x) \, dx \\
 &= \sum_{j}
 \int_{B_{j}} \br{\sum_{k = 1}^{\infty} \int_{(k + 1)B_{\mu,x} \setminus
   k B_{\mu,x}} K_{\mu,\alpha}(x,y) \abs{f(y)} \, dy}^{\nu}
 w^{\nu / p}(x) \, dx.
\end{split}\end{align*}
For $j \in \N$, $k \geq 1$ and $x \in B_{j}$, Lemma \ref{lem:Balls} gives
$$
\br{k B_{\mu,x}}^{c} \subset B(x, \frac{k}{\beta A_{0}}
\rho_{\mu}(x))^{c} \cup B(x,2 \rho_{\mu}(x))^{c}.
$$
Therefore, for $y \in (k + 1) B_{\mu,x} \setminus k B_{\mu,x}$,
$$
K_{\mu,\alpha}(x,y) \lesssim \frac{e^{- \varepsilon d_{\mu}(x,y)}}{\abs{x -
    y}^{d - \alpha}} \lesssim e^{- \delta k} \rho_{\mu}(x)^{-(d - \alpha)},
$$
where $\delta := \varepsilon A_{0}^{-1}$. Since $x \in B_{j}$
we then have by Lemma \ref{lem:CriticalMeasure},
\begin{equation}
  \label{eqtn:RieszPotKernel}
K_{\mu,\alpha}(x,y) \lesssim e^{- \delta k} \rho(x_{j})^{-(d - \alpha)}.
\end{equation}
Recall that it was proved in the proof of Theorem
\ref{thm:RieszGlobal} that the inclusion $(k + 1)B_{\mu,x} \subset (k
+ 1)B_{\mu,j}$ holds for any $j \in \N$, $x \in B_{j}$ and $k \geq 1$.
Applying the kernel estimate \eqref{eqtn:RieszPotKernel}, followed by
this inclusion and finally H\"{o}lder's inequality gives
\begin{align*}\begin{split}  
 & \norm{I^{\alpha}_{\mu,glob}f}_{L^{\nu}(w^{\nu/p})}^{\nu} \lesssim \sum_{j}
 \int_{B_{j}} \br{\sum_{k = 1}^{\infty} e^{- \delta k} \rho_{\mu}(x_{j})^{-(d - \alpha)} \int_{(k + 1)B_{\mu,x} \setminus
   k B_{\mu,x}} \abs{f(y)} \, dy}^{\nu}
w^{\nu / p}(x) \, dx \\
& \qquad \lesssim \sum_{j}
 \br{\sum_{k = 1}^{\infty} e^{- \delta k} \rho_{\mu}(x_{j})^{-(d -
     \alpha)} \norm{f}_{L^{p}((k + 1)B_{\mu,j},w)}
   w^{-\frac{1}{p - 1}}((k + 1)B_{\mu,j})^{\frac{p - 1}{p}}}^{\nu}
w^{\frac{\nu}{p}}(B_{j}).
 \end{split}\end{align*}
Since $B_{j} \subset B_{\mu,j}$ for each $j \in \N$, we then obtain
\begin{align}\begin{split}  
  \label{eqtn:RieszPotential1}
& \norm{I^{\alpha}_{\mu,glob}f}_{L^{\nu}(w^{\nu/p})}^{\nu} \\ & \quad  \lesssim \sum_{j}
 \br{\sum_{k = 1}^{\infty} e^{- \delta k} \rho_{\mu}(x_{j})^{-(d -
     \alpha)} \norm{f}_{L^{p}((k + 1)B_{\mu,j},w)}
   w^{-\frac{1}{p - 1}}((k + 1)B_{\mu,j})^{\frac{p - 1}{p}}
   w^{\frac{\nu}{p}}(B_{\mu,j})^{\frac{1}{\nu}}}^{\nu}.
 \end{split}\end{align}
 Using the assumption that $w^{\frac{\nu}{p}} \in S^{\mu}_{1 +
   \frac{\nu}{p'},c}$,
 \begin{align*}\begin{split}  
 w^{-\frac{1}{p-1}}\br{(k + 1)B_{\mu,j}}^{\frac{1}{p'}}
 w^{\frac{\nu}{p}}(B_{\mu,j})^{\frac{1}{\nu}} &\leq w^{-\frac{1}{p-1}}\br{(k + 1)B_{\mu,j}}^{\frac{1}{p'}}
 w^{\frac{\nu}{p}}((k + 1)B_{\mu,j})^{\frac{1}{\nu}} \\
 &\lesssim \abs{(k + 1)B_{\mu,j}}^{\br{1 + \frac{\nu}{p'}}
   \frac{1}{\nu}} e^{c \br{1 + \frac{\nu}{p'}}(k +
   1)\frac{1}{\nu}(\beta + A_{0}^{-1})}.
\end{split}\end{align*}
The estimate \eqref{eqtn:BallMeasure} then gives
\begin{align*}\begin{split}  
& w^{-\frac{1}{p-1}}\br{(k + 1)B_{\mu,j}}^{\frac{1}{p'}}
 w^{\frac{\nu}{p}}(B_{\mu,j})^{\frac{1}{\nu}} \\  & \qquad \qquad \lesssim (k + 1)^{d
   (k_{0} + 1) \br{1 + \frac{\nu}{p'}} \frac{1}{\nu}} \rho_{\mu}(x_{j})^{d
   \br{1 + \frac{\nu}{p'}}\frac{1}{\nu}} e^{c \br{1 +
     \frac{\nu}{p'}}(k + 1)\frac{1}{\nu}(\beta + A_{0}^{-1})} \\
 & \qquad \qquad \lesssim \rho_{\mu}(x_{j})^{d \br{1 + \frac{\nu}{p'}} \frac{1}{\nu}} e^{c'
   k} \\
 & \qquad \qquad = \rho_{\mu}(x_{j})^{d - \alpha} e^{c'
   k},
\end{split}\end{align*}
where $c' := 2 c \br{1 + \frac{\nu}{p'}} \frac{\beta}{\nu}$ and the
final line follows from the equality $d \br{1 +
  \frac{\nu}{p'}}\frac{1}{\nu} = d - \alpha$. Note that this equality
also implies that $c' = 2 c \beta \br{\frac{d - \alpha}{d}}$ and
therefore $c'$ is independent of $p$. Applying
this estimate to \eqref{eqtn:RieszPotential1} implies
\begin{align*}\begin{split}  
 \norm{I^{\alpha}_{\mu,glob}f}_{L^{\nu}(w^{\nu/p})} &\lesssim \br{\sum_{j}
 \br{\sum_{k = 1}^{\infty} e^{(c' - \delta)
       k}  \norm{f}_{L^{p}((k + 1)B_{\mu,j},w)}
   }^{\nu} }^{\frac{1}{\nu}} \\
 &\lesssim \sum_{k = 1}^{\infty} e^{(c' - \delta)k} \br{\sum_{j}
   \norm{f}_{L^{p}\br{(k + 1)B_{\mu,j},w}}^{\nu}}^{\frac{1}{\nu}} \\
 &\leq \sum_{k = 1}^{\infty} e^{(c' - \delta)k} \br{\sum_{j} \norm{f}^{\nu}_{L^{p}(B_{j,k},w)}}^{\frac{1}{\nu}},
\end{split}\end{align*}
where the inclusion $(k + 1)B_{\mu,j} \subset B_{j,k}$, as proved in
Theorem \ref{thm:RieszGlobal}, was used to obtain the final line.
The bounded overlap property of the balls $\lb B_{j} \rb_{j \in \N}$ and the fact that $\nu \geq p$ will then
complete our proof provided that we set $c$ small enough so that $c' <
\delta$.
\end{prof}

\subsection{The Heat Maximal Operator}
\label{subsec:Heat}

Let's now move onto the boundedness of the heat maximal operator for $L_{V}$ and
the proof of Theorem \ref{thm:Heat}. Let $k_{t}^{V} : \R^{d} \xx
\R^{d} \rightarrow \R$ denote the kernel of the operator $e^{- t
  \mathcal{L}_{V}}$ so that
$$
e^{-t \mathcal{L}_{V}}f(x) = \int_{\R^{d}} k_{t}(x,y) f(y) \, dy
$$
for all $f \in L^{1}_{loc}\br{\R^{d}}$ and $x \in \R^{d}$.
This function is called the heat kernel for the
operator $\mathcal{L}_{V}$.
We will require the following pointwise
estimate for the heat kernel proved by K. Kurata in
\cite{kurata2000estimate}.

\begin{prop}[{\cite[Thm.~1]{kurata2000estimate}}]
  \label{prop:KurataClassic}
  Suppose that $V \in RH_{\frac{d}{2}}$. There exist constants $D_{0}, \, D_{1}, \, D_{2} > 0$  such
that
$$
0 \leq k_{t}^{V}(x,y) \leq D_{0} \cdot e^{- D_{1} \br{1 +
    \frac{\sqrt{t}}{\rho_{V}(x)}}^{\frac{1}{(k_{0} + 1)}}}
\br{\frac{1}{t^{\frac{d}{2}}} e^{- D_{2} \frac{\abs{x - y}^{2}}{t}}}
$$
for all $x, \, y \in \R^{d}$ and $t > 0$, where $k_{0}$ is the
constant from Lemma \ref{lem:Shen0} corresponding to $\rho_{V}$.
\end{prop}

Notice that the previous proposition implies the pointwise estimate
$$
T^{*}_{V}f(x) \leq T^{*}_{0}f(x)
$$
for all $f \in L^{1}_{loc}\br{\R^{d}}$ and $x \in \R^{d}$. Therefore,
in order to prove the boundedness of the operator $T^{*}_{V,loc}$ on
$L^{p}(w)$ for $w \in H_{p,c}^{V,m}$ with $m, \, c > 0$ it is sufficient to prove the
boundedness of $T^{*}_{0,loc}$ on $L^{p}(w)$. This follows on noting
that $H_{p,c}^{V,m} \subset A_{p}^{V,loc}$ by Propositions
\ref{prop:LocalApVc} and \ref{prop:Inclusion} and that $T^{*,}_{0,loc}$ is bounded on $L^{p}(w)$
for weights in $A_{p}^{V,loc}$ by {\cite[Thm.~1]{bongioanni2011classes}}.
It remains to establish the boundedness of $T^{*}_{V,glob}$ on $L^{p}(w)$.

\begin{prop} 
 \label{prop:GlobalHeat} 
 There exists $c > 0$ such that the global part of the heat maximal operator
 $T^{*}_{V,glob}$ is bounded on $L^{p}(w)$ for any $w \in
 H_{p,c}^{V,m_{0}}$ where $1 < p < \infty$ and $m_{0} := (2(k_{0} + 1))^{-1}$.
\end{prop}

\begin{prof}  
  Let $c > 0$, $w \in H_{p,c}^{V,m_{0}}$ and fix $f \in L^{p}(w)$. Let's
 first obtain a pointwise estimate for $T^{*}_{V,glob}f$. Let $\lb
 B_{j} \rb_{j \in \N} = \lb B(x_{j}, \rho_{V}(x_{j})) \rb_{j \in \N}$
 be a cover of balls of $\R^{d}$ as given in Proposition \ref{prop:Cover}. Let $B_{x}$
 denote the critical ball $B_{x} := B(x,\rho_{V}(x))$. For $x \in
 B_{j}$ with $j \in \N$, Proposition \ref{prop:KurataClassic} implies that
 \begin{align*}\begin{split}  
 T^{*}_{V,glob}f(x) &\leq \sup_{t > 0} \int_{B_{x}^{c}} k_{t}^{V}(x,y)
 \abs{f(y)} \, dy \\
 &\lesssim \sup_{t > 0} \int_{B_{x}^{c}} \Phi_{2 m_{0}, D_{1}}^{V}(\sqrt{t},x)^{-1} \frac{e^{-
     D_{2} \frac{\abs{x - y}^{2}}{t}}}{t^{\frac{d}{2}}} \abs{f(y)} \, dy.
\end{split}\end{align*}
Notice that since $x \in B_{j}$ we have by Lemma \ref{lem:Shen0},
$$
\Phi_{2 m_{0}, D_{1}'}^{V}(\sqrt{t},x_{j}) \lesssim \Phi_{2 m_{0},D_{1}}^{V} (\sqrt{t},x)
$$
for some $D_{1}' > 0$, for all $t > 0$. This then gives
\begin{align*}\begin{split}  
 T^{*}_{V,glob}f(x) &\lesssim \sup_{t > 0} \Phi_{2 m_{0}, D_{1}'}^{V}(\sqrt{t},x_{j})^{-1}
 \sum_{k = 1}^{\infty}
 \int_{(k + 1)B_{x} \setminus k B_{x}} \frac{e^{- D_{2} \frac{\abs{x -
         y}^{2}}{t}}}{t^{\frac{d}{2}}} \abs{f(y)} \, dy \\
 &\lesssim \sup_{t > 0} \sum_{k = 1}^{\infty}
 \Phi_{2 m_{0}, D_{1}'}^{V}(\sqrt{t},x_{j})^{-1} \frac{e^{- D_{2}
     \frac{k^{2}\rho_{V}(x)^{2}}{t}}}{t^{\frac{d}{2}}} \int_{(k+1) B_{x}
   \setminus k B_{x}} \abs{f(y)} \, dy.
\end{split}\end{align*}
Since $x \in B_{j}$ we must have $\rho_{V}(x_{j}) \simeq
\rho_{V}(x)$ by Lemma \ref{lem:Shen0} and
therefore
$$
e^{-D_{2} \frac{k^{2} \rho_{V}(x)^{2}}{t}} \lesssim e^{- D_{2}'
  \frac{k^{2} \rho_{V}(x_{j})^{2}}{t}}
$$
for some $D_{2}' > 0$. On successively applying this, the inclusion
$B_{x} \subset \tilde{B}_{j} := 2 \sigma B_{j}$ where $\sigma := \beta
2^{\frac{k_{0}}{k_{0} + 1}}$ and H\"{o}lder's inequality,
\begin{align*}\begin{split}  
 T^{*}_{V,glob}f(x) &\lesssim \sup_{t > 0} \sum_{k = 1}^{\infty}
 \Phi_{2 m_{0},D_{1}'}^{V}(\sqrt{t},x_{j})^{-1} \frac{e^{- D_{2}'
     \frac{k^{2} \rho_{V}(x_{j})^{2}}{t}}}{t^{\frac{d}{2}}} \int_{(k + 1) B_{x}
   \setminus k B_{x}} \abs{f(y)} \, dy \\
 &\lesssim \sup_{t > 0} \sum_{k = 1}^{\infty}
 \Phi_{2 m_{0},D_{1}'}^{V}(\sqrt{t},x_{j})^{-1} \frac{e^{- D_{2}'
     \frac{k^{2} \rho_{V}(x_{j})^{2}}{t}}}{t^{\frac{d}{2}}} \int_{(k + 1) \tilde{B}_{j}} \abs{f(y)} \, dy \\
 &\lesssim \sup_{t > 0} \sum_{k = 1}^{\infty}
 \Phi_{2 m_{0},D_{1}'}^{V}(\sqrt{t},x_{j})^{-1} \frac{e^{- D_{2}'
     \frac{k^{2} \rho_{V}(x_{j})^{2}}{t}}}{t^{\frac{d}{2}}} \norm{f}_{L^{p}((k +
   1) \tilde{B}_{j},w)} w^{-\frac{1}{p-1}} \br{(k +
   1)\tilde{B}_{j}}^{\frac{p - 1}{p}}.
\end{split}\end{align*}
This pointwise estimate then allows us to estimate our norm by
\begin{align}\begin{split}  
    \label{eqtn:HeatNormExpansion}
    &\norm{T^{*}_{V,glob}f}^{p}_{L^{p}(w)}  \leq \sum_{j} \int_{B_{j}}
    T^{*}_{V,glob}f(x)^{p} w(x) \, dx  \\
    & \ \lesssim \sum_{j} \br{\sup_{t > 0} \sum_{k = 1}^{\infty}
 \Phi_{2 m_{0},D_{1}'}^{V}(\sqrt{t},x_{j})^{-1} \frac{e^{- D_{2}'
     \frac{k^{2} \rho_{V}(x_{j})^{2}}{t}}}{t^{\frac{d}{2}}} \norm{f}_{L^{p}((k +
   1) \tilde{B}_{j},w)} w^{-\frac{1}{p-1}} \br{(k +
   1)\tilde{B}_{j}}^{\frac{p - 1}{p}} w(B_{j})^{\frac{1}{p}}}^{p}. 
 \end{split}\end{align}
  The condition $w \in H_{p,c}^{V,m_{0}}$ implies that
  \begin{align*}\begin{split}  
 w (B_{j})^{\frac{1}{p}} w^{-\frac{1}{p - 1}} \br{(k +
   1)\tilde{B}_{j}}^{\frac{p - 1}{p}} &\leq w\br{(k +
   1)\tilde{B}_{j}}^{\frac{1}{p}}  w^{-\frac{1}{p - 1}} \br{(k +
   1)\tilde{B}_{j}}^{\frac{p - 1}{p}} \\
 &\lesssim e^{c (4\sigma)^{m_{0}}(k + 1)^{m_{0}}} (k +
1)^{d} \rho_{V}(x_{j})^{d} \\ &\lesssim e^{c (8\sigma)^{m_{0}} k^{m_{0}}} k^{d} \rho_{V}(x_{j})^{d}.
 \end{split}\end{align*}
Applying this estimate to \eqref{eqtn:HeatNormExpansion} then gives
\begin{align}\begin{split}  
    \label{eqtn:HeatNormExpansion2}
&\norm{T^{*}_{V,glob}f}^{p}_{L^{p}(w)} \\ & \ \lesssim \sum_{j} \br{\sup_{t > 0} \sum_{k = 1}^{\infty}
 \Phi_{2 m_{0},D_{1}'}^{V}(\sqrt{t},x_{j})^{-1} \frac{e^{- D_{2}'
     \frac{k^{2} \rho_{V}(x_{j})^{2}}{t}}}{t^{\frac{d}{2}}} \norm{f}_{L^{p}((k +
   1) \tilde{B}_{j},w)} e^{c (8\sigma)^{m_{0}} k^{m_{0}}} k^{d} \rho_{V}(x_{j})^{d}  }^{p}.    
 \end{split}\end{align}
Define, for $t > 0$, $j \in \N$ and $k \in \N^{*}$,
\begin{align*}\begin{split}  
 F(t,j,k) &:= \Phi_{2 m_{0},D_{1}'}^{V}(\sqrt{t},x_{j})^{-1}e^{- D_{2}'
  \frac{k^{2} \rho_{V}(x_{j})^{2}}{t}}e^{c (8\sigma)^{m_{0}} k^{m_{0}}} \br{\frac{k \rho_{V}(x_{j})}{t^{\frac{1}{2}}}}^{d} \\
&= \br{\frac{k \rho_{V}(x_{j})}{t^{\frac{1}{2}}}}^{d} \exp \br{ -
  D_{1}' \br{1 + \frac{\sqrt{t}}{\rho_{V}(x_{j})}}^{2 m_{0}} - D_{2}'
  \frac{k^{2}\rho_{V}(x_{j})^{2}}{t} + c (8\sigma)^{m_{0}} k^{m_{0}}}.
 \end{split}\end{align*}
It will now be proved that if $c > 0$ is small enough then there will
exist $\varepsilon > 0$ for which
\begin{equation}
  \label{eqtn:UniformEstimate}
F(t,j,k) \lesssim e^{-\varepsilon k^{m_{0}}}
\end{equation}
for all $t > 0$, $j \in \N$ and $k \in \N^{*}$.
First note
that from the estimate $x^{d} \lesssim e^{a x^{2}}$ for any $a > 0$
\begin{align*}\begin{split}  
 F(t,j,k) &\lesssim \exp \br{ -
  D_{1}' \br{1 + \frac{\sqrt{t}}{\rho_{V}(x_{j})}}^{2 m_{0}} - D_{2}'
  \frac{k^{2}\rho_{V}(x_{j})^{2}}{2 t} + c (8\sigma)^{m_{0}} k^{m_{0}}} \\
&\lesssim \exp \br{ -
  D_{1}' \br{\frac{t}{\rho_{V}(x_{j})^{2}}}^{m_{0}} - D_{2}'
  \frac{k^{2}\rho_{V}(x_{j})^{2}}{2 t} + c (8\sigma)^{m_{0}} k^{m_{0}}}.
\end{split}\end{align*}
Consider the case $t \geq k \rho_{V}(x_{j})^{2}$. For this case, we will
have the estimate
$$
F(t,j,k) \lesssim \exp \br{- D_{1}' k^{m_{0}} + c (8\sigma)^{m_{0}} k^{m_{0}}}.
$$
If we let $c$ be small enough, namely $c < \frac{D_{1}'}{(8\sigma)^{m_{0}}}$, then
we will have
\begin{equation}
  \label{eqtn:PartUnifEstimate}
F(t,j,k) \lesssim e^{- \varepsilon_{1} k^{m_{0}}}
\end{equation}
for some $\varepsilon_{1} > 0$.
Next, consider the case $t < k \rho_{V}(x_{j})^{2}$. In this situation, we
will have
$$
F(t,j,k) \lesssim \exp \br{- \frac{D_{2} ' k}{2} + c (8\sigma)^{m_{0}} k^{m_{0}}}.
$$
Since $m_{0} < 1$ we will then have
$$
F(t,j,k) \lesssim e^{- \varepsilon_{2} k}
$$
for some $\varepsilon_{2} > 0$. Putting this together with
\eqref{eqtn:PartUnifEstimate} then gives \eqref{eqtn:UniformEstimate}
for all $t > 0$, $j \in \N$ and $k \in \N^{*}$.

\vspace*{0.1in}

Referring back to \eqref{eqtn:HeatNormExpansion2}, we can apply
\eqref{eqtn:UniformEstimate} to obtain
\begin{align*}\begin{split}  
 \norm{T^{*}_{glob}f}_{L^{p}(w)} &\lesssim \br{\sum_{j} \br{ \sum_{k = 1}^{\infty}
 e^{- \varepsilon k^{m_{0}}} \norm{f}_{L^{p}((k +
   1) \tilde{B}_{j},w)} }^{p}}^{\frac{1}{p}} \\
&\lesssim \sum_{k = 1}^{\infty} e^{- \varepsilon k^{m_{0}}} \br{\sum_{j} \norm{f}_{L^{p}((k + 1)\tilde{B}_{j},w)}^{p}}^{\frac{1}{p}}.
\end{split}\end{align*}
It then follows from the bounded overlap property of the balls
$B_{j}$, Proposition \ref{prop:Cover}, that
\begin{align*}\begin{split}  
 \norm{T^{*}_{glob}f}_{L^{p}(w)} &\lesssim \br{\sum_{k = 1}^{\infty} k^{\frac{N_{1}}{p}}
   e^{- \varepsilon k^{m_{0}}}} \norm{f}_{L^{p}(\R^{d},w)} \\
 &\lesssim \norm{f}_{L^{p}(w)}.
 \end{split}\end{align*}

\end{prof}

\section{Uniformly Elliptic Operators with Potential}
\label{sec:UniformlyElliptic}

 Let $A$ be a $d
\xx d$ matrix-valued function with real-valued coefficients in
$L^{\infty}\br{\R^{d}}$. Suppose that $A$ satisfies the
ellipticity condition
\begin{equation}
  \label{eqtn:sec:Ellipticity}
\lambda \abs{\xi}^{2} \leq \langle A(x) \xi, \xi \rangle \leq
\Lambda \abs{\xi}^{2}
\end{equation}
for some $\lambda, \, \Lambda > 0$, for all $\xi \in \R^{d}$ and for almost every
$x \in \R^{d}$.
In this section we consider the uniformly elliptic operator with
potential $V \in RH_{\frac{d}{2}}$,
$$
L_{A,V} := - \mathrm{div} A \nabla + V.
$$
This operator is defined through its corresponding sesquilinear form as an unbounded operator on
$L^{2}\br{\R^{d}}$ with maximal domain.
Similar to the perturbation free case, appropriate analogues of the
usual operators can be defined,
$$
R_{A,V} := \nabla L_{A,V}^{-\frac{1}{2}}, \qquad R_{A,V}^{*} :=
L_{A^{*},V}^{-\frac{1}{2}} \nabla, \qquad I^{\alpha}_{A,V} := L_{A,V}^{-\frac{\alpha}{2}}
$$
for $0 < \alpha \leq 2$ and
$$
T^{*}_{A,V} f(x) := \sup_{t > 0} e^{- t L_{A,V}}\abs{f}(x) \qquad for \ f \in
L^{1}_{loc}\br{\R^{d}}, \ x \in \R^{d}.
$$
In this section we will prove weighted estimates for these operators
that are analogous to Theorems \ref{thm:Riesz} and
\ref{thm:Heat}. Before attempting to do so we must first discuss
exponential decay estimates for the associated kernels.

\subsection{Exponential Decay Estimates}

In a similar manner to the perturbation free case of Section
\ref{sec:Schrodinger}, the weighted results for $L_{A,V}$ will be
proved using exponential decay estimates for the relevant kernels.

Let $\Gamma_{A,V}$ denote the fundamental solution of the operator
$L_{A,V}$. This is a function defined on $\lb (x,y) \in \R^{d} \xx
\R^{d} : x \neq y \rb$ with the properties that $\Gamma_{A,V}(\cdot,y) \in
L^{1}_{loc}\br{\R^{d}}$ and $L_{A,V} \Gamma_{A,V}(\cdot,y) =
\delta_{y}$ for each $y \in \R^{d}$, where $\delta_{y}$ is the Dirac delta
distribution with pole at $y$. Refer to \cite{davey2018fundamental}
for the construction of this object. The following exponential decay estimates for
$\Gamma_{A,V}$ were proved by S. Mayboroda and B. Poggi in
\cite{mayboroda2019exponential}. These estimates are a generalisation
of the perturbation free case proved by Shen in
\cite{shen1999fundamental} (see Theorem \ref{thm:MeasureExp}).

\begin{thm}[{\cite[Cor.~6.16]{mayboroda2019exponential}}]
  \label{thm:MayborodaSharp}
 Let $A \in
 L^{\infty}\br{\R^{d};\mathcal{L}\br{\C^{d}}}$ have real-valued
 coefficients and assume that it
 satisfies the ellipticity condition
 \eqref{eqtn:sec:Ellipticity}. Fix $V \in RH_{\frac{d}{2}}$.
 There exist constants
 $C_{1}, \, C_{2}, \, \varepsilon_{1}, \, \varepsilon_{2} > 0$ for which
 $$
C_{1} \frac{e^{- \varepsilon_{1} d_{V}(x,y)}}{\abs{x - y}^{d - 2}} \leq
\Gamma_{A,V}(x,y) \leq C_{2} \frac{e^{- \varepsilon_{2} d_{V}(x,y)}}{\abs{x - y}^{d - 2}}
$$
for all $x, \, y \in \R^{d}$. The constants $C_{1}, \, C_{2}, \,
\varepsilon_{1}$ and $\varepsilon_{2}$ will depend on $V$ only through $\brs{V}_{RH_{\frac{d}{2}}}$.
\end{thm}

Using the previous theorem, it is then
possible to estimate the derivative of
$\Gamma_{A,V}$ from above provided that $A$ is H\"{o}lder continuous. 

\begin{prop}
  \label{prop:EllipticDerivative}
  Let $A \in
 L^{\infty}\br{\R^{d};\mathcal{L}\br{\C^{d}}}$ have real-valued
 coefficients and assume that it
 satisfies the ellipticity condition
 \eqref{eqtn:sec:Ellipticity}. Suppose also that $A$ is $\gamma$-H\"{o}lder continuous for some
 $\gamma \in (0,1)$.
  Fix $V \in RH_{q}$ for some $q > \frac{d}{2}$.

\begin{enumerate}[(i)]
\item  There exist
  constants $C_{1}, \, \varepsilon_{1} > 0$ for which
  $$
 \abs{\nabla \Gamma_{A,V}(x,y)} \leq C_{1} \frac{e^{- \varepsilon_{1}
    d_{V}(x,y)}}{\abs{x - y}^{d - 2}} \br{1 + \frac{1}{\abs{x - y}} +
  \int_{B(x,\abs{x - y}/2)} \frac{V(z)}{\abs{x - z}^{d - 1}} \, dz}
  $$
for all $x, \, y \in \R^{d}$.

\item Suppose that $q \geq d$. Then $C_{2}, \,
\varepsilon_{2} > 0$ can be chosen so that
$$
\abs{\nabla \Gamma_{A,V}(x,y)} \leq C_{2}
\frac{e^{-\varepsilon_{2}d_{V}(x,y)}}{\abs{x - y}^{d - 2}} \br{1 +
  \frac{1}{\abs{x - y}}}
$$
for all $x, \, y \in \R^{d}$.
\end{enumerate}
The constants $\varepsilon_{1}$ and $\varepsilon_{2}$ will
depend on $V$ only through $\brs{V}_{RH_{\frac{d}{2}}}$.
\end{prop}

\begin{prof}  
 The estimate for $q \geq d$ was proved in
 {\cite[Thm.~3.2]{bailey2019unbounded}} which itself was an adaptation
 of the perturbation free  arguments from {\cite[Lem.~2.20]{shen1999fundamental}}. Let's apply this argument to
 the case $q > \frac{d}{2}$. Fix $x, \, y \in \R^{d}$ with $x \neq y$ and define $R :=
 \abs{x - y} / 2$.  Let $r \leq R$. Define
 $$
u(\xi) := \Gamma_{A,V} (\xi,y)
$$
for $\xi \in B(x,R)$. Clearly $u$ is a weak solution to $L_{A,V} u = 0$
on the ball $B(x,R)$. Next define
$$
v(\xi) := u(\xi) + \int_{B(x,r)} \Gamma_{A,0}(\xi,z) u(z) V(z) \, d z
$$
for $\xi \in B(x,R)$. From reasoning identical to that of  {\cite[Thm.~3.2]{bailey2019unbounded}}, $v$ will be
$L_{A,0}$-harmonic in $B := B(x,r)$. Theorem 2.1 of
\cite{bailey2019unbounded} then implies that
$$
\norm{\nabla v}_{L^{\infty}(B/2)} \lesssim \frac{1}{r} \norm{v}_{L^{\infty}(B)}.
$$
Therefore,
\begin{align*}\begin{split}  
 \abs{\nabla \Gamma_{A,V}(x,y)} = \abs{\nabla u(x)} &\leq \abs{\nabla v(x)} + \int_{B} \abs{\nabla
   \Gamma_{A,0}(x,z)} \abs{u(z)} V(z) \, dz \\
 &\lesssim \frac{1}{r} \norm{v}_{L^{\infty}(B)} +
 \norm{u}_{L^{\infty}(B)} \int_{B} \abs{\nabla \Gamma_{A,0}(x,z)}
 V(z) \, dz \\
 &\lesssim \frac{1}{r} \norm{u}_{L^{\infty}(B)} \br{1 + \sup_{\xi \in B}
   \int_{B} \abs{\Gamma_{A,0}(\xi,z)} V(z) \, dz} \\ & \qquad \qquad  + \norm{u}_{L^{\infty}(B)} \int_{B} \abs{\nabla \Gamma_{A,0}(x,z)}
 V(z) \, dz .
\end{split}\end{align*}
For $\xi \in B$, by {\cite[Lem.~2.2]{bailey2019unbounded}}, Lemma \ref{lem:RHn2}
 and the inclusion $B \subset
B(\xi,2r) \subset B(x,4r)$,
\begin{align*}\begin{split}  
 \int_{B} \abs{\Gamma_{A,0}(\xi,z)} V(z) \, dz &\lesssim \int_{B}
 \frac{V(z)}{\abs{\xi - z}^{d - 2}} \, dz 
 \leq \int_{B(\xi,2r)} \frac{V(z)}{\abs{\xi - z}^{d - 2}} \, dz \\
 &\lesssim \frac{V(B(\xi,2r))}{r^{d - 2}}  \lesssim
 \frac{V(B(x,4r))}{r^{d - 2}} \\
 &\lesssim \frac{V(B(x,r))}{r^{d - 2}},
\end{split}\end{align*}
where the last line follows from the doubling property of $V$.
Lemma 2.3 of \cite{bailey2019unbounded} also gives
\begin{align*}\begin{split}  
 \int_{B} \abs{\nabla \Gamma_{A,0}(x,z)} V(z) \, dz &\lesssim \int_{B}
\frac{V(z)}{\abs{x - z}^{d - 1}} \, dz +\int_{B}
\frac{V(z)}{\abs{x - z}^{d - 2}} \, dz \\
&\lesssim \int_{B} \frac{V(z)}{\abs{x - z}^{d - 1}} \, dz +  \frac{V(B(x,r))}{r^{d - 2}}.
 \end{split}\end{align*}
Putting everything together implies that for any $r \leq \abs{x -
  y}/2$,
$$
  \abs{\nabla \Gamma_{A,V}(x,y)} \lesssim
  \norm{u}_{L^{\infty}(B(x,r))} \br{\frac{1}{r} +
    \frac{V(B(x,r))}{r^{d - 1}} + \frac{V(B(x,r))}{r^{d - 2}} +
    \int_{B(x,r)} \frac{V(z)}{\abs{z - x}^{d - 1}} \, dz}.
$$
Then, by setting $r = \abs{x - y}/2$ for the case $\abs{x - y} \leq 2
\rho_{V}(x)$ and $r = \rho_{V}(x)$ for the case $2 \rho_{V}(x) <
\abs{x - y}$ as in {\cite[Thm.~3.2]{bailey2019unbounded}}, this estimate will
imply part (i) of our proposition. 
 \end{prof}

Consider the fundamental solution of the potential free operator $- \mathrm{div} A
\nabla$, $\Gamma_{A,0}$.
An interesting consequence of the presence of the perturbation $A$ is
that, in contrast to the Laplacian case, the derivative of the
fundamental solution is no longer guaranteed to be bounded universally from above
by a constant multiple of $\abs{x - y}^{-(d - 1)}$. Instead, $\abs{\nabla \Gamma_{A,0}(x,y)}$ is only
guaranteed to satisfy this estimate locally. At a global scale, we can
only assert that
\begin{equation}
  \label{eqtn:PotentialFree}
  \abs{\nabla \Gamma_{A,0}(x,y)} \lesssim \frac{1}{\abs{x - y}^{d -
      2}} \qquad for \ all \ x, \, y \in \R^{d} \ s.t. \ \abs{x - y}
  \geq 1.
\end{equation}
See \cite{bailey2019unbounded} for a proof of this bound.
A corollary of this decrease in the strength of the global decay
is that the Riesz transform operator $R_{A,0} := \nabla (-
\mathrm{div} A \nabla)^{-\frac{1}{2}}$ will no longer necessarily be
Calder\'{o}n-Zygmund and therefore it will not necessarily be bounded
on $L^{p}(w)$ for all weights $w$ in the classical Muckenhoupt class
$A_{p}$ (see \cite{shen2005bounds}, Remarks $1.7$ and $1.8$). The perturbation $A$ and
potential $V$ therefore have two directly opposing effects on the
underlying Muckenhoupt class. The inclusion of the perturbation $A$
decreases the size of the associated weight class to weights that have
less global decay, while the potential
$V$ will increase the size of the class to include weights of greater
global decay. These two effects can be
played against each other and, for large enough $V$, the
effect of the perturbation $A$ can be effectively cancelled out by the
potential. This interaction betweeen the perturbation $A$ and the
potential $V$ is precisely what is responsible for the validity of the below
corollary. This corollary, and therefore the cancellation effect, will
play an important role when we come to consider the boundedness of the
Riesz transforms $R_{A,V}$.

\begin{cor}
  \label{cor:BoundedCriticalRadius}
Let $A \in
 L^{\infty}\br{\R^{d};\mathcal{L}\br{\C^{d}}}$ have real-valued
 coefficients and assume that it
 satisfies the ellipticity condition
 \eqref{eqtn:sec:Ellipticity}. Suppose also that $A$ is $\gamma$-H\"{o}lder continuous for some
 $\gamma \in (0,1)$.
 Fix $V \in RH_{q}$ for some $q > \frac{d}{2}$. Suppose
  that there exists $D_{V} > 0$ for which $\rho_{V}(x) \leq D_{V}$ for
  all $x \in \R^{d}$.

 \begin{enumerate}[(i)]
\item    There must exist constants $C_{1}, \, \varepsilon_{1} > 0$ for which
  $$
\abs{\nabla \Gamma_{A,V}(x,y)} \leq C_{1} \frac{e^{- \varepsilon_{1}
    d_{V}(x,y)}}{\abs{x - y}^{d - 2}} \br{\int_{B(x,\abs{x - y}/2)}
  \frac{V(z)}{\abs{z - x}^{d - 1}} \, dz + \frac{1}{\abs{x - y}}}
$$
for all $x, \, y \in \R^{d}$.

\item Assume that $q \geq d$.  Then
$C_{2}, \, \varepsilon_{2} > 0$ can be chosen so that 
$$
\abs{\nabla \Gamma_{A,V}(x,y)} \leq C_{2} \frac{e^{- \varepsilon_{2}
    d_{V}(x,y)}}{\abs{x - y}^{d - 1}}
$$
for all $x, \, y \in \R^{d}$.
\end{enumerate}

The constants $\varepsilon_{1}$ and $\varepsilon_{2}$ will depend on
$V$ only through $\brs{V}_{RH_{\frac{d}{2}}}$.
\end{cor}

\begin{prof}
  Let $x, \, y \in \R^{d}$. Suppose first that $\abs{x - y} \leq 2
  \rho_{V}(x)$. Then
  \begin{align*}\begin{split}  
 \frac{e^{- \varepsilon
     d_{V}(x,y)}}{\abs{x - y}^{d - 2}} 
 &= \abs{x - y} \frac{e^{- \varepsilon d_{V}(x,y)}}{\abs{x - y}^{d -
     1}} \\
 &\leq 2 \rho_{V}(x) \frac{e^{- \varepsilon d_{V}(x,y)}}{\abs{x - y}^{d -
     1}} \\
 &\lesssim_{D_{V}} \frac{e^{- \varepsilon d_{V}(x,y)}}{\abs{x - y}^{d - 1}}.
 \end{split}\end{align*}
  Next, suppose that $\abs{x - y} > 2 \rho_{V}(x)$. For this case,
  \begin{align*}\begin{split}  
  \frac{e^{- \varepsilon
     d_{V}(x,y)}}{\abs{x - y}^{d - 2}} 
 &= \rho_{V}(x) \frac{\abs{x - y}}{\rho_{V}(x)}  \frac{e^{- \varepsilon d_{V}(x,y)}}{\abs{x - y}^{d -
     1}} \\
 &\lesssim_{D_{V}} \frac{\abs{x - y}}{\rho_{V}(x)}  \frac{e^{- \varepsilon
     d_{V}(x,y)}}{\abs{x - y}^{d - 1}} \\
 &\lesssim d_{V}(x,y)^{k_{0} + 1} \frac{e^{- \varepsilon
     d_{V}(x,y)}}{\abs{x - y}^{d - 1}} \\
 &\lesssim \frac{e^{- \varepsilon' d_{V}(x,y)}}{\abs{x - y}^{d - 1}},
\end{split}\end{align*}
for any $\varepsilon' \in (0,\varepsilon)$, where we applied Lemma
\ref{lem:Shen} in the second to last last line. The result then
follows from Proposition \ref{prop:EllipticDerivative}.
\end{prof}

 The condition $\rho_{V} \leq D_{V}$ ensures that the size of the potential is
sufficiently large to cancel out the negative effect of the
perturbation $A$. Without this size condition the corollary will not
necessarily be valid and as a result the operator $R_{A,V}$ might not
be bounded on $L^{p}(w)$ for all $w \in A_{p}$.

Finally, the following exponential decay estimate for the heat kernel
$k_{t}^{A,V}$ of $L_{A,V}$ will be used in the proof of Theorem \ref{thm:EllipticHeat}.

\begin{thm}[{\cite[Thm.~1]{kurata2000estimate}}]
  \label{thm:Kurata}
Suppose that $V \in RH_{\frac{d}{2}}$. Let $A \in
L^{\infty}\br{\R^{d};\mathcal{L}\br{\C^{d}}}$ have real-valued
coefficients. Assume that $A$ satisfies the ellipticity condition
\eqref{eqtn:sec:Ellipticity} and $A = A^{*}$. There exist constants $D_{0}, \, D_{1}, \, D_{2} > 0$  such
that
$$
0 \leq k_{t}^{A,V}(x,y) \leq D_{0} \cdot e^{- D_{1} \br{1 +
    \frac{\sqrt{t}}{\rho_{V}(x)}}^{\frac{1}{(k_{0} + 1)}}}
\br{\frac{1}{t^{\frac{d}{2}}} e^{- D_{2} \frac{\abs{x - y}^{2}}{t}}}
$$
for all $x, \, y \in \R^{d}$ and $t > 0$, where $k_{0}$ is the
constant from Lemma \ref{lem:Shen} corresponding to $\rho_{V}$.
\end{thm}

\subsection{Riesz Transform}
\label{subsec:UniformlyEllipticRiesz}

For the Riesz transforms $R_{A,V}$ and their adjoints $R_{A,V}^{*}$,
our main result is as follows.

\begin{thm} 
 \label{thm:UniformlyElliptic} 
 Let $A \in
 L^{\infty}\br{\R^{d};\mathcal{L}\br{\C^{d}}}$ have real-valued
 coefficients and assume that it
 satisfies the ellipticity condition
 \eqref{eqtn:sec:Ellipticity}. Suppose also that $A$ is $\gamma$-H\"{o}lder continuous for some
 $\gamma \in (0,1)$. That is, there
exists some $\tau > 0$ so that for any $x, \, y \in \R^{d}$
$$
\abs{A(x) - A(y)} \leq \tau \abs{x - y}^{\gamma}.
$$
 Fix $V \in RH_{q}$ for some $q > \frac{d}{2}$. Assume that $\rho_{V}(x) \leq D_{V}$ for all $x \in \R^{d}$, for some $D_{V} > 0$.
The following
 statements are true.

\vspace*{0.1in}
 
 \begin{enumerate}[(i)]
 \item If $q \geq d$ then there exists $c_{1} > 0$ for which both $R_{A,V}$ and $R_{A,V}^{*}$
  are bounded on $L^{p}(w)$ for all $w \in S_{p,c_{1}}^{V}$ with $1 <
  p < \infty$.

  \vspace*{0.1in}

  \item Suppose instead that $\frac{d}{2} < q < d$ and let $s$ be
    defined through $\frac{1}{s} = \frac{1}{q} - \frac{1}{d}$. Then there
    exists a constant $c_{2} > 0$ for which the operator
    $R_{A,V}^{*}$ is bounded on $L^{p}(w)$ for $s' < p < \infty$ when $w
    \in S_{p/s',c_{2}}^{V}$ and the operator $R_{A,V}$ is bounded on $L^{p}(w)$
    for $1 < p < s$ when $w^{-\frac{1}{p - 1}} \in S_{p'/s',c_{2}}^{V}$.
\end{enumerate}
In each of the above statements, the constants $c_{1}$ and $c_{2}$
are independent of $p$ and depend
on $V$ only through $[V]_{RH_{\frac{d}{2}}}$ and $D_{V}$.
\end{thm}

\begin{rmk}
The condition $\rho_{V}(x) \leq D_{V}$ for all $x \in \R^{d}$ is a
size condition on the potential. From equation (0.13) of
\cite{shen1995lp}, we know that if $V = P$ is a non-negative
polynomial of degree $k \in \N$ then
$$
\rho_{V}(x)^{-1} = \sum_{\abs{\alpha} \leq k} \abs{\partial^{\alpha}_{k}
  P(x)}^{\frac{1}{\abs{\alpha} + 2}}.
$$
Therefore any non-negative polynomial will satisfy the size
condition $\rho_{V} \leq D_{V}$.
  \end{rmk}

 Throughout this
section fix $A \in
L^{\infty}\br{\R^{d};\mathcal{L}\br{\C^{d}}}$ with real-valued coefficients
 that satisfies \eqref{eqtn:sec:Ellipticity} and assume that $A$ is
 $\gamma$-H\"{o}lder continuous for some $\gamma \in (0,1)$.
The notation $K_{A,0}$ and $K_{A,V}$ will be used to denote the singular kernels
of the operators $R_{A,0}$ and $R_{A,V}$ respectively. Similarly, let
$K_{A,0}^{*}$ and $K_{A,V}^{*}$ be the kernels for $R_{A,0}^{*}$ and
$R_{A,V}^{*}$ respectively.
  The operator
  $L_{A,V}^{-\frac{1}{2}}$ can be expressed as
  $$
L_{A,V}^{-\frac{1}{2}}f(x) = \frac{1}{\pi} \int^{\infty}_{0}
\lambda^{-\frac{1}{2}} L_{A,V+\lambda}^{-1} f(x) \, d \lambda.
$$
Refer to {\cite[pg.~281]{kato1980perturbation}} for a proof of this
formula. Fubini's Theorem then implies
\begin{align*}\begin{split}  
 R_{A,V}f(x) &= \frac{1}{\pi} \int^{\infty}_{0} \lambda^{-\frac{1}{2}}
 \nabla L_{A,V + \lambda}^{-1}f(x) \, d \lambda \\
 &= \frac{1}{\pi} \int^{\infty}_{0} \lambda^{-\frac{1}{2}} \nabla
 \int_{\R^{d}} \Gamma_{A,V + \lambda}(x,y)f(y) \, dy \, d \lambda \\
 &= \int_{\R^{d}} \br{\frac{1}{\pi}\int^{\infty}_{0}
   \lambda^{-\frac{1}{2}} \nabla \Gamma_{A,V + \lambda}(x,y) \, d \lambda}
 f(y) \, dy.
\end{split}\end{align*}
This leads to the expression 
  \begin{equation}
    \label{eqtn:LocallyCZ1}
K_{A,V}(x,y) = \frac{1}{\pi} \int^{\infty}_{0} \lambda^{-\frac{1}{2}}
\nabla \Gamma_{A,V + \lambda}(x,y) \, d \lambda
\end{equation}
for all $x, \, y \in \R^{d}$.
It follows from duality that
$$
R_{A,V}^{*}f(x) = \int_{\R^{d}} K_{A,V}^{*}(x,y) f(y) \, dy =
\int_{\R^{d}} \overline{K_{A,V}(y,x)} f(y) \, dy.
$$
Then, since $K_{A,V}$ is real-valued,
$$
K_{A,V}^{*}(x,y) = \overline{K_{A,V}(y,x)} = K_{A,V}(y,x)
$$
for all $x, \, y \in \R^{d}$. Similarly we will have
$$
K_{A,0}^{*}(x,y) = K_{A,0}(y,x)
$$
for all $x, \, y \in \R^{d}$.

 As usual,
to prove the boundedness of the operators $R_{A,V}$ and $R_{A,V}^{*}$
it is sufficient to consider the local and global behaviour
separately. For the local behaviour, a few preliminary results must
first be proved.

\begin{prop}
  \label{prop:Regularity}
  There exists a constant $\varepsilon > 0$ such that for fixed $R > 0$ we
  have
  $$
\abs{\nabla \Gamma_{A,\lambda}(x,y) - \nabla \Gamma_{A,
    \lambda}(x',y)} \lesssim_{R}  \br{\frac{\abs{x -
      x'}}{\abs{x - y}}}^{\gamma} \frac{e^{- \varepsilon \sqrt{\lambda}
  \abs{x -  y}}}{\abs{x - y}^{d - 1}}
$$
for all $\lambda > 0$, $y\in \R^{d}$ and $x, \, x' \in B(y,R)$ with
$\abs{x - x'} \leq \frac{1}{2} \abs{x - y}$. 
  An
identical estimate will hold for  $\abs{\nabla \Gamma_{A, \lambda}(y,x) - \nabla
  \Gamma_{A,\lambda}(y,x')}$.
\end{prop}

\begin{prof}  
 Fix $R > 0$. In {\cite[Prop.~3.3]{bailey2019unbounded}}, it was
 proved that 
 \begin{equation}
   \label{eqtn:Regularity1}
    \abs{\nabla \Gamma_{A,\lambda}(x,y) - \nabla \Gamma_{A,
    \lambda}(x',y)} \lesssim_{R}  \frac{\abs{x -
      x'}^{\gamma}}{\abs{x - y}^{\gamma + 1}}
\sup_{\xi \in B(x,\frac{3}{4}\abs{x - y})}
\abs{\Gamma_{A,\lambda}(\xi,y)} \br{1 +
  \lambda \abs{x - y}^{2}}
   \end{equation}
for all $\lambda > 0$, $y\in \R^{d}$ and $x, \, x' \in B(y,R)$ with
$\abs{x - x'} \leq \frac{1}{2} \abs{x - y}$. Theorem
\ref{thm:MayborodaSharp} then implies that there exists $\varepsilon >
0$ such that for any $\xi \in
B(x,\frac{3}{4}\abs{x - y})$
\begin{align*}\begin{split}  
 \abs{\Gamma_{A,\lambda}(\xi,y)} &\lesssim \frac{e^{-\varepsilon \sqrt{\lambda}
     \abs{\xi - y}}}{\abs{\xi - y}^{d - 2}} \\
 &\lesssim \frac{e^{- \frac{\varepsilon}{4} \sqrt{\lambda} \abs{x -
       y}}}{\abs{x - y}^{d - 2}},
\end{split}\end{align*}
where we have used the identity $d_{\lambda}(\xi,y) =
\lambda^{\frac{1}{2}}\abs{\xi - y}$.
 This together with
\eqref{eqtn:Regularity1} produces our result.
\end{prof}

\begin{cor}
  \label{cor:LocallyCZ}
  The kernel $K_{A,0}$ is locally
  Calder\'{o}n-Zygmund. That is, for each $R > 0$
  $$
\abs{K_{A,0}(x,y)} \lesssim_{R} \abs{x - y}^{-d}
$$
for all $x, \, y \in \R^{d}$ with $\abs{x - y} \leq R$ and
$$
\abs{K_{A,0}(x,y) - K_{A,0}(x',y)} + \abs{K_{A,0}(y,x) -
  K_{A,0}(y,x')} \lesssim \frac{\abs{x - x'}^{\gamma}}{\abs{x - y}^{d + \gamma}} 
$$
for all $x, \, x', \, y \in \R^{d}$ with $x, \, x' \in B(y,R)$ and
$2 \abs{x - x'} \leq \abs{x - y}$.
\end{cor}

\begin{prof}
 For the size estimates, Proposition
\ref{prop:EllipticDerivative} combined with the expression
\eqref{eqtn:LocallyCZ1} lead to
\begin{align*}\begin{split}  
\abs{K_{A,0}(x,y)}  &\lesssim \int^{\infty}_{0} \lambda^{-\frac{1}{2}}
\abs{\nabla \Gamma_{A,\lambda}(x,y)} \, d \lambda \\
&\lesssim_{R} \int^{\infty}_{0} \lambda^{-\frac{1}{2}} \frac{e^{-
    \varepsilon \sqrt{\lambda} \abs{x - y}}}{\abs{x - y}^{d -
    1}} \, d \lambda \\
&\lesssim  \frac{1}{\abs{x - y}^{d}}.
 \end{split}\end{align*}
For the regularity estimate, the
expression \eqref{eqtn:LocallyCZ1} implies that
$$
 \abs{K_{A,0}(x,y) - K_{A,0}(x',y)} \lesssim
 \int^{\infty}_{0} \lambda^{-\frac{1}{2}} \abs{\nabla
   \Gamma_{A,\lambda}(x,y) - \nabla
   \Gamma_{A,\lambda}(x',y)} \, d \lambda.
$$
Taken with Proposition \ref{prop:Regularity}, this gives
\begin{align*}\begin{split}  
 \abs{K_{A,0}(x,y) - K_{A,0}(x',y)} &\lesssim
 \int^{\infty}_{0} \lambda^{-\frac{1}{2}} \abs{\nabla
   \Gamma_{A,\lambda}(x,y) - \nabla
   \Gamma_{A,\lambda}(x',y)} \, d \lambda \\
 &\lesssim_{R} \br{\frac{\abs{x - x'}}{\abs{x -
       y}}}^{\gamma} \frac{1}{\abs{x - y}^{d - 1}}
 \int^{\infty}_{0} \lambda^{-\frac{1}{2}} e^{- \varepsilon
   \sqrt{\lambda} \abs{x - y}} \, d \lambda \\
 &\lesssim  \br{\frac{\abs{x - x'}}{\abs{x -
       y}}}^{\gamma} \frac{1}{\abs{x - y}^{d}}.
\end{split}\end{align*}
The difference $\abs{K_{A,0}(y,x) - K_{A,0}(y,x')}$ is
estimated similarly.
\end{prof}

\begin{lem}
  \label{lem:Flatness}
Let $V \in RH_{q}$ for some $q > \frac{d}{2}$. For fixed $R > 0$, 
  \begin{equation}
    \label{eqtn:Flatness1}
\abs{K_{A,V}^{*}(x,y) - K_{A,0}^{*}(x,y)} \lesssim_{R} \frac{1}{\abs{x - y}^{d - 1}}
\int_{B(x,\abs{x - y}/2)} \frac{V(z)}{\abs{z - x}^{d - 1}} \, dz +
\br{\frac{\abs{x - y}}{\rho_{V}(x)}}^{2 - \frac{d}{q}} \frac{1}{\abs{x
  - y}^{d}}.
\end{equation}
for all $x, \, y \in \R^{d}$ with $\abs{x - y} \leq \mathrm{min} \br{\rho_{V}(x),R}$.
\end{lem}

\begin{prof}
  The perturbation free case of $- \Delta + V$ was shown in
  the proof of
  {\cite[Lem.~7.13]{shen1999fundamental}}. This argument will be generalised to our setting.
It is straightforward to prove using the uniqueness property of the
fundamental solution and the property $L_{A}^{V}
\Gamma_{A,V}(x,y) = \delta_{y}$ for all $y \in \R^{d}$, where
$\delta_{y}$ is the Dirac delta distribution, that
$$
\Gamma_{A,\lambda}(y,x) = \Gamma_{A,V + \lambda}(y,x) + \int_{\R^{d}}
\Gamma_{A,\lambda}(y,z) \Gamma_{A,V + \lambda}(z,x) V(z) \, dz
$$
for almost every $(x,y) \in \R^{d} \xx \R^{d}$, for any $\lambda > 0$.
Theorem \ref{thm:MayborodaSharp}, Proposition
\ref{prop:EllipticDerivative}.(ii) and \eqref{eqtn:SumAgmon} then imply
\begin{align}\begin{split}
    \label{eqtn:Flatness2}
 \abs{\nabla \Gamma_{A,V + \lambda}(y,x) - \nabla
   \Gamma_{A,\lambda}(y,x)} &\leq \int_{\R^{d}}
 \abs{\nabla \Gamma_{A^{*},\lambda}(z,y)} \Gamma_{A,V + \lambda}(z,x) V(z) \,
 dz \\
 &\lesssim \int_{\R^{d}} \frac{e^{-\varepsilon \sqrt{\lambda} \abs{y -
     z}}}{\abs{y - z}^{d - 2}} \br{1 + \frac{1}{\abs{y - z}}} \frac{e^{- \varepsilon \sqrt{\lambda}\abs{x - z}} e^{- \varepsilon
   d_{V}(z,x)}}{\abs{x - z}^{d - 2}} V(z) \, dz \\
&\lesssim_{R} \int_{B(y,R)} \frac{e^{-\varepsilon \sqrt{\lambda}
    \abs{y - z}}}{\abs{y - z}^{d - 1}}  \frac{e^{- \varepsilon \sqrt{\lambda}\abs{x - z}} e^{- \varepsilon
    d_{V}(z,x)}}{\abs{x - z}^{d - 2}} V(z) \, dz \\
& \qquad + \int_{\R^{d} \setminus B(y,R)} \frac{e^{- \varepsilon
    \sqrt{\lambda} \abs{y - z}}}{\abs{y - z}^{d - 2}}  \frac{e^{- \varepsilon \sqrt{\lambda}\abs{x - z}} e^{- \varepsilon
   d_{V}(z,x)}}{\abs{x - z}^{d - 2}} V(z) \, dz \\
 &:= J_{1} + J_{2}.
\end{split}\end{align}
The term $J_{1}$ is precisely the quantity obtained in
the proof of Lemma 7.13 of \cite{shen1999fundamental} and can
therefore be estimated from above by
\begin{equation}
  \label{eqtn:Flatness3}
J_{1} \lesssim e^{- \varepsilon \sqrt{\lambda}
  r} \br{\frac{1}{r^{d - 2}} \int_{B(y,r)} \frac{V(z) \, dz}{\abs{z -
      y}^{d - 1}} + \br{\frac{r}{\rho_{V}(x)}}^{2 - \frac{d}{q}}
  \frac{1}{r^{d - 1}}},
\end{equation}
where $r := \abs{x - y}/2$. As for the term $J_{2}$,  if $z \in \R^{d}
\setminus B(y,R)$ we will have $r := \abs{x - y}  \leq \abs{z - y}$. Therefore,
$$
J_{2} \lesssim e^{-
    \varepsilon \sqrt{\lambda} r} \int_{\R^{d} \setminus B(y,R)} \frac{1}{\abs{y - z}^{d - 2}}
\frac{e^{- \varepsilon \sqrt{\lambda}\abs{x - z}} e^{- \varepsilon
    d_{V}(z,x)}}{\abs{x - z}^{d - 2}} V(z) \, dz.
$$
Using the argument from Lemma 4.8 of
\cite{shen1999fundamental}, this term can be estimated by
$$
J_{2} \lesssim
\frac{e^{- \varepsilon \sqrt{\lambda} r}}{r^{d - 2}}
\br{\frac{r}{\rho_{V}(x)}}^{2 - \frac{d}{q}}
\lesssim_{R} \frac{e^{- \varepsilon \sqrt{\lambda} r}}{r^{d - 1}}
\br{\frac{r}{\rho_{V}(x)}}^{2 - \frac{d}{q}}.
$$
Combining this with \eqref{eqtn:Flatness3} produces
$$
\abs{\nabla \Gamma_{A,V + \lambda}(y,x) - \nabla \Gamma_{A,\lambda}(y,x)} \lesssim_{R}  e^{- \varepsilon \sqrt{\lambda}
  r} \br{\frac{1}{r^{d - 2}} \int_{B(y,r)} \frac{V(z) \, dz}{\abs{z -
      y}^{d - 1}} + \br{\frac{r}{\rho_{V}(x)}}^{2 - \frac{d}{q}}
  \frac{1}{r^{d - 1}}},
$$
for some $\varepsilon > 0$. This estimate can then be used together
with the expressions
$$
K_{A,V}^{*}(x,y) = \frac{1}{\pi} \int^{\infty}_{0}
\lambda^{-\frac{1}{2}} \nabla \Gamma_{A,V + \lambda}(y,x) \, d \lambda, \quad K_{A,0}^{*}(x,y) = \frac{1}{\pi} \int^{\infty}_{0}
\lambda^{-\frac{1}{2}} \nabla \Gamma_{A,\lambda}(y,x) \, d \lambda
$$
to obtain \eqref{eqtn:Flatness1}.
\end{prof}

Let $\lb B_{j} \rb_{j \in \N} = \lb B(x_{j}, \rho_{V}(x_{j})) \rb_{j
  \in \N}$ be a collection of balls as given in
Proposition \ref{prop:Cover}. Define $\sigma = \beta
2^{\frac{k_{0}}{k_{0} + 1}}$, with $k_{0}$ as in Lemma
\ref{lem:Shen0} and $\beta$ as in \eqref{eqtn:Beta}. Set $\tilde{B}_{j} := 2 \sigma B_{j}$ for each $j \in \N$.
 For each $j \in \N$, let $\chi_{j} : \R^{d} \rightarrow [0,1]$ be a smooth function that is
  identically equal to one on $\tilde{B}_{j}$, vanishes outside of $2
  \tilde{B}_{j}$ and $\norm{\nabla \chi_{j}}_{\infty} \lesssim 
  \frac{1}{\sigma \rho_{V}(x_{j})}$. For each $j \in \N$, define the operator
$$
R_{A,0}^{*,j}f(x) := \chi_{j}(x) \cdot R_{A,0}^{*} \br{f
  \chi_{j}}(x) \qquad for \ f \in L^{1}_{loc}\br{\R^{d}}, \, x \in \R^{d}.
$$

\begin{lem}
  \label{lem:CutOffRiesz}
Let $V \in RH_{q}$ for some $q > \frac{d}{2}$ and suppose that $\rho_{V}(x) \leq D_{V}$ for some $D_{V} > 0$, for all
 $x \in \R^{d}$. For $1 < p < \infty$ and $w \in A_{p}$, 
  $$
\norm{R_{A,0}^{*,j}f}_{L^{p}(w)} \lesssim  \brs{w}_{A_{p}}^{\max
  \br{1,1/(p - 1)}} \norm{f}_{L^{p}(w)}
$$
for all $f \in L^{p}(w)$ and $j \in \N$, where the implicit constant
is independent of $j \in \N$ and $w$.
\end{lem}

\begin{prof}
This will be proved by demonstrating that $R^{*,j}_{A,0}$ is a
Calder\'{o}n-Zygmund operator that is bounded on $L^{2}\br{\R^{d}}$
with constant independent of $j \in \N$ and whose kernel satisfies size and regularity
estimates with constants independent of $j$.
Corollary \ref{cor:LocallyCZ} states that the kernel of $R_{A,0}$ is locally
Calder\'{o}n-Zygmund. This implies that $K_{A,0}^{*}$ is itself locally Calder\'{o}n-Zygmund since 
$K_{A,0}^{*}(x,y) = K_{A,0}(y,x)$ is satisfied for all $x, \, y \in
\R^{d}$.  Set $R := 20 \sigma D_{V}$. Then
$$
\abs{K_{A,0}^{*}(x,y)} \lesssim_{R} \abs{x - y}^{-d}
$$
for all $x, \, y \in \R^{d}$ with $\abs{x - y} \leq R$ and
$$
\abs{K_{A,0}^{*}(x,y) - K_{A,0}^{*}(x',y)} \lesssim_{R} \frac{\abs{x -
    x'}^{\gamma}}{\abs{x - y}^{d + \gamma}}
$$
for all $x, \, x', \, y \in \R^{d}$ with $x, \, x' \in B(y,R)$ and
$\abs{x - x'} \leq \frac{1}{2}\abs{x - y}$.

Let's first prove the size
estimate for $K_{A,0}^{*,j}(x,y) := \chi_{j}(x) K_{A,0}^{*}(x,y) \chi_{j}(y)$. Fix $x, \, y \in \R^{d}$.
If either $x
\notin 2 \tilde{B}_{j}$ or $y \notin 2 \tilde{B}_{j}$ then $K_{A,0}^{*,j}(x,y)$
vanishes completely and the size estimate will be trivially satisfied.
Suppose that $x \in 2 \tilde{B}_{j}$ and $y \in 2 \tilde{B}_{j}$. It then follows
from the boundedness of the critical radius function that
\begin{align*}\begin{split}  
    \abs{x - y} &\leq \abs{x - x_{j}} + \abs{x_{j} - y} \\
    &\leq 4 \sigma \rho_{V}(x_{j}) + 4 \sigma \rho_{V}(x_{j}) \\
    &\leq 8 \sigma D_{V} \\
    &\leq R.
  \end{split}\end{align*}
Therefore,
$$
K_{A,0}^{*,j}(x,y) =  \chi_{\tilde{B}_{j}}(x) K_{A,0}^{*}(x,y)
\chi_{\tilde{B}_{j}}(y) \lesssim_{R} \abs{x - y}^{-d}.
$$
Next, let's prove the regularity estimate for $K_{A,0}^{*,j}$. Let
$x, \, x', \, y \in \R^{d}$ with $\abs{x - x'} \leq
\frac{1}{2}\abs{x - y}$. If either $y \notin 2 \tilde{B}_{j}$ or $x$ and
$x' \notin 2 \tilde{B}_{j}$ then the regularity estimate will be
trivially satisfied. It can therefore be assumed that $y \in 2
\tilde{B}_{j}$ and either $x$ or $x' \in 2 \tilde{B}_{j}$. This
will imply that both $x$ and $x'$ are contained in $8
\tilde{B}_{j}$ and therefore $x, \, x' \in B(y,R)$. We have
\begin{align*}\begin{split}  
 \abs{K_{A,0}^{*,j}(x,y) - K_{A,0}^{*,j}(x', y)} &\leq
 \abs{\chi_{j}(x) K_{A,0}^{*}(x,y) -
   \chi_{j}(x')K_{A,0}^{*}(x',y)} \\
 &\leq \abs{\chi_{j}(x) - \chi_{j}(x')}
 \abs{K_{A,0}^{*}(x,y)} + \abs{K_{A,0}^{*}(x,y) -
   K_{A,0}^{*}(x',y)} \\
 &\lesssim_{R} \norm{\nabla \chi_{j}}_{\infty} \abs{x - x'}
 \abs{K_{A,0}^{*}(x,y)} +  \frac{\abs{x -
     x'}^{\gamma}}{\abs{x - y}^{d + \gamma}} \\
 &\lesssim_{R} \frac{\abs{x - x'}}{\sigma
   \rho_{V}(x_{j})\abs{x - y}^{d}} +  \frac{\abs{x -
     x'}^{\gamma}}{\abs{x - y}^{d + \gamma}} \\
 &\lesssim \frac{\abs{x - x'}}{\sigma
   \abs{x - y}^{d + 1}} +  \frac{\abs{x -
     x'}^{\gamma}}{\abs{x - y}^{d + \gamma}} \\
 &\lesssim  \frac{\abs{x - x'}^{\gamma}}{\abs{x - y}^{d + \gamma}},
\end{split}\end{align*}
where the second to last line follows from the fact that $\abs{x - y}
\lesssim \rho_{V}(x_{j})$. Similar reasoning can be applied to obtain the
estimate
$$
\abs{K_{A,0}^{*,j}(y,x) - K_{A,0}^{*,j}(y,x')} \lesssim_{R}  \frac{\abs{x -
  x'}^{\gamma}}{\abs{x - y}^{d + \gamma}}
$$
for all $x, \, y \in \R^{d}$. Since it is obvious that the operators $R_{A,0}^{*,j}$ are
bounded on $L^{2}\br{\R^{d}}$ with constant independent of $j \in \N$,
it follows that the operators $R_{A,0}^{*,j}$ are
Calder\'{o}n-Zygmund with constants independent of $j \in \N$. Our result then follows from the well-known
$A_{2}$-conjecture that was proved in \cite{hytonen2012sharp}.
 \end{prof}

\begin{prop} 
 \label{prop:LocalEllipticRiesz} 
Let $V \in RH_{q}$ for some $q > \frac{d}{2}$ and assume there
exists $D_{V} > 0$ for which $\rho_{V}(x) \leq D_{V}$ for all $x \in \R^{d}$.
 \begin{enumerate}[(i)]
\item Suppose that $q \geq d$ and let $1 < p < \infty$. The operators
  $R_{A,V}^{loc}$ and $R_{A,V}^{*,loc}$ are bounded on $L^{p}(w)$ for
  any $w \in A_{p}^{V,loc}$.
  \item Suppose that $\frac{d}{2} < q < d$ and let $s$ be defined
    through $\frac{1}{s} = \frac{1}{q} - \frac{1}{d}$. The operator
    $R_{A,V}^{*,loc}$ is bounded on $L^{p}(w)$ for any $s' < p <
    \infty$ with $w \in A_{p/s'}^{V,loc}$ and the operator
    $R_{A,V}^{loc}$ is bounded on $L^{p}(w)$ for any $1 < p < s$ with
    $w^{-\frac{1}{p - 1}} \in A_{p'/s'}^{V,loc}$.
   \end{enumerate}
 \end{prop}

 \begin{prof}
For $1 < p < \infty$, weight $w$ on $\R^{d}$ and $f \in L^{p}(w)$, the triangle inequality allows us to estimate
the $L^{p}(w)$-norm of $R^{*,loc}_{A,V}f$ from above by
\begin{equation}
  \label{eqtn:LocalEllipticRiesz1}
  \norm{R^{*,loc}_{A,V}f}_{L^{p}(w)} \lesssim  \norm{\br{R_{A,V}^{*,loc} -
     R_{A,0}^{*,loc}}f}_{L^{p}(w)} + \norm{R_{A,0}^{*,loc}f}_{L^{p}(w)}.
  \end{equation}
In the proof of {\cite[Thm.~3]{bongioanni2011classes}}, the authors
prove the boundedness of the operator difference $R_{V,loc}^{*} -
R_{0,loc}^{*}$ using only the kernel estimate
$$
\abs{K_{V}^{*}(x,y) - K_{0}^{*}(x,y)} \lesssim \frac{1}{\abs{x - y}^{d - 1}}
\int_{B(x,\abs{x - y}/2)} \frac{V(z)}{\abs{z - x}^{d - 1}} \, dz +
\br{\frac{\abs{x - y}}{\rho_{V}(x)}}^{2 - \frac{d}{q}} \frac{1}{\abs{x
  - y}^{d}},
$$
for all $x, \, y \in \R^{d}$ with $\abs{x - y} \leq \rho_{V}(x)$.
By setting $R = D_{V}$ in Lemma \ref{lem:Flatness}, it is clear that
this estimate also holds for the difference $\abs{K_{A,V}^{*}(x,y) -
  K_{A,0}^{*}(x,y)}$. Therefore, an argument identical to that of 
{\cite[Thm.~3]{bongioanni2011classes}} can be used. This will imply
the boundedness of the difference $R_{A,V}^{*,loc} - R_{A,0}^{*,loc}$ on
$L^{p}(w)$ for any $w \in A_{p}^{V,loc}$ when $q \geq d$ and for $w
\in A_{p/s'}^{V,loc}$ with $s' < p < \infty$ when $\frac{d}{2} < q
< d$. 

\vspace*{0.1in}

To complete the proof of our proposition, it is then sufficient to
show that the operator $R_{A,0}^{*,loc}$ is bounded on
$L^{p}(w)$ for $1 < p < \infty$ and $w \in A_{p}^{V,loc}$ when $q >
\frac{d}{2}$. 
 Assume that $q > \frac{d}{2}$ and fix $1 < p < \infty$ and $w \in A_{p}^{V,loc}$. We have
\begin{align}\begin{split}  
    \label{eqtn:LocalEllipticRiesz2}
\norm{R_{A,0}^{*,loc}f}^{p}_{L^{p}(w)}    &\leq \sum_{j}
\int_{B_{j}} \abs{R_{A,0}^{*,loc}f}^{p}_{L^{p}(w)} w(x) \, dx \\
&\lesssim \sum_{j} \int_{B_{j}} \abs{R_{A,0}^{*,loc}f(x) -
  R_{A,0}^{*,j}f(x)}^{p} w(x) \, dx +
\int_{B_{j}} \abs{R_{A,0}^{*,j} f(x)}^{p}
w(x) \, dx.
\end{split}\end{align}
For any $x \in B_{j}$, it follows from Lemma \ref{lem:Shen0} that
$B(x,\rho_{V}(x)) \subset \tilde{B}_{j}$. The
fact that $K_{A,0}^{*}$ is locally Calder\'{o}n-Zygmund implies that for any
$x \in B_{j}$,
\begin{align*}\begin{split}  
 \abs{R^{*,loc}_{A,0}f(x) - R^{*,j}_{A,0}f(x)} &\leq \int_{2
   \tilde{B}_{j} \setminus B(x,\rho_{V}(x))} \abs{K_{A,0}^{*}(x,y)}
 \abs{f(y)} \, dy \\
 &\lesssim_{D_{V}} \int_{2 \tilde{B}_{j} \setminus B(x,\rho_{V}(x))}
 \frac{\abs{f(y)}}{\abs{x - y}^{d}} \, dy \\
 &\lesssim \frac{1}{\abs{2 \tilde{B}_{j}}} \int_{2 \tilde{B}_{j}}
 \abs{f(y)} \, dy,
\end{split}\end{align*}
where the last line follows from Lemma \ref{lem:Shen0}.
This will lead to
\begin{align*}\begin{split}  
 &\sum_{j} \int_{B_{j}} \abs{R^{*,loc}_{A,0}f(x) - R^{*,j}_{A,0} f(x)}^{p} w(x) \, dx \lesssim_{D_{V}} \sum_{j} \int_{B_{j}}
 \br{\frac{1}{\abs{2 \tilde{B}_{j}}} \int_{2 \tilde{B}_{j}} \abs{f(y)}
   \, dy}^{p} w(x) \, dx \\
 & \qquad \qquad \qquad \lesssim \sum_{j} w(2 \tilde{B}_{j}) w^{-\frac{1}{p - 1}} \br{2
   \tilde{B}_{j}}^{p - 1} \frac{1}{ \abs{2 \tilde{B}_{j}}^{p}}
 \int_{ 2 \tilde{B}_{j}} \abs{f(y)}^{p} w(y) \, dy.
\end{split}\end{align*}
On exploiting the property that $w \in A_{p}^{\rho_{V},loc} = A_{p}^{4
  \sigma \rho_{V},loc}$ (c.f. {\cite[Cor.~1]{bongioanni2011classes}}) and the bounded overlap property of the
balls $\lb B_{j} \rb_{j \in \N}$,
\begin{align*}\begin{split}  
 \sum_{j} \int_{B_{j}} \abs{R^{*,loc}_{A,0}f(x) - R^{*}_{A,0}(f
   \chi_{j})(x)}^{p} w(x) \, dx  &\lesssim \sum_{j}
 \int_{2 \tilde{B}_{j}} \abs{f(y)}^{p} w(y) \, dy \\
 &\lesssim \int_{\R^{d}} \abs{f(y)}^{p} w(y) \, dy.
\end{split}\end{align*}
For $j \in \N$, let $w_{j} \in A_{p}$ denote the extension of $w \vert_{2
  \tilde{B}_{j}}$ to all of $\R^{d}$ with $\brs{w_{j}}_{A_{p}} \leq \brs{w}_{A_{p}^{V,loc}}$. The existence of such a weight
is given in {\cite[Lem.~1]{bongioanni2011classes}}. 
For the second term in \eqref{eqtn:LocalEllipticRiesz2}, Lemma
\ref{lem:CutOffRiesz} implies
\begin{align*}\begin{split}  
 \sum_{j} \int_{B_{j}} \abs{R_{A,0}^{*,j}f(x)}^{p} w(x) \, dx &\leq \sum_{j}
 \int_{\R^{d}} \abs{R_{A,0}^{*,j}(f \mathbbm{1}_{2 \tilde{B}_{j}})(x)}^{p} w_{j}(x) \, dx \\
 &\lesssim \sum_{j} \brs{w_{j}}^{\mathrm{max}(1,1/(p -
   1))}_{A_{p}}\int_{2 \tilde{B}_{j}} \abs{f(x)}^{p} w_{j}(x) \, dx
 \\
 &\leq \sum_{j} \brs{w}^{\mathrm{max}(1,1/(p -
   1))}_{A_{p}^{V,loc}} \int_{2 \tilde{B}_{j}} \abs{f(x)}^{p} w(x) \, dx \\
 &\lesssim \int_{\R^{d}} \abs{f(x)}^{p} w(x) \, dx.
\end{split}\end{align*}
This completes our proof of the $L^{p}(w)$-boundedness of $R_{A,V}^{*,loc}$. The
boundedness of $R_{A,V}^{loc}$ follows from this by duality.
\end{prof}

It remains to consider the global behaviour of $R_{A,V}$ and $R_{A,V}^{*}$.

\begin{prop} 
  \label{prop:GlobalElliptic}
  Let $V \in RH_{q}$ for some $q > \frac{d}{2}$ and
 assume that $\rho_{V}(x) \leq D_{V}$ for all $x \in \R^{d}$, for some $D_{V} > 0$.
  \begin{enumerate}[(i)]
\item Suppose that $q \geq d$. There must exist a constant $c_{1} > 0$ for which the operators $R_{A,V}^{glob}$ and $R_{A,V}^{*,glob}$ are bounded on
  $L^{p}(w)$ for all $w \in S_{p,c_{1}}^{V}$ with $1 < p < \infty$.
  \item Suppose instead that $\frac{d}{2} < q < d$ and let $s$ be defined
    through $\frac{1}{s} = \frac{1}{q} - \frac{1}{d}$. There must
    exist $c_{2} > 0$ for which the operator
    $R_{A,V}^{*,glob}$ is bounded on $L^{p}(w)$ for any $s' < p <
    \infty$ with $w \in S_{p/s',c_{2}}^{V}$ and the operator
    $R_{A,V}^{glob}$ is bounded on $L^{p}(w)$ for any $1 < p < s$ with
    $w^{-\frac{1}{p - 1}} \in S_{p'/s',c_{2}}^{V}$.
    \end{enumerate}
 The
  constants $c_{1}$ and $c_{2}$ will be independent of $p$ and depend on $V$ only
  through $\brs{V}_{RH_{\frac{d}{2}}}$ and $D_{V}$.
 \end{prop}

 \begin{prof}
 From identical reasoning as in the proof of
   Lemma \ref{lem:RieszKernelEst}, Corollary \ref{cor:BoundedCriticalRadius} implies that the the singular kernel of $R_{A,V}^{*}$
   will satisfy
   $$
\abs{K^{*}_{A,V}(x,y)} \lesssim \frac{e^{-\varepsilon d_{V}(x,y)}}{\abs{x -
    y}^{d - 1}} \br{\int_{B(y,\abs{x - y}/2)} \frac{V(z)}{\abs{z -
      x}^{d - 1}} \, dz + \frac{1}{\abs{x - y}}}
$$
for all $x, \, y \in \R^{d}$, for some $\varepsilon > 0$. The proof of the perturbation free analogue of
the statement that we wish to prove, Theorem
   \ref{thm:RieszGlobal}, relied entirely on the estimate provided by
   Lemma \ref{lem:RieszKernelEst}. Since this estimate is also true in the perturbation
   dependent case, it follows that the proof from Theorem
   \ref{thm:RieszGlobal} can be repeated verbatim to give us our
   result.
   
 \end{prof}

 Combining Propositions \ref{prop:LocalEllipticRiesz} and \ref{prop:GlobalElliptic} completes the proof of Theorem \ref{thm:UniformlyElliptic}.

\subsection{Riesz Potentials}
\label{subsec:UniformlyEllipticRieszPot}

For the Riesz potentials $I^{\alpha}_{A,V}$ we have the following
theorem.

\begin{thm} 
  \label{thm:EllipticRieszPot}
  Fix $V \in RH_{\frac{d}{2}}$. Let $A \in
 L^{\infty}\br{\R^{d};\mathcal{L}\br{\C^{d}}}$ have real-valued
 coefficients and assume that it
 satisfies the ellipticity condition
 \eqref{eqtn:sec:Ellipticity}.
 For any $0 < \alpha \leq 2$, there
      must exist $c > 0$ for which the operator
      $I^{\alpha}_{A,V}$ is bounded from $L^{p}(w)$ to
$L^{\nu}(w^{\nu/p})$ for $w^{\nu/p} \in S_{1 + \frac{\nu}{p'},c}^{V}$
and $1 < p < \frac{d}{\alpha}$,
where $\frac{1}{\nu} = \frac{1}{p} - \frac{\alpha}{d}$. Moreover, the
constant $c$ is independent of $p$ and depends on $V$ only through $\brs{V}_{RH_{\frac{d}{2}}}$.
\end{thm}

\begin{prof}
 Let $1 < p < \frac{d}{\alpha}$ and $\frac{1}{\nu} =
\frac{1}{p} - \frac{\alpha}{d}$.
Let's
first prove the boundedness of the local component $I^{\alpha,loc}_{A,V}$. Notice that the pointwise estimate
$$
\abs{I_{A,V}^{\alpha}f(x)} \leq I_{0}^{\alpha}\abs{f}(x)
$$
holds for all $f \in L^{1}_{loc}\br{\R^{d}}$ and $x \in
\R^{d}$. Therefore, in order to prove the boundedness of the operator
$I_{A,V}^{\alpha,loc}$ from $L^{p}(w)$ to $L^{\nu}(w^{\nu/p})$ for
$w^{\nu/p} \in S_{1 + \frac{\nu}{p'},c}^{V}$ and $c > 0$ it is sufficient to prove the boundedness of
$I_{0}^{\alpha,loc}$. This follows on noting that $S^{V}_{1 +
  \frac{\nu}{p'},c} \subset A_{1 + \frac{\nu}{p'}}^{V,loc}$ for any $c
> 0$ by Proposition \ref{prop:LocalApVc} and
$I^{\alpha,loc}_{0}$ is bounded from $L^{p}(w)$ to $L^{p}(w^{\nu/p})$
for $w^{\nu/p} \in A_{1 + \frac{\nu}{p'}}^{V,loc}$ by 
{\cite[Thm.~1]{bongioanni2011classes}}.

It remains to consider the boundedness of the global component
$I^{\alpha,glob}_{A,V}$. The proof of Theorem
\ref{thm:RieszPotentials} relied entirely on the estimate
$$
\Gamma_{V}(x,y) \lesssim  \frac{e^{- \varepsilon d_{V}(x,y)}}{\abs{x -
  y}^{d - 2}}.
$$
Since an estimate of this form holds for the case $L_{A,V}$ by
Theorem \ref{thm:MayborodaSharp}, it follows
that the proof of Theorem \ref{thm:RieszPotentials} can be repeated
verbatim. This argument will show that there exists $c > 0$,
independent of $p$ and depending on $V$ only through $\brs{V}_{RH_{\frac{d}{2}}}$, so that $I^{\alpha,glob}_{A,V}$ is bounded from $L^{p}(w)$ to $L^{\nu}(w^{\nu/p})$ for any
$w^{\nu/p} \in S^{V}_{1 + \frac{\nu}{p'},c}$.
\end{prof}

 \subsection{Heat Maximal Operator}
 \label{subsec:UniformEllipHeat}

For the heat maximal operator $T^{*}_{A,V}$, our weighted result is
given below.

\begin{thm} 
  \label{thm:EllipticHeat}
  Fix $V \in RH_{\frac{d}{2}}$. Let $A \in
 L^{\infty}\br{\R^{d};\mathcal{L}\br{\C^{d}}}$ have real-valued
 coefficients and assume that it
 satisfies the ellipticity condition
 \eqref{eqtn:sec:Ellipticity}. Suppose also that
 $A = A^{*}$. Then there exists $c > 0$ for which the heat maximal
 operator $T^{*}_{A,V}$ is bounded on $L^{p}(w)$ for all $w \in
 H_{p,c}^{V,m_{0}}$ with $m_{0} := (2 (k_{0} + 1))^{-1}$ and $1 < p <
 \infty$. The constant $c$ is independent of $p$.
 \end{thm}

\begin{prof}  
 Let $1 < p < \infty$. Consider the heat maximal
 operator for $L_{A,V}$,
 $$
T^{*}_{A,V}f(x) := \sup_{t > 0} e^{- t L_{A,V}}\abs{f}(x) \qquad f \in
L^{1}_{loc}\br{\R^{d}}, \ x \in \R^{d}.
$$
From Theorem \ref{thm:Kurata}, we have the pointwise estimate
$$
T^{*}_{A,V} f(x) \leq T^{*}_{0}f(x),
$$
for any $f \in L^{1}_{loc}\br{\R^{d}}$ and $x \in \R^{d}$. Therefore,
in order to prove the boundedness of the operator $T^{*,loc}_{A,V}$ on
$L^{p}(w)$ for $w \in H_{p,c}^{V,m_{0}}$ and $c > 0$ it is sufficient to prove the boundedness of $T^{*,loc}_{0}$ on
$L^{p}(w)$. This follows on noting that $H_{p,c}^{V,m_{0}} \subset
A_{p}^{V,loc}$ for any $c > 0$ by Propositions \ref{prop:LocalApVc} and \ref{prop:Inclusion} and that $T^{*,loc}_{0}$ is bounded on
$L^{p}(w)$ for weights in $A_{p}^{V,loc}$ by {\cite[Thm~1]{bongioanni2011classes}}.

It remains to establish the boundedness of $T^{*,glob}_{A,V}$. The proof of Theorem \ref{thm:Heat} in Section
\ref{subsec:Heat} relied entirely on the pointwise estimate
 $$
k_{t}^{V}(x,y) \lesssim \frac{1}{t^{\frac{d}{2}}} e^{- D_{1} \br{1 +
    \frac{\sqrt{t}}{\rho_{V}(x)}}^{\frac{1}{k_{0} + 1}}}  e^{- D_{2}
  \frac{\abs{x - y}^{2}}{t}}.
$$
Since, by Theorem \ref{thm:Kurata}, this estimate also holds for the
heat kernel $k^{A,V}_{t}$, the entire proof 
of Proposition \ref{prop:GlobalHeat} can be repeated verbatim. This
will allow us to show that there exists $c > 0$, independent of $p$,
such that $T^{*,glob}_{A,V}$ is bounded on $L^{p}(w)$ for any $w \in H_{p,c}^{V,m_{0}}$.
 \end{prof}

\section{Magnetic Schr\"{o}dinger Operators}
\label{sec:Magnetic}

 The final form of
 Schr\"{o}dinger operator that will be considered in this
 article is the magnetic Schr\"{o}dinger operator. 
Let $a = (a_{1}, \cdots, a_{d})$ be a vector of
real-valued functions in $C^{1}\br{\R^{d}}$ that will be
referred to as the magnetic potential. The magnetic
field, denoted by $B : \R^{d} \rightarrow \R^{d}$, is then defined through
$$
B := \mathrm{curl} \, a.
$$
The magnetic Schr\"{o}dinger operator with electric potential $V$ and
magnetic potential $a$ is the operator 
$$
L_{V}^{a} := \br{\nabla - i a}^{*} \br{\nabla - i a} + V.
$$
 The standard operators associated with $L_{V}^{a}$ are defined through
$$
R_{V}^{a} := \br{\nabla - i a} \br{L_{V}^{a}}^{-\frac{1}{2}}, \qquad R_{V}^{a,*} :=
(L_{V}^{a})^{-\frac{1}{2}} \br{\nabla + i a}, \qquad I_{V}^{a,\alpha} := (L_{V}^{a})^{-\frac{\alpha}{2}}
$$
for $0 < \alpha \leq 2$ and
$$
T^{a,*}_{V}f(x) := \sup_{t > 0} e^{-t L_{V}^{a}}\abs{f}(x), \qquad f
\in L^{1}_{loc}\br{\R^{d}}, \, x \in \R^{d}.
$$
As usual, in order to prove
the weighted boundedness of our operators, we will require exponential
decay estimates for the associated fundamental solution,
denoted by $\Gamma^{a}_{V}$ (refer to \cite{davey2018fundamental} for the
construction of $\Gamma_{V}^{a}$), 
and heat kernel, denoted by $h_{t}^{a,V}$ (see \cite{kurata2000estimate}).

 \begin{thm}[{\cite[Cor.~6.16]{mayboroda2019exponential}}]
   \label{thm:MayborodaSharpMagnetic}
   Suppose that $V + \abs{B} \in RH_{\frac{d}{2}}$
and that there exists $D, \, D' > 0$ for which
\begin{equation}
  \label{eqtn:MagneticConds}
0 \leq V \leq D \cdot \rho_{V + \abs{B}}^{-2} \quad and \quad \abs{\nabla B}
\leq D' \cdot \rho_{V + \abs{B}}^{-3}.
\end{equation}
 There exist 
 constants $\varepsilon, \, C > 0$ for which
 $$
\abs{\Gamma_{V}^{a}(x,y)} \leq C \frac{e^{- \varepsilon
    d_{V + \abs{B}}(x,y)}}{\abs{x - y}^{d - 2}}
 $$
 for all $x, \, y \in \R^{d}$. The constant $\varepsilon$
 will depend on $V$ and $a$ only through $\brs{V + \abs{B}}_{RH_{\frac{d}{2}}}$, $D$ and $D'$.
\end{thm}

Let $\nabla^{a} := \nabla - i a$. Using Theorem \ref{thm:MayborodaSharpMagnetic} and the work of B. Ben Ali from \cite{ali2010maximal},
it is then not too difficult to prove exponential decay estimates for
the derivative of the fundamental solution of $L^{a}_{V}$.

\begin{prop} 
 \label{prop:MagneticDerivative} 
Suppose that $V + \abs{B} \in RH_{d}$
and that there exists $D, \, D' > 0$ for which
\eqref{eqtn:MagneticConds} is satisfied.
 There exist constants $\varepsilon, \, C > 0$ for which
 $$
\abs{\nabla^{a} \Gamma_{V}^{a}(x,y)} \leq C \frac{e^{- \varepsilon
    d_{V + \abs{B}}(x,y)}}{\abs{x - y}^{d - 1}}
$$
for all $x, \, y \in \R^{d}$. The constant $\varepsilon$ will
depend on $V$ and $a$ only through $\brs{V + \abs{B}}_{RH_{d}}$, $D$ and $D'$.
\end{prop}

\begin{prof}
Let $x, \, y \in \R^{d}$ and $0 < R \leq \abs{x - y} / 2 \sqrt{d}$. Consider the
cube $Q := Q(x,R)$ centered at $x$ with side-length $R$ and define the function $u$ on $Q$ through
$$
u(\xi) := \Gamma_{V}^{a}(\xi,y)
$$
for $\xi \in Q$. It is obvious that $u$ is a weak solution to $L^{a}_{V}
u = 0$ on $Q$.

\vspace*{0.1in}

\textit{Case 1:} Suppose that $\rho_{V + \abs{B}}(x) \leq \abs{x - y} / 2$. Set $R :=
\rho_{V + \abs{B}}(x) / \sqrt{d}$. {\cite[Lem.~4.10]{ali2010maximal}} and Theorem
\ref{thm:MayborodaSharpMagnetic} then imply that
\begin{align*}\begin{split}  
    \abs{\nabla^{a} \Gamma^{a}_{V}(x,y)} &\leq \sup_{\xi \in Q / 2} \abs{\nabla^{a} u(\xi)} \\
    &\lesssim \frac{1}{\rho_{V + \abs{B}}(x)} \sup_{\xi \in Q} \abs{u(\xi)} \\
    &\lesssim \frac{1}{\rho_{V + \abs{B}}(x)} \sup_{\xi \in  Q} \frac{e^{-
        \varepsilon d_{V + \abs{B}}(\xi,y)}}{\abs{\xi - y}^{d - 2}},
  \end{split}\end{align*}
where $Q/2 = Q(x,R / 2)$ is the cube of side-length $R/2$ centered at $x$.
For $\xi \in Q$ we have $\abs{x - y}
\lesssim \abs{\xi - y}$ and $d_{V + \abs{B}}(x,y) \leq d_{V + \abs{B}}(\xi,y) + C$
for some $C > 0$, the latter inequality following from the triangle
inequality and Lemma \ref{lem:LocalAgmon}. Lemma \ref{lem:Shen} then implies
\begin{align*}\begin{split}  
\abs{ \nabla^{a} \Gamma^{a}_{V}(x,y)} &\lesssim \frac{1}{\rho_{V + \abs{B}}(x)}
 \frac{e^{-\varepsilon d_{V + \abs{B}}(x,y)}}{\abs{x - y}^{d - 2}} \\
 &= \frac{\abs{x - y}}{\rho_{V + \abs{B}}(x)} \frac{e^{- \varepsilon
     d_{V + \abs{B}}(x,y)}}{\abs{x - y}^{d - 1}} \\
 &\lesssim d_{V + \abs{B}}(x,y)^{k_{0} + 1} \frac{e^{-\varepsilon d_{V
       + \abs{B}}(x,y)}}{\abs{x - y}^{d - 1}} \\
 &\lesssim \frac{e^{-\frac{\varepsilon}{2} d_{V + \abs{B}}(x,y)}}{\abs{x - y}^{d
   - 1}}.
 \end{split}\end{align*}

\vspace*{0.1in}

\textit{Case 2}: Suppose that $\rho_{V + \abs{B}}(x) > \abs{x - y} / 2$. Set $R :=
\abs{x - y} / 2 \sqrt{d}$. {\cite[Lem.~4.10]{ali2010maximal}} and Theorem
\ref{thm:MayborodaSharpMagnetic} then imply that
\begin{align*}\begin{split}  
   \abs{ \nabla^{a} \Gamma^{a}_{V}(x,y)} &\leq \sup_{\xi \in Q / 2} \abs{\nabla^{a} u(\xi)} \\
    &\lesssim \frac{1}{\abs{x - y}} \sup_{\xi \in Q} \abs{u(\xi)} \\
    &\lesssim \frac{1}{\abs{x - y}} \sup_{\xi \in Q} \frac{e^{-
        \varepsilon d_{V + \abs{B}}(\xi,y)}}{\abs{\xi - y}^{d - 2}}.
  \end{split}\end{align*}
The estimates $\abs{x - y} \lesssim \abs{\xi - y}$ and $d_{V +
  \abs{B}}(x,y) \leq d_{V + \abs{B}}(\xi,y) + C$ for some $C > 0$ and for
all $\xi \in Q$ then allow us to conclude the proof of our proposition.
\end{prof}

\begin{thm}[{\cite[Thm.~1]{kurata2000estimate}}]
  \label{thm:Kurata2}
 Suppose that $V + \abs{B} \in RH_{\frac{d}{2}}$ and
 \eqref{eqtn:MagneticConds} is satisfied with constants $D, \, D' >
 0$. There exist
constants $D_{0}, \, D_{1}, \, D_{2} > 0$ such that
$$
\abs{h_{t}^{a,V}(x,y)} \leq D_{0} \cdot e^{- D_{1} \br{1 +
    \frac{\sqrt{t}}{\rho_{V + \abs{B}(x)}}}^{\frac{1}{k_{0} + 1}}}
\br{\frac{1}{t^{\frac{d}{2}}} e^{- D_{2} \frac{\abs{x - y}^{2}}{t}}}
$$
for all $x, \, y \in \R^{d}$ and $t > 0$.
\end{thm}

For the global components of the Riesz transforms $R^{a}_{V}$ and
their adjoints $R^{a,*}_{V}$ our weighted result is as given below.

\begin{thm} 
 \label{thm:GlobalMagneticRiesz} 
 Let $1 < p < \infty$.  Suppose that $V + \abs{B} \in RH_{d}$ and
 \eqref{eqtn:MagneticConds} is satisfied with constants $D, \, D' >
 0$. Then there exists $c > 0$ for which
$R_{V,glob}^{a}$ and $R_{V,glob}^{a,*}$
are bounded on $L^{p}(w)$ for any $w \in S_{p,c}^{V + \abs{B}}
$. Moreover, $c$ will depend on $V$ and $a$ only through $\brs{V +
  \abs{B}}_{RH_{\frac{d}{2}}}$, $D$ and $D'$.
\end{thm}

\begin{prof}  
. The proof of Theorem \ref{thm:RieszGlobal} is entirely reliant
 on the estimate
 $$
\abs{\nabla \Gamma_{V}(x,y)} \lesssim \frac{e^{- \varepsilon
    d_{V}(x,y)}}{\abs{x - y}^{d - 1}}.
$$
Since the appropriate magnetic analogue for this estimate is true by
Proposition \ref{prop:MagneticDerivative}, the proof of Theorem
\ref{thm:RieszGlobal} can be reapplied verbatim to this case. This will be
enough to tell us that there exists $c > 0$ such that $R_{V,glob}^{a,*}$ is
bounded on $L^{p}(w)$ for all $w \in S_{p,c}^{V + \abs{B}}$.
 The boundedness of $R_{V,glob}^{a}$ follows by duality. 
 \end{prof}

In order to prove that the operators
$R_{V}^{a}$ and $R_{V}^{a,*}$ are bounded on $L^{p}(w)$ for weights in
the class $S_{p,c}^{V + \abs{B}}$ it is then sufficient to prove the
boundedness of the local components. A condition that can be used to
guarantee this is the rather strong
hypothesis that $R_{V}^{a}$ and $R_{V}^{a,*}$ are both bounded on
$L^{p}(w)$ for any weight $w$ in the standard Muckenhoupt class
$A_{p}$. Investigating when exactly this occurs is unfortunately
beyond the scope of this article.

\begin{thm}
\label{thm:MagneticRiesz}
Let $1 < p < \infty$.  Suppose that $V + \abs{B} \in RH_{d}$ and
 \eqref{eqtn:MagneticConds} is satisfied with constants $D, \, D' >
 0$. Suppose also that for any $w \in A_{p}$
$$
\norm{R_{V}^{a}}_{L^{p}(w)}, \, \norm{R^{a,*}_{V}}_{L^{p}(w)} \lesssim \brs{w}_{A_{p}}^{l}
$$
for some $l \geq 1$. Then there exists $c > 0$ for which
$R_{V}^{a}$ and $R_{V}^{a,*}$
are bounded on $L^{p}(w)$ for any $w \in S_{p,c}^{V + \abs{B}}
$. Moreover, $c$ will depend on $V$ and $a$ only through $\brs{V +
  \abs{B}}_{RH_{\frac{d}{2}}}$, $D$ and $D'$.
\end{thm}

\begin{prof}  
Let $c > 0$ be as given in Theorem \ref{thm:GlobalMagneticRiesz} and fix $w \in S_{p,c}^{V + \abs{B}}$. For the local
boundedness of $R_{V}^{a,*}$, let $B_{j} = B (x_{j}, \rho_{V +
  \abs{B}}(x_{j}))$ for $j \in \N$ be as given in Proposition
\ref{prop:Cover} and set $\tilde{B}_{j} := 2 \sigma B_{j}$.  We have
\begin{align*}\begin{split}  
 \norm{R^{a,*}_{V,loc}f}^{p}_{L^{p}(w)} &\leq \sum_{j}
 \int_{B_{j}} \abs{R^{a,*}_{V,loc}f(x)}^{p}_{L^{p}(w)} w(x) \, dx \\
 &\lesssim \sum_{j} \int_{B_{j}} \abs{R^{a,*}_{V,loc}f(x) -
   R^{a,*}_{V}(f \cdot \mathbbm{1}_{2 \tilde{B}_{j}})}^{p} w(x) \, dx \\ &
 \qquad \qquad + \int_{B_{j}}
 \abs{R^{a,*}_{V} (f \cdot \mathbbm{1}_{2 \tilde{B}_{j}})}^{p} w(x) \, dx.
\end{split}\end{align*}
The first term can be handled in an identical manner to how the
difference term from \eqref{eqtn:LocalEllipticRiesz2} was
handled by making use of the derivative estimates, Proposition \ref{prop:MagneticDerivative}. The second term can also be handled in a similar manner to
how the second term was handled in
\eqref{eqtn:LocalEllipticRiesz2}. Namely, let $w_{j} \in A_{p}$ denote
the extension of $w \vert_{2 \tilde{B}_{j}}$ to all of $\R^{d}$ with
$\brs{w_{j}}_{A_{p}} \leq \brs{w}_{A_{p}^{V + \abs{B},loc}}$ (c.f. {\cite[Lem.~1]{bongioanni2011classes}}). Then on
applying our hypothesis,
\begin{align*}\begin{split}  
 \sum_{j} \int_{B_{j}} \abs{R^{a,*}_{V}(f \mathbbm{1}_{2
     \tilde{B}_{j}})(x)}^{p} w(x) \, dx &\leq \sum_{j} \int_{\R^{d}}
 \abs{R^{a,*}_{V}(f \mathbbm{1}_{2 \tilde{B}_{j}})(x)}^{p} w_{j}(x) \,
 dx \\
 &\lesssim \sum_{j} \brs{w_{j}}^{l}_{A_{p}} \int_{2 \tilde{B}_{j}}
 \abs{f(x)}^{p} w_{j}(x) \, dx \\
 &\leq \sum_{j} \brs{w}_{A_{p}^{V + \abs{B},loc}}^{l} \int_{2
   \tilde{B}_{j}} \abs{f(x)}^{p} w(x) \, dx \\
 &\lesssim \int_{\R^{d}} \abs{f(x)}^{p} w(x) \, dx.
\end{split}\end{align*}
Duality will imply that $R_{V,loc}^{a}$ is also bounded on $L^{p}(w)$.
 \end{prof}

Our weighted result for the magnetic Riesz potentials and heat maximal
operator is as follows.

\begin{thm} 
 \label{thm:MagneticHeatRiesz} 
Let $1 < p < \infty$. Suppose that $V + \abs{B} \in RH_{\frac{d}{2}}$
and that there exists $D, \, D' > 0$ for which \eqref{eqtn:MagneticConds} is satisfied.
The following statements are true.

\begin{enumerate}[(i)]
\item For $0 < \alpha \leq 2$, there
      must exist a constant $c_{1} > 0$ for which the operator
      $I^{a,\alpha}_{V}$ is bounded from $L^{p}(w)$ to
$L^{\nu}(w^{\nu/p})$ for $w^{\nu/p} \in S_{1 + \frac{\nu}{p'},c_{1}}^{V}$
and $1 < p < \frac{d}{\alpha}$,
where $\frac{1}{\nu} = \frac{1}{p} - \frac{\alpha}{d}$.
\item There exists $c_{2} > 0$ for which $T^{a,*}_{V}$ is
 bounded on $L^{p}(w)$ for all $w \in H_{p,c_{2}}^{V + \abs{B},m_{0}}$ with
 $m_{0} := (2(k_{0} + 1))^{-1}$.
\end{enumerate}
\end{thm}

\begin{prof}  
The weighted boundedness of the operators $I^{a,\alpha}_{V}$ and $T^{a,*}_{V}$ is proved in an identical manner to Theorem
\ref{thm:EllipticRieszPot} and \ref{thm:EllipticHeat} respectively by
making use of Theorems \ref{thm:MayborodaSharpMagnetic} and \ref{thm:Kurata2}. 
 \end{prof}

\section{Necessary Conditions}
\label{sec:Optimality}

Let us now investigate necessary conditions for a weight $w$ to satisfy in
order for $R_{V}$ and $T^{*}_{V}$ to be bounded on $L^{p}(w)$.
We will begin by proving that $w \in A_{p}^{V,loc}$ if
$\norm{R_{V,loc}}_{L^{p}(w)} < \infty$ for $V \in RH_{d}$, thereby
demonstrating that at a local scale the operator $R_{V}$ characterises $A_{p}^{V,loc}$.
Following this, it will be shown that the first chain of inclusions in  Conjecture
\ref{conj:Upper} is true for constant potentials and the second chain
of inclusions is true for any potential that is bounded both from
above and below. Finally, we will prove that the second chain of
inclusions in Conjecture \ref{conj:Weak} is true for the harmonic
oscillator potential $V(x) = \abs{x}^{2}$.

\subsection{Characterisation of $A_{p}^{V,loc}$ in terms of $R_{V,loc}$}

In this section it will be proved that the condition $w \in
A_{p}^{V,loc}$ is necessary in order for the operator $R_{V,loc}$ to
be bounded on
$L^{p}(w)$. To keep the result as general as possible, we will
prove this statement for generalised Schr\"{o}dinger operators with measure potentials $\mu$ in the sense of Section
\ref{subsec:SchrodingerMeasure}. When taken together with Proposition \ref{prop:MeasureLocalized}, 
this result will complete the proof of a localized
Hunt-Muckenhoupt-Wheeden type
theorem for measure potentials $\mu$ satisfying
\eqref{eqtn:Measure1} and \eqref{eqtn:Measure2} with constants
$C_{\mu}, \, D_{\mu}$ and $\delta > 1$. In particular, the below theorem proves that the class
$A_{p}^{\mu,loc}$ is characterised completely by the boundedness of the
operator $R_{\mu,loc}$ on $L^{p}(w)$.

\begin{thm} 
  \label{thm:LocalHMW}
  Let $\mu$ be a non-negative Radon measure on $\R^{d}$ that satisfies
  \eqref{eqtn:Measure1} and \eqref{eqtn:Measure2} with constants
  $C_{\mu}, \, D_{\mu} > 0$ and $\delta_{\mu} > 1$. Fix $1 < p < \infty$ and
  weight $w$ on $\R^{d}$.
  Suppose that the operator $R_{\mu,loc}$ is bounded on
 $L^{p}(w)$. Then it must be true that $w \in A_{p}^{\mu,loc}$.
\end{thm}

\begin{prof}  
Assume that
$\norm{R_{\mu,loc}}_{L^{p}(w)} < \infty$. Let $0 < \kappa_{0} < 1$ be some constant whose value will be determined
 at a later time. Since $\rho_{\mu}$ is a critical radius function in
 the sense of Definition \ref{def:CriticalRadius}, it follows easily
 that the function $\kappa_{0} \cdot \rho_{\mu}$ will also satisfy
 \eqref{eqtn:Shen0} and therefore will also be a critical radius
 function. Moreover, it follows from {\cite[Cor.~1]{bongioanni2011classes}}
 that $A_{p}^{\kappa_{0} \cdot \rho_{\mu},loc} =
 A_{p}^{\rho_{\mu},loc}$. Therefore, in order to prove that $w \in
 A_{p}^{\rho_{\mu},loc}$ it is sufficient to prove that $w \in
 A_{p}^{\kappa_{0}\cdot \rho_{\mu},loc}$. It must therefore be proved that
\begin{equation}
  \label{eqtn:LocalHMW12}
w(B)^{\frac{1}{p}} w^{-\frac{1}{p - 1}}(B)^{\frac{p - 1}{p}} \lesssim \abs{B}
\end{equation}
for all Euclidean balls $B = B(c,r) \subset \R^{d}$ with $r \leq
\kappa_{0} \cdot \rho_{\mu}(c)$.

 Consider the operator
$$
R^{(1)}_{\mu}f(x) := \partial_{1} (-\Delta + \mu)^{-\frac{1}{2}}f(x).
$$
On combining the expression
$$
 \Gamma_{\mu + \lambda}(x,y) =  \Gamma_{\lambda}(x,y) -
\int_{\R^{d}}  \Gamma_{\lambda}(x,z) \Gamma_{\mu + \lambda}(z,y) \, 
d \mu(z) \quad \forall \ x, \, y \in \R^{d}
$$
with \eqref{eqtn:KernelExpression} we obtain the identity
\begin{equation}
  \label{eqtn:LocalHMW1}
K_{\mu}^{(1)}(x,y) = K_{0}^{(1)}(x,y) -  \frac{1}{\pi} \int^{\infty}_{0}
\lambda^{-\frac{1}{2}} \int_{\R^{d}} \partial_{1} \Gamma_{\lambda}(x,z) \Gamma_{\mu + \lambda}(z,y)
\, d \mu(z) \, d \lambda
\end{equation}
for all $x, \, y \in \R^{d}$, where $K_{\mu}^{(1)}$ and $K_{0}^{(1)}$ are the first components of the
singular kernels for $R_{\mu}$ and $R_{0}$ respectively.

Fix $B = B(c,r) \subset \R^{d}$ with
$r \leq \kappa_{0} \cdot c$. Let $B'$ be the
Euclidean ball in $\R^{d}$ that has center $c' = c + 4 r = (c_{1}
+ 4 r, \cdots, c_{d} + 4 r)$ and radius $r$. Let's
estimate the size of the second term in
\eqref{eqtn:LocalHMW1} for $x \in B'$ and $y \in B$. Theorem \ref{thm:MeasureExp} states that there
must exist some $\varepsilon > 0$ such that
\begin{align}\begin{split}
    \label{eqtn:LocalHMW11}
& \int_{\R^{d}} \abs{\partial_{1} \Gamma_{\lambda}(x,z)} \Gamma_{\mu +
   \lambda}(z,y) \, d \mu(z)  \lesssim \int_{\R^{d}} \frac{e^{-
     \varepsilon \sqrt{\lambda} \abs{x - z}}}{\abs{x - z}^{d - 1}} \frac{e^{-
     \varepsilon \sqrt{\lambda} \abs{y - z}} e^{-\varepsilon
     d_{\mu}(z,y)}}{\abs{z - y}^{d - 2}} \, d \mu(z) \\
 & \qquad \qquad \lesssim e^{- \frac{\varepsilon}{2} \sqrt{\lambda} \abs{x - y}}
 \br{\frac{1}{\abs{x - y}^{d - 2}}
   \int_{B(x,\abs{x - y}/2)} \frac{d \mu(z)}{\abs{z - x}^{d - 1}} +
   \br{\frac{\abs{x - y}}{\rho_{\mu}(y)}}^{\delta_{\mu}} \frac{1}{\abs{x - y}^{d - 1}}},
\end{split}\end{align}
where the argument used to obtain equation
(7.15) from \cite{shen1999fundamental} was used to obtain the
final line.
 Since $x \in B$ and $y \in B'$, we will clearly have $\abs{x - y} \leq 8 \sqrt{d} r \leq 8
\sqrt{d} \kappa_{0} \rho_{\mu}(c) \leq 8 \sqrt{d} \rho_{\mu}(c)$. On
succesively applying Lemma \ref{lem:RHn2}, \eqref{eqtn:Measure1} and
the inclusion $B(x, 4 \sqrt{d} \rho_{\mu}(c)) \subset B(c, 12 \sqrt{d}
\rho_{\mu}(c))$, the estimate \eqref{eqtn:LocalHMW11}
will lead to
\begin{align*}\begin{split}  
& \int_{\R^{d}} \abs{\partial_{1} \Gamma_{\lambda}(x,z)} \Gamma_{\mu +
   \lambda}(z,y) \, d \mu(z)  \\ & \qquad \qquad \lesssim
 e^{- \frac{\varepsilon}{2} \sqrt{\lambda} \abs{x - y}} \br{\frac{1}{\abs{x - y}^{d - 2}}
   \frac{\mu(B(x,\abs{x - y} / 2))}{\abs{x - y}^{d - 1}} +
   \br{\frac{\abs{x - y}}{\rho_{\mu}(y)}}^{\delta_{\mu}} \frac{1}{\abs{x -
       y}^{d - 1}}} \\
 & \qquad \qquad \lesssim
 e^{- \frac{\varepsilon}{2}\sqrt{\lambda} \abs{x - y}} \br{\frac{1}{\abs{x - y}^{d - 1}}
   \br{\frac{\abs{x - y}}{\rho_{\mu}(c)}}^{\delta_{\mu}}
   \frac{\mu(B(x,4 \sqrt{d}  \rho_{\mu}(c)))}{\rho_{\mu}(c)^{d - 2}} +
   \br{\frac{\abs{x - y}}{\rho_{\mu}(y)}}^{\delta_{\mu}} \frac{1}{\abs{x -
       y}^{d - 1}}} \\
  & \qquad \qquad \lesssim
 e^{- \frac{\varepsilon}{2}\sqrt{\lambda} \abs{x - y}} \br{\frac{1}{\abs{x - y}^{d - 1}}
   \br{\frac{\abs{x - y}}{\rho_{\mu}(c)}}^{\delta_{\mu}}
   \frac{\mu(B(c,12 \sqrt{d}  \rho_{\mu}(c)))}{\rho_{\mu}(c)^{d - 2}} +
   \br{\frac{\abs{x - y}}{\rho_{\mu}(y)}}^{\delta_{\mu}} \frac{1}{\abs{x -
       y}^{d - 1}}} \\
 & \qquad \qquad \lesssim e^{-\frac{\varepsilon}{2}\sqrt{\lambda} \abs{x - y}}
 \br{\br{\frac{\abs{x - y}}{\rho_{\mu}(c)}}^{\delta_{\mu}} \frac{1}{\abs{x -
       y}^{d - 1}} + \br{\frac{\abs{x - y}}{\rho_{\mu}(y)}}^{\delta_{\mu}} \frac{1}{\abs{x -  y}^{d - 1}}},
\end{split}\end{align*}
where the final line follows from \eqref{eqtn:Measure2} and Remark
\ref{rmk:Critical}. Lemma \ref{lem:CriticalMeasure} and the bound
$\abs{x - y} \leq 8 \sqrt{d} \kappa_{0} \rho_{\mu}(c)$ then imply
$$
\int_{\R^{d}} \abs{\partial_{1} \Gamma_{\lambda}(x,z)} \Gamma_{\mu +
   \lambda}(z,y) \, d \mu(z) \lesssim  e^{-\frac{\varepsilon}{2} \sqrt{\lambda} \abs{x - y}}
  \frac{\kappa_{0}^{\delta_{\mu}}}{\abs{x - y}^{d-1}}
  $$
and therefore there must exist some $C > 0$, independent of $\kappa_{0}$, for which 
\begin{equation}
  \label{eqtn:LocalHMW2}
 \frac{1}{\pi} \int^{\infty}_{0} \lambda^{-\frac{1}{2}} \int_{\R^{d}}
 \abs{\partial_{1} \Gamma_{\lambda}(x,z)} \Gamma_{\mu + \lambda}(z,y)
 \, d \mu(z) \, d \lambda \leq 
 \frac{C \kappa_{0}^{\delta_{\mu}}}{\abs{x - y}^{d}}.
\end{equation}
For $x \in B'$ and $y \in B$, it is clear that
$x_{1} \geq y_{1}$ and therefore
$$
K_{0}^{(1)}(x,y) = - c \frac{\br{x_{1} - y_{1}}}{\abs{x - y}^{d + 1}} \leq 
- 
\frac{2 c'}{\abs{x - y}^{d}},
$$
for some $c, \, c' > 0$.
This can be combined with \eqref{eqtn:LocalHMW1} and
\eqref{eqtn:LocalHMW2} to produce the estimate
 $$
K_{\mu}^{(1)}(x,y) \leq 
\frac{\br{ - 2 c' +  C \kappa_{0}^{\delta_{\mu}}}}{\abs{x - y}^{d}}
$$
If we now set $\kappa_{0} = \br{\frac{c'}{C}}^{1 / \delta}$ we will
then obtain the estimate
$$
K_{\mu}^{(1)}(x,y) \leq - \frac{c'}{\abs{x - y}^{d}}
$$
for all $x \in B'$ and $y \in B$. Let $f \in L^{1}_{loc}(\R^{d})$ be a non-negative function with
support contained in $B$ that satisfies
$$
\dashint_{B} f > 0.
$$
Then for $x \in B'$,
\begin{align*}\begin{split}  
    \abs{R_{\mu}^{(1)}f(x)} &= \abs{\int_{B} K_{\mu}^{(1)}(x,y) f(y) \, dy} \\
    &\geq c' \int_{B} \frac{f(y)}{\abs{x - y}^{d}} \, dy \\
    &\geq c'' \dashint_{B} f(y) \, dy,
  \end{split}\end{align*}
for some $c'' > 0$.
This proves that for all $0 < \alpha < c'' \dashint_{B} f$ we must
have
$$
B' \subset \lb x \in \R^{d} : \abs{R_{\mu}^{(1)}f(x)} > \alpha \rb.
$$
The $L^{p}(w)$-boundedness of the operator $R_{\mu}^{(1)}$ can then be
exploited in order to obtain the estimate
$$
w(B') \lesssim \frac{1}{\alpha^{p}} \int_{B} f(x)^{p} w(x) \, dx
$$
for all $\alpha < c'' \dashint_{B} f$. On letting $\alpha \rightarrow
c'' \dashint_{B} f$,
\begin{equation}
  \label{eqtn:LocalHMW3}
\br{\dashint_{B} f}^{p} \lesssim \frac{1}{w(B')} \int_{B} f(x)^{p}
w(x) \, dx.
\end{equation}
By reversing the roles of $B$ and $B'$ in the previous argument, we
will also find that the estimate
$$
\br{\dashint_{B'} g}^{p} \lesssim \frac{1}{w(B)} \int_{B'} g(x)^{p}
w(x) \, dx
$$
must be valid for all non-negative $g$ supported in $B'$ with $\dashint_{B'} g >
0$. By setting $g =
\mathbbm{1}_{B'}$ in this estimate we find that $w(B) \lesssim
w(B')$. On applying this to \eqref{eqtn:LocalHMW3},
$$
\br{\dashint_{B} f}^{p} \lesssim \frac{1}{w(B)} \int_{B} f(x)^{p} w(x)
\, dx.
$$
Setting $f := (w + \varepsilon)^{-\frac{1}{p - 1}} \mathbbm{1}_{B}$
then leads to 
$$
w(B) \br{\dashint_{B} (w + \varepsilon)^{-\frac{1}{p - 1}}}^{p} \lesssim \int_{B}
 (w + \varepsilon)^{-\frac{p}{p - 1}}(x) w(x) \, dx \leq \int_{B} (w +
 \varepsilon)^{-\frac{1}{p - 1}}(x) \, dx.
$$
By letting $\varepsilon \rightarrow 0$, the monotone convergence
theorem then gives us \eqref{eqtn:LocalHMW12}.
 \end{prof}

\subsection{Constant Potentials}

Throughout this section, set $V \equiv N$ for some $N > 0$. Notice that
for this case the Agmon distance will be given by
 $$
d_{V}(x,y) = N^{\frac{1}{2}} \abs{x - y} \qquad \forall \, x, \, y \in \R^{d}.
$$
This straightforward expression for the Agmon distance leads to a
simpler characterisation of our weight class $S_{p,c}^{V}$.

\begin{prop} 
 \label{prop:ConstantAp} 
Let $1 < p < \infty$ and $c > 0$. A weight $w$ will be contained in
the class $S_{p,c}^{N}$ if and only if
\begin{equation}
  \label{eqtn:ConstantAlternate}
\sup_{B} \br{\frac{1}{\abs{B} e^{c N^{\frac{1}{2}} r}} \int_{B} w}^{\frac{1}{p}}
\br{\frac{1}{\abs{B}e^{c N^{\frac{1}{2}} r}} \int_{B} w^{-\frac{1}{p -
      1}}}^{\frac{p - 1}{p}} < \infty,
\end{equation}
where the supremum is taken over all Euclidean balls $B = B(x,r) \in
\R^{d}$ with $x \in \R^{d}$ and $r > 0$.
\end{prop}

\begin{prof}  
  The proof follows entirely from the fact that
  $$
B(x,r) = B_{V}(x,N^{\frac{1}{2}}r) 
$$
for all $x \in \R^{d}$ and $r > 0$. Indeed, first suppose that $w \in
S_{p,c}^{N}$. Let $B(x,r) \subset \R^{d}$ be a ball in the Euclidean
metric. Then
\begin{align*}\begin{split}  
 w (B(x,r))^{\frac{1}{p}} w^{-\frac{1}{p - 1}} \br{B(x,r)}^{\frac{p - 1}{p}} &=
 w \br{B_{V}(x,N^{\frac{1}{2}}r)}^{\frac{1}{p}} w^{-\frac{1}{p - 1}}
 \br{B_{V}(x,N^{\frac{1}{2}} r)}^{\frac{p - 1}{p}} \\
 &\lesssim e^{c N^{\frac{1}{2}} r} \abs{B_{V}(x,N^{\frac{1}{2}} r)} \\
 &= e^{c N^{\frac{1}{2}}r} \abs{B(x,r)}.
\end{split}\end{align*}
Let's now prove the converse. Suppose that
\eqref{eqtn:ConstantAlternate} is finite. Then for any ball
$B_{V}(x,r) \subset \R^{d}$ in the metric $d_{V}$,
\begin{align*}\begin{split}  
 w (B_{V}(x,r))^{\frac{1}{p}} w^{-\frac{1}{p - 1}} \br{B_{V}(x,r)}^{\frac{p -
     1}{p}} &= w \br{B(x,N^{-\frac{1}{2}}r)}^{\frac{1}{p}}
 w^{-\frac{1}{p - 1}} \br{B(x,N^{-\frac{1}{2}}r)}^{\frac{p - 1}{p}} \\
 &\lesssim e^{c r} \abs{B(x,N^{-\frac{1}{2}}r)} \\
 &= e^{c r} \abs{B_{V}(x,r)},
\end{split}\end{align*}
which proves that $w \in S_{p,c}^{N}$.
 \end{prof}

In order to prove the
first chain of inclusions of Conjecture \ref{conj:Upper}, we need to
know more about the behaviour of the singular
kernel of $R_{V} = R_{N}$. This is provided by the following lemma.

\begin{lem}
  \label{lem:Bessel}
For $1 \leq j \leq d$,  the singular kernel of the Riesz transform $R_{N}^{(j)} := \partial_{j} \br{- \Delta + N}^{-\frac{1}{2}}$ is given by
$$
K_{N}^{(j)}(x,y) = - c_{N} \frac{\br{x_{j} - y_{j}}}{\abs{x - y}} e^{- N^{\frac{1}{2}}\abs{x - y}}
s(N^{\frac{1}{2}} \abs{x - y}),
$$
for all $x, \, y \in \R^{d}$, where $c_{N} > 0$ and $s : (0,\infty)
\rightarrow [0,\infty)$ is the function defined through
$$
s(a) := \int^{\infty}_{0} e^{- a t} \br{t +
  \frac{t^{2}}{2}}^{\frac{d - 2}{2}} \, dt + \int^{\infty}_{0} t e^{-
  a t} \br{t + \frac{t^{2}}{2}}^{\frac{d - 2}{2}} \, dt.
$$
\end{lem}

\begin{prof}  
  The operator $\br{- \Delta + N}^{-\frac{1}{2}}$ is given by
  convolution with the function
  $$
G^{N}(z) = c_{N} \int^{\infty}_{0} e^{- N t} e^{- \frac{\abs{z}^{2}}{4 t}}
t^{\frac{1 - d}{2}} \, \frac{dt}{t},
$$
for some constant $c_{N} > 0$ (c.f. {\cite[pg.~7]{grafakos2009modern}}). Through a change of variables, it is
easy to see that
\begin{equation}
  \label{eqtn:Bessel1}
G^{N}(z) \simeq G^{1}(N^{\frac{1}{2}}z)
\end{equation}
for all $z \in \R^{d}$. The function $G^{1}$ is the well-known
Bessel kernel of order one and has the representation
$$
G^{1}(z) = c e^{- \abs{z}} \int^{\infty}_{0} e^{- \abs{z} t} \br{t +
  \frac{t^{2}}{2}}^{\frac{d - 2}{2}} \, \frac{dt}{t}
$$
for all $z \in \R^{d}$. The proof of this representation can be found
in \cite{aronszajn1961theory}. This, together with \eqref{eqtn:Bessel1} and
simple differentiation then proves our lemma.
\end{prof}

Notice that for the constant potential $V \equiv N$, the associated
Hardy-Littlewood operator for $c > 0$, as defined in Definition
\ref{def:AdaptedAveraging}, is given by
$$
M_{N,c}f(x) := \sup_{t > 0} A_{t,c}^{N} \abs{f}(x) := \sup_{t > 0}
\frac{1}{\abs{B(x,N^{-\frac{1}{2}} t)} e^{c t}}
\int_{B(x,N^{-\frac{1}{2}} t)} \abs{f(y)} \, dy,
$$
where we are using the shorthand notation $M_{N,c} = M_{\rho_{N},c}$
and $A_{t,c}^{N} = A_{t,c}^{\rho_{N}}$.
Let $R_{0,loc}$, $T^{*}_{0,loc}$ and $M_{0,loc}$ denote the $\rho_{N}$-localized
parts of the classical Riesz transform, heat and Hardy-Littlewood
maximal operators
respectively. That is,
$$
R_{0,loc}f(x) := \nabla \br{-
  \Delta}^{-\frac{1}{2}} \br{f \cdot \mathbbm{1}_{B(x,\rho_{N}(x))}}(x), \qquad T^{*}_{0,loc}f(x) := \sup_{t > 0} e^{t \Delta} \abs{f \cdot
  \mathbbm{1}_{B(x,\rho_{N}(x))}}(x)
$$
and
$$
M_{0,loc}f(x) := \sup_{t >
  0} A_{t}^{0} \abs{f \mathbbm{1}_{B(x,\rho_{N}(x))}}(x) := \sup_{t > 0}
\frac{1}{\abs{B(x,t)}} \int_{B(x,t)} \abs{f \cdot \mathbbm{1}_{B(x,\rho_{N}(x))}}.
$$
The following proposition is an embodiment of the idea that the
operators attached to the Schr\"{o}dinger operator $-\Delta+N$ should
behave, at a local scale, like their classical counterparts. In
particular, a Muckenhoupt-type equivalence will hold at a local level
between all of the operators listed.
  
\begin{prop} 
 \label{prop:Local} 
 For any $c > 0$, the following statements are equivalent:
 \begin{enumerate}[(i)]
 \item $w \in A_{p}^{\rho_{N},loc}$;
 \item $T^{*}_{0,loc}$ is bounded on $L^{p}(w)$;
 \item $T^{*}_{N,loc}$ is bounded on $L^{p}(w)$;
 \item $M_{0,loc}$ is bounded on $L^{p}(w)$;
 \item $M_{N,c}^{loc}$ is bounded on $L^{p}(w)$;
 \item $R_{0,loc}$ is bounded on $L^{p}(w)$;
   \item $R_{N,loc}$ is bounded on $L^{p}(w)$.
   \end{enumerate}
 \end{prop}

 \begin{prof}  
 $(i) \Rightarrow (ii)$. This is proved in {\cite[Thm.~1]{bongioanni2011classes}}.

\vspace*{0.1in}

$(ii) \Rightarrow (iii)$. This follows trivially from the fact that the
heat kernel of the operator $-\Delta + N$ is pointwise bounded from
above by the heat kernel of  $-\Delta$.

\vspace*{0.1in}

$(iii) \Rightarrow (iv)$. Fix a ball $B := B(x,r) \subset \R^{d}$ with $r
\leq \rho_{N}(x) = N^{-\frac{1}{2}}$. For
any other point $y \in B$ we will clearly have
$$
\exp \br{ - \frac{\abs{x - y}^{2}}{4 r^{2}}} \simeq 1.
$$
Also note that since $r \leq N^{-\frac{1}{2}}$,
$$
e^{-N r^{2}} \simeq 1.
$$
This gives
$$
 \frac{1}{\abs{B}} \int_{B} \abs{f(y)} \, dy \lesssim \frac{1}{e^{N r^{2}}} \frac{1}{r^{n}} \int_{B} \exp\br{-\frac{\abs{x -
       y}^{2}}{4 r^{2}}} \abs{f(y)} \, dy.
 $$
 Set $t = r^{2}$ to obtain
 \begin{align*}\begin{split}  
 \frac{1}{\abs{B}} \int_{B} \abs{f(y)} \, dy &\lesssim \frac{1}{e^{N
     t} t^{\frac{n}{2}}} \int_{B} \exp \br{- \frac{\abs{x - y}^{2}}{4
     t}} \abs{f(y)} \, dy \\
 &\leq e^{- t(N - \Delta)} \abs{f \mathbbm{1}_{B(x,\rho_{N}(x))}}(x) \\
 &\lesssim T^{*}_{N,loc} f(x).
\end{split}\end{align*}
This demonstrates that $\norm{T^{*}_{N,loc}}_{L^{p}(w)} < \infty$ implies
$\norm{M_{0,loc}}_{L^{p}(w)} < \infty$.

\vspace*{0.1in}

$(iv) \Rightarrow (v)$. This implication follows trivially from the
fact that the kernel of $A_{t,c}^{N}$ is pointwise bounded from above by
the kernel of $A_{N^{-\frac{1}{2}} t}^{0}$ for all $t > 0$.

\vspace*{0.1in}

$(v) \Rightarrow (i)$. Suppose that $M_{N,c}^{loc}$ is bounded on
$L^{p}(w)$. First note that from the argument at the beginning of
Theorem \ref{thm:LocalHMW}, in order to show that $w \in
A_{p}^{\rho_{N},loc}$ it is sufficient to prove that $w \in
A_{p}^{\frac{1}{2} \rho_{N},loc}$. The proof that $w \in A_{p}^{\frac{1}{2}\rho_{N},loc}$
is essentially identical to the
classical proof that can be found in {\cite[pg.~280]{grafakos2009modern}}. Fix $B := B(x,r) \subset \R^{d}$ with $r \leq \frac{1}{2}\rho_{N}(x)
= \frac{1}{2} N^{-\frac{1}{2}}$. Let $\mathcal{M}_{N,c} =
\mathcal{M}_{\rho_{N},c}$ be the uncentered Hardy-Littlewood operator
as defined in Definition \ref{def:AdaptedAveraging}. Then since $B
\subset B(y,\rho_{N}(y))$ for all $y \in B$,
\begin{align*}\begin{split}
    w(B) \br{\frac{1}{\abs{B}} \int_{B} \abs{f}}^{p} &\lesssim
 w(B) \br{\frac{1}{e^{c N^{\frac{1}{2}}r} \abs{B}} \int_{B}
   \abs{f}}^{p} \\
&= \int_{B} \br{\frac{1}{e^{c N^{\frac{1}{2}} r} \abs{B}} \int_{B}
  \mathbbm{1}_{B(y,\rho_{N}(y))} \abs{f}}^{p} w(y) \, dy \\
 &\leq \int_{B} \mathcal{M}_{N,c}^{loc}(f
 \mathbbm{1}_{B})(y)^{p} w(y) \, dy \\
 &\lesssim \int_{B} M_{N,c'}^{loc} (f \mathbbm{1}_{B})(y)^{p} w(y) \, dy,
\end{split}\end{align*}
for $c' < \frac{c}{2}$, where the last line follows from Proposition
\ref{prop:ConverseCenter}. It is not difficult to see that since the
operator $M_{N,c'}^{loc}$ is localized we must have
$$
M_{N,c'}^{loc}(f \mathbbm{1}_{B})(y) \lesssim M_{N,c}^{loc}(f \mathbbm{1}_{B})(y).
$$
The boundedness of $M_{N,c}^{loc}$ then implies that
\begin{align*}\begin{split}  
  w(B) \br{\frac{1}{\abs{B}} \int_{B} \abs{f}}^{p} 
  &\lesssim \int_{B} M_{N,c}^{loc} (f \mathbbm{1}_{B})(y)^{p} w(y) \,
  dy \\
  &\lesssim \int_{B} \abs{f(y)}^{p} w(y) \, dy.
 \end{split}\end{align*}
 For $\epsilon > 0$, take $f := \br{w +
  \epsilon}^{-\frac{1}{p-1}}$ in the above inequality to obtain
$$
\frac{1}{\abs{B}^{p}} w(B) \br{\int_{B} \br{w
    + \epsilon}^{-\frac{1}{p-1}}}^{p} \lesssim  \br{\int_{B} \br{w +
    \epsilon}^{-\frac{p}{p-1}} w } \leq  \br{\int_{B} (w + \varepsilon)^{-\frac{1}{p-1}}}.
$$
Which leads to
$$
\br{\frac{1}{ \abs{B}} \int_{B} w} \br{
  \frac{1}{\abs{B}}\int_{B} \br{w
    + \epsilon}^{-\frac{1}{p-1}}}^{p - 1} \leq C,
$$
for some constant $C > 0$. The monotone convergence theorem then completes the proof of this
implication.

\vspace*{0.1in}

$(i) \Rightarrow (vi), \, (i) \Rightarrow (vii)$. The
proof of these two implication is essentially contained in
{\cite[Thm.~1]{bongioanni2011classes}} and {\cite[Thm.~3]{bongioanni2011classes}}.

\vspace*{0.1in}

$(vi) \Rightarrow (i)$. This implication can be proved using a
classical argument that can be found, for example, in
{\cite[Thm.~9.4.8]{grafakos2009modern}}.

\vspace*{0.1in}

$(vii) \Rightarrow (i)$. This is proved in Theorem \ref{thm:LocalHMW}.

 \end{prof}

Let's now prove the first chain of inclusions in Conjecture
\ref{conj:Upper} for the case $V \equiv N$.

\begin{thm} 
 \label{thm:ConstantSharpRiesz} 
Let $1 < p < \infty$. There exists $c_{1}, \, c_{2} > 0$, independent
of $p$ and $N$, such that
 $$
S_{p,c_{1}}^{N} \subset \lb w : \norm{R_{N}}_{L^{p}(w)} < \infty \rb
\subset S_{p,c_{2}}^{N}
 $$
 \end{thm}

 \begin{prof}  
 The first inclusion has already been proved in Theorem \ref{thm:Riesz} and so it suffices to
 consider the second inclusion. Suppose that $w$ is a weight on
 $\R^{d}$ for which $\norm{R_{N}}_{L^{p}(w)} < \infty$. We will adapt the classical
 proof of {\cite[Thm.~9.4.9]{grafakos2009modern}}.

\vspace*{0.1in}
 
 Consider the operator $Wf(x) := \sum_{j = 1}^{d} R_{N}^{(j)}f(x)$. Let $B =
 B(c, r)$ be
 a ball in $\R^{d}$ and $f \in L^{1}_{loc}(\R^{d})$ a non-negative function with support
 contained in $B$ that satisfies
 $$
\dashint_{B} f > 0.
$$
Let $B'$ be the ball in $\R^{d}$ that has center $c' = c + 2 r
= (c_{1} + 2r, \cdots, c_{d} + 2 r)$ and radius $r$. Clearly
$B'$ will satisfy $x_{j} \geq y_{j}$ for each $1
\leq j \leq d$ when $x = (x_{1}, \cdots, x_{d}) \in B'$ and $y =
(y_{1}, \cdots, y_{d}) \in B$. Lemma \ref{lem:Bessel} then implies
that for $x \in B'$,
\begin{align*}\begin{split}  
 \abs{Wf(x)} &\simeq \sum_{j = 1}^{d} \int_{B} \frac{x_{j} - y_{j}}{\abs{x
     - y}} e^{- N^{\frac{1}{2}} \abs{x - y}} s(N^{\frac{1}{2}}\abs{x - y}) f(y) \, dy \\
 &\geq \int_{B} e^{- N^{\frac{1}{2}} \abs{x - y}} s(N^{\frac{1}{2}} \abs{x - y}) f(y) \, dy.
\end{split}\end{align*} 
It isn't too difficult to check that the function $s$ satisfies the
lower bound
$$
s(a) \gtrsim \frac{1}{a^{d}}
$$
for all $a > 0$. Indeed, this follows from
\begin{align*}\begin{split}  
 s(a) &\geq \int^{\infty}_{0} t e^{- a t} \br{t +
   \frac{t^{2}}{2}}^{\frac{d - 2}{2}} \, dt \\
 &= \int^{\infty}_{0} \frac{t}{a} e^{- t} \br{\frac{t}{a} +
   \frac{t^{2}}{2 a^{2}}}^{\frac{d - 2}{2}} \, \frac{dt}{a} \\
 &\geq \frac{1}{a^{2}} \int^{\infty}_{0} t e^{-t} \br{\frac{t^{2}}{2
     a^{2}}}^{\frac{d - 2}{2}} \, dt \\
 &\simeq \frac{1}{a^{d}} \int^{\infty}_{0} t^{d - 1} e^{-t} \, dt \\
 &\simeq \frac{1}{a^{d}}.
 \end{split}\end{align*}
Therefore, for $x \in B'$ we will have
\begin{align*}\begin{split}  
 \abs{Wf(x)} &\geq D' \int_{B} \frac{e^{-N^{\frac{1}{2}} \abs{x - y}}}{\abs{x - y}^{d}}
 f(y) \, dy \\
 &\geq D e^{- 4 \sqrt{d} r N^{\frac{1}{2}}} \dashint_{B} f(y) \, dy,
 \end{split}\end{align*}
for some constants $D, \, D' > 0$. This implies that for any $ 0 < \alpha < D
e^{- 4 \sqrt{d} r N^{\frac{1}{2}}}\dashint_{B}f$ we will have
$$
B' \subset \lb x \in \R^{d} : \abs{Wf(x)} > \alpha \rb.
$$
The $L^{p}(w)$-boundedness of the operator $W$ will imply
$$
w(B') \lesssim \frac{1}{\alpha^{p}} \int_{B} f(x)^{p} w(x) \, dx
$$
for all $\alpha < D e^{- 4 \sqrt{d} r N^{\frac{1}{2}}} \dashint_{B}
f$ which then gives
\begin{equation}
  \label{eqtn:ConstantSharpRiesz1}
\br{\dashint_{B} f}^{p} e^{- 4 p \sqrt{d} r N^{\frac{1}{2}}} \lesssim
\frac{1}{w(B')} \int_{B} f(x)^{p} w(x) \, dx.
\end{equation}
The roles of $B$ and $B'$ can be reversed to obtain
$$
\br{\dashint_{B'} g}^{p} e^{- 4 p \sqrt{d} r N^{\frac{1}{2}}} \lesssim
\frac{1}{w(B)} \int_{B'} g(x)^{p} w(x) \, dx
$$
for all non-negative $g$ supported in $B'$ with $\dashint_{B'} g >
0$. Setting $g = \mathbbm{1}_{B'}$ in the above estimate then gives
$w(B) e^{- 4 p \sqrt{d} r N^{\frac{1}{2}}} \lesssim w(B')$. Applying this to
\eqref{eqtn:ConstantSharpRiesz1},
$$
w(B) \br{\dashint_{B} f}^{p} e^{- 8 p \sqrt{d} r N^{\frac{1}{2}}} \lesssim
 \int_{B} f(x)^{p} w(x) \, dx.
$$
For $\varepsilon > 0$, set $f := (w + \varepsilon)^{-\frac{1}{p - 1}}
\mathbbm{1}_{B}$ in the above estimate to obtain
\begin{align*}\begin{split}  
 w(B) \br{\frac{1}{\abs{B}} \int_{B} \br{w + \varepsilon}^{-\frac{1}{p
     - 1}}}^{p} e^{- 8 \sqrt{d} p r N^{\frac{1}{2}}} &\lesssim \int_{B} (w + \varepsilon)^{-\frac{p}{p
   - 1}} w \\
&\leq \int_{B} (w + \varepsilon)^{-\frac{1}{p - 1}}.
\end{split}\end{align*}
The monotone convergence theorem then yields
$$
w \br{B}^{\frac{1}{p}} w^{-\frac{1}{p - 1}}\br{B}^{\frac{p - 1}{p}}
\lesssim e^{8 \sqrt{d} N^{\frac{1}{2}} r} \abs{B}. 
$$
Proposition \ref{prop:ConstantAp} then implies that $w \in S_{p,8 \sqrt{d}}^{N}$.
 \end{prof}

  \subsection{Potentials Bounded from Above and Below}

In this section, it will be proved that the second chain of inclusions
in Conjecture \ref{conj:Upper} holds for potentials that are bounded both from
above and from below. Suppose that there exists $N, \, M > 0$ for
which $M \leq V(x) \leq N$ for all $x \in \R^{d}$. We will require the
following lemma.

\begin{lem}
\label{lem:AboveBelow}
The heat kernel for the Schr\"{o}dinger operator $-\Delta + V$
satisfies the estimate
$$
e^{-N t} \frac{e^{-\frac{\abs{x - y}^{2}}{4 t}}}{t^{\frac{d}{2}}} \lesssim k_{t}^{V}(x,y) \lesssim e^{-M t} \frac{e^{-\frac{\abs{x - y}^{2}}{4 t}}}{t^{\frac{d}{2}}}
$$
for a.e. $x, \, y \in \R^{d}$ and $t > 0$.
\end{lem}

\begin{prof}  
 Let $V_{1}$ and $V_{2}$ be two potentials with $0 \leq V_{1}(x) \leq
 V_{2}(x)$ for all $x \in \R^{d}$. To prove our lemma it is sufficient
 to show that $k_{t}^{V_{2}}(x,y) \leq k_{t}^{V_{1}}(x,y)$ for each $t
 > 0$ and a.e. $x, \, y \in \R^{d}$. Theorem 2.24 of
 \cite{ouhabaz2005analysis} states that the semigroup of $(V_{2} -
 \Delta)$ will be dominated by the semigroup of $(V_{1} -
 \Delta)$. That is,
 \begin{equation}
   \label{eqtn:Domination}
\abs{e^{-t (V_{2} - \Delta)}f} \leq e^{-t (V_{1} - \Delta)}\abs{f}
\end{equation}
for all $f \in L^{2}(\R^{d})$. Fix $(x,t) \in \R^{d} \xx (0,\infty)$
and suppose that there exists compactly supported $E \subset \R^{d}$ such that
$$
k^{V_{1}}_{t}(x,y) < k^{V_{2}}_{t}(x,y)
$$
for all $y \in E$. Then applying \eqref{eqtn:Domination} to the
function $f(y) = \mathbbm{1}_{E}$ yields
$$
\int_{E} k^{V_{2}}_{t}(x,y) \, dy \leq \int_{E} k^{V_{1}}_{t}(x,y) \, dy,
$$
implying $\abs{E} = 0$. Therefore $k^{V_{2}}_{t}(x,y) \leq
k^{V_{1}}_{t}(x,y)$ for each $t > 0$, for almost every $x, \, y \in
\R^{d}$. 
 \end{prof}

\begin{lem}
  \label{lem:ClassAB}
  For any $c > 0$ and $1 < p < \infty$, the following chain of inclusions holds,
  $$
S_{p,c}^{M} \subset S_{p,c}^{V} \subset S_{p,c}^{N}.
  $$
\end{lem}

\begin{prof}
  It is clear from the definition of the Agmon distance that the
  inequality $M \leq V(x) \leq N$ will imply
  $$
d_{M}(x,y) \leq d_{V}(x,y) \leq d_{N}(x,y)
$$
for any $x, \, y \in \R^{d}$. Therefore,
$$
B_{N}(x,r) \subset B_{V}(x,r) \subset B_{M}(x,r)
$$
for all $x \in \R^{d}$ and $r > 0$.

Suppose that $w \in S_{p,c}^{M}$. Then for any $x \in \R^{d}$
and $r > 0$,
 \begin{align*}\begin{split} 
 w \br{B_{V}(x,r)}^{\frac{1}{p}} w^{-\frac{1}{p - 1}}
 \br{B_{V}(x,r)}^{\frac{p - 1}{p}} &\leq w
 \br{B_{M}(x,r)}^{\frac{1}{p}} w^{-\frac{1}{p - 1}}
 \br{B_{M}(x,r)}^{\frac{p - 1}{p}} \\
 &\lesssim e^{c r} \abs{B_{M}(x,r)} \\
 &\simeq e^{c r} \abs{B_{N}(x,r)} \\
 &\leq e^{c r} \abs{B_{V}(x,r)}.
\end{split}\end{align*}
Suppose instead that $w \in S_{p,c}^{V}$. Then for any $x \in \R^{d}$
and $r > 0$,
\begin{align*}\begin{split}  
 w (B_{N}(x,r))^{\frac{1}{p}} w^{-\frac{1}{p - 1}}
 (B_{N}(x,r))^{\frac{p - 1}{p}} &\leq w (B_{V}(x,r))^{\frac{1}{p}}
 w^{-\frac{1}{p - 1}}(B_{V}(x,r))^{\frac{p-1}{p}} \\
 &\lesssim e^{c r} \abs{B_{V}(x,r)} \\
 &\leq e^{c r} \abs{B_{M}(x,r)} \\
 &\simeq e^{c r} \abs{B_{N}(x,r)}.
 \end{split}\end{align*}
 \end{prof}

 \begin{prop} 
 \label{prop:ConstantPotentialMuckenhoupt} 
 Let $1 < p < \infty$. We have
 the following chain of inclusions for some $c_{1}, \, c_{2}, \, c_{3}
 > 0$ with $c_{1} < c_{2} < c_{3}$,
 $$
S_{p,c_{1}}^{V} \subset \lb w : \norm{T^{*}_{V}}_{L^{p}(w)} < \infty \rb \subset
\lb w : \norm{M_{N,c_{2}}}_{L^{p}(w)} < \infty \rb \subset S_{p,c_{3}}^{V}.
$$
The constants $c_{1}, \, c_{2}$ and $c_{3}$ will be independent of $p$.
\end{prop}

\begin{prof}
  First let's consider the last inclusion in the above chain.
 The implication that $\norm{M_{N,c_{2}}}_{L^{p}(w)} < \infty$ gives
 $w \in S_{p,3 c_{2}}^{N}$ is asserted by Proposition \ref{prop:HardyLittlewood}.  It is obvious that there
 must exist some $c_{3} > 3 c_{2}$ for which $S_{p,3 c_{2}}^{N} \subset S_{p,c_{3}}^{M}$. The last inclusion of our
 proposition then follows from Lemma \ref{lem:ClassAB}.

 \vspace*{0.1in}

Next, let's prove the second inclusion. Let $w$ be a weight on
$\R^{d}$ for which $\norm{T^{*}_{V}}_{L^{p}(w)} < \infty$. 
 Set $c_{2} = N^{\frac{1}{2}} +
\frac{N^{-\frac{1}{2}}}{4}$ and fix $f \in L^{p}(w)$. Then for $0 < t < N^{-\frac{1}{2}}$ we have
\begin{align*}\begin{split}  
 A_{N^{\frac{1}{2}} t,c_{2}}^{N} \abs{f}(x) &= \frac{1}{\abs{B(x,t)}e^{c_{2} N^{\frac{1}{2}}
     t}} \int_{B(x,t)} \abs{f(y)} \, dy \\
 &\leq \frac{1}{t^{d}} \int_{B(x,t)} \abs{f(y)} \, dy \\
 &\lesssim \frac{1}{t^{d} e^{N t^{2}}} \int_{B(x,t)} \abs{f(y)} \, dy.
\end{split}\end{align*}
For $y \in B(x,t)$ we have
$$
\exp \br{-\frac{\abs{x - y}^{2}}{4 t^{2}}} \simeq 1
$$
and therefore, on applying Lemma \ref{lem:AboveBelow},
\begin{align*}\begin{split}  
  A_{N^{\frac{1}{2}} t,c_{2}}^{N} \abs{f}(x) &\lesssim \frac{1}{t^{d}
    e^{N t^{2}}} \int_{B(x,t)} \exp \br{- \frac{\abs{x - y}^{2}}{4 t^{2}}} \abs{f(y)} \, dy \\
  &\lesssim \int_{B(x,t)} k_{t^{2}}^{V}(x,y)  \abs{f(y)} \, dy \\
  &\lesssim T^{*}_{V}f(x).
\end{split}\end{align*}
Next, consider the case $t \geq N^{-\frac{1}{2}}$. We have
\begin{align*}\begin{split}  
 A_{N^{\frac{1}{2}}t,c_{2}}^{N} \abs{f}(x) &= \frac{1}{\abs{B(x,t)}e^{c_{2}
     N^{\frac{1}{2}} t}} \int_{B(x,t)} \abs{f(y)} \, dy \\
 &\simeq \frac{1}{t^{d} e^{N t} e^{\frac{t}{4}}} \int_{B(x,t)}
 \abs{f(y)} \, dy.
 \end{split}\end{align*}
For $y \in B(x,t)$ we must have
$$
\exp \br{-\frac{t}{4}} \leq \exp \br{- \frac{\abs{x - y}^{2}}{4 t}}.
$$
Therefore,
\begin{equation}
  \label{eqtn:AtNEst}
A_{N^{\frac{1}{2}} t,c_{2}}^{N} \abs{f}(x) \lesssim \frac{1}{t^{d} e^{N t}} \int_{B(x,t)}
e^{-\frac{\abs{x - y}^{2}}{4 t}} \abs{f(y)} \, dy.
\end{equation}
Notice that since $t \geq
N^{-\frac{1}{2}}$, \eqref{eqtn:AtNEst} gives
\begin{align*}\begin{split}  
 A_{N^{\frac{1}{2}} t,c_{2}}^{N} \abs{f}(x) 
 &\lesssim \frac{1}{t^{\frac{d}{2}} e^{N t}} \int_{B(x,t)}
 e^{-\frac{\abs{x - y}^{2}}{4 t}} \abs{f(y)} \, dy \\
 &\lesssim  e^{-t \br{V - \Delta}} \abs{f}(x) \\
 &\lesssim T^{*}_{V}f(x),
\end{split}\end{align*}
where Lemma \ref{lem:AboveBelow} was applied in the second line.
This proves that
$$
M_{N,c_{2}}f(x) \lesssim T^{*}_{V}f(x)
$$
for all $x \in \R^{d}$ and thus completes the proof of the second inclusion.

 \vspace*{0.1in}

 Finally, let's prove the first inclusion. Suppose that $w \in
 S_{p,c_{1}}^{V}$ for some $c_{1} > 0$, $1 < p < \infty$ and fix $f
\in L^{p}(w)$. Proposition \ref{prop:LocalApVc} tells us that
$S_{p,c_{1}}^{V} \subset A_{p}^{V,loc}$. It follows from the fact that
$T^{*}_{V}$ is pointwise bounded from above by the heat maximal
operator for the Laplacian $T^{*}_{0}$ and
{\cite[Thm~1]{bongioanni2011classes}} that
$\norm{T^{*}_{V,loc}}_{L^{p}(w)} < \infty$. It remains to show that
the global part of $T^{*}_{V}$ is bounded on $L^{p}(w)$.

 Let $B_{j} := B(x_{j},\rho_{V}(x_{j}))$ be as given in
Proposition \ref{prop:Cover} and set $B_{x} := B(x,\rho_{V}(x))$ for
each $x \in \R^{d}$.  On expanding
the $L^{p}(w)$-norm of $T^{*}_{V,glob}f$ and applying Lemma \ref{lem:AboveBelow},
\begin{align}\begin{split}
    \label{eqtn:AboveBelow1}
 \norm{T^{*}_{V,glob}f}_{L^{p}(w)} &\lesssim \br{\sum_{j}
   \int_{B_{j}} \br{\sup_{t > 0} \frac{e^{-M t}}{t^{\frac{d}{2}}} \int_{B_{x}^{c}} 
  e^{-\frac{\abs{x - y}^{2}}{4 t}}    \abs{f(y)} \, dy}^{p} w(x) \,
dx}^{\frac{1}{p}} \\
&= \br{\sum_{j}
   \int_{B_{j}} \br{\sup_{t > 0} \sum_{k = 1}^{\infty} \frac{e^{-M
         t}}{t^{\frac{d}{2}}} \int_{(k + 1)B_{x} \setminus k B_{x}} 
  e^{-\frac{\abs{x - y}^{2}}{4 t}}    \abs{f(y)} \, dy}^{p} w(x) \,
dx}^{\frac{1}{p}}.
\end{split}\end{align}
Observe that the critical radius function satisfies
$N^{-\frac{1}{2}}
\leq \rho_{V}(x)
\leq M^{-\frac{1}{2}}$ for all $x \in \R^{d}$. Therefore, for $j \in
\N$, $k \geq 1$, $x \in B_{j}$ and $y \notin k B_{x}$,
$$
\abs{x - y} \geq k \rho_{V}(x) \geq k N^{-\frac{1}{2}}.
$$
Applying this bound to \eqref{eqtn:AboveBelow1} gives
\begin{equation}
  \label{eqtn:AboveBelow2}
\norm{T^{*}_{V,glob}f}_{L^{p}(w)} \lesssim \br{\sum_{j}
   \int_{B_{j}} \br{\sup_{t > 0} \sum_{k = 1}^{\infty} \frac{e^{-M
         t}}{t^{\frac{d}{2}}} e^{-\frac{k^{2}}{4 N t}} \int_{(k + 1)B_{x} \setminus k B_{x}} \abs{f(y)} \, dy}^{p} w(x) \,
dx}^{\frac{1}{p}}.
\end{equation}
For $y \in (k + 1)B_{x}$, Lemma \ref{lem:Shen0} implies that
\begin{align*}\begin{split}  
    \abs{x_{j} - y} &\leq \abs{x - y} + \abs{x_{j} - x} \\
    &\leq (k + 1) \rho_{V}(x) + \rho_{V}(x_{j}) \\
    &\leq \beta 2^{\frac{k_{0}}{k_{0} + 1}}(k + 1) \rho_{V}(x_{j}) +
    \rho_{V}(x_{j}) \\
    &\leq 4 \sigma k \rho_{V}(x_{j}),
  \end{split}\end{align*}
where $\sigma := \beta 2^{\frac{k_{0}}{k_{0} + 1}}$. Therefore $(k +
1) B_{x} \subset 4 \sigma k B_{j}$. This inclusion together with $B_{j}
\subset B(x_{j},M^{-\frac{1}{2}})$ and
H\"{o}lder's inequality then leads to 
\begin{align}\begin{split}  
    \label{eqtn:AboveBelow3}
    &\norm{T^{*}_{V,glob}f}_{L^{p}(w)} \lesssim \\
     &\br{\sum_{j} \br{\sup_{t > 0} \sum_{k = 1}^{\infty} \frac{e^{-M
         t}}{t^{\frac{d}{2}}} e^{-\frac{k^{2}}{4 N t}} w^{-\frac{1}{p
         - 1}}(B(x_{j}, 4 \sigma k  M^{-\frac{1}{2}}))^{\frac{p - 1}{p}}
     w(B(x_{j}, M^{-\frac{1}{2}}))^{\frac{1}{p}}\norm{f}_{L^{p}(4 \sigma k \cdot B_{j},w)}}^{p} }^{\frac{1}{p}}.
 \end{split}\end{align}
 Lemma \ref{lem:ClassAB} tells us that $w \in
 S_{p,c_{1}}^{N}$. This condition can be applied to the term involving the
 weights to give
 \begin{align*}\begin{split}  
& w (B(x_{j}, 4 \sigma k M^{-\frac{1}{2}}))^{\frac{1}{p}} w^{-\frac{1}{p - 1}}
\br{ B(x_{j}, 4 \sigma k M^{-\frac{1}{2}})}^{\frac{p - 1}{p}} \\
& \qquad \qquad = w
 \br{B_{N}(x_{j}, 4 \sigma k  N^{\frac{1}{2}}
   M^{-\frac{1}{2}})}^{\frac{1}{p}} w^{-\frac{1}{p - 1}} \br{B_{N}(x_{j},
  4 \sigma k N^{\frac{1}{2}}M^{-\frac{1}{2}})}^{\frac{p - 1}{p}} \\
 & \qquad \qquad \lesssim e^{c_{1} 4 \sigma k N^{\frac{1}{2}} M^{-\frac{1}{2}}} \abs{
   B_{N}(x_{j}, 4 \sigma k N^{\frac{1}{2}}M^{-\frac{1}{2}})} \\
 & \qquad \qquad \lesssim k^{d} e^{c_{1}' k},
\end{split}\end{align*}
where $c_{1}' := 4 \sigma c_{1}  N^{\frac{1}{2}} M^{-\frac{1}{2}}$.
Applying this bound to \eqref{eqtn:AboveBelow3},
\begin{equation} 
 \label{eqtn:Const1} 
 \norm{T^{*}_{V,glob}f}_{L^{p}(w)} \lesssim \br{\sum_{j} \br{\sup_{t > 0} \sum_{k = 1}^{\infty}
     \frac{e^{-M t}}{t^{\frac{d}{2}}} e^{-\frac{k^{2}}{4 N t}} e^{c_{1}'
       k} k^{d}  \norm{f}_{L^{p}(4 \sigma k \cdot B_{j},w)}}^{p}}^{\frac{1}{p}} 
 \end{equation}
Define, for $t > 0$ and $k \in \N^{*}$, the quantity
$$
F(k,t) := \frac{e^{-M t}}{t^{\frac{d}{2}}} e^{-\frac{k^{2}}{4 N t}} e^{
  c_{1}' k} k^{d}.
$$
It will be proved that if we set $c_{1}$ small enough then there
must exist some $C, \, \epsilon > 0$ such that
\begin{equation}
  \label{eqtn:ConstUnifEst}
F(k,t) \leq C \cdot e^{- \epsilon k}
\end{equation}
for all $k \in \N^{*}$ and $t > 0$. First note that $k^{d} /
t^{\frac{d}{2}} \lesssim e^{\delta k^{2} / t}$ for any $\delta >
0$. We will therefore have
$$
F(k,t) \lesssim e^{-M t} e^{-\frac{k^{2}}{5 N t}} e^{c_{1}' k}
$$
for all $k \in \N^{*}$ and $t > 0$.

Suppose first that $t \geq k$. Then
$$
F(k,t) \lesssim e^{-M k} e^{c_{1}' k}.
$$
If we let $c_{1}$ be small enough so that $c_{1}' = 4 \sigma c_{1} 
N^{\frac{1}{2}} M^{-\frac{1}{2}} < M$ then
\eqref{eqtn:ConstUnifEst} will hold for this case.

Next, suppose that $t < k$. Then
$$
F(k,t) \lesssim e^{-\frac{k}{5 N}} e^{c_{1}' k}.
$$
Setting $c_{1}$ small enough so that $c_{1}' = 4 \sigma c_{1} N^{\frac{1}{2}}
M^{-\frac{1}{2}}< \frac{1}{5 N}$ will ensure that
\eqref{eqtn:ConstUnifEst} holds for this case as well.

Applying estimate \eqref{eqtn:ConstUnifEst} to
\eqref{eqtn:Const1} gives
\begin{align}\begin{split}
    \label{eqtn:Const2}
 \norm{T^{*}_{V,glob}}_{L^{p}(w)} &\lesssim \br{\sum_{j} \br{ \sum_{k=1}^{\infty} e^{- \epsilon k}
     \norm{f}_{L^{p}(4 \sigma k \cdot B_{j},w)}}^{p}}^{\frac{1}{p}} \\
 &\leq \sum_{k = 1}^{\infty} e^{- \epsilon k} \br{\sum_{j}
   \norm{f}^{p}_{L^{p}(4 \sigma k \cdot B_{j},w)}}^{\frac{1}{p}}.
\end{split}\end{align}
From the bounded overlap property of the balls $B_{j}$, Proposition \ref{prop:Cover}, there exists $N_{1} > 0$ for which
$$
\sum_{j} \norm{f}^{p}_{L^{p}(4 \sigma kB_{j};w)} \lesssim
k^{N_{1}} \norm{f}_{L^{p}(w)}^{p}.
$$
Applying this to \eqref{eqtn:Const2} then gives us our result.
\end{prof}

\subsection{The Harmonic Oscillator}
\label{subsec:Harmonic}

In this section we consider the harmonic oscillator potential, $V(x) =
\abs{x}^{2}$, and prove the second chain of inclusions in Conjecture
\ref{conj:Weak}. We will require the following lemma. It states the exact form of the heat kernel
corresponding to this potential. Its proof can be found in
\cite{simon1979functional} in dimension $d = 1$. Higher dimensions follow
from this case by taking tensor products of Hermite functions.

\begin{lem} 
 \label{lem:Kernel} 
For $t > 0$, define the map $\tilde{k}_{t} : \R^{d} \xx \R^{d} \rightarrow \R$ through
\begin{equation}
\label{eq:kernel}
\tilde{k}_{t}(x,y) =   \frac{1}{\br{2 \pi t}^{d/2}} \exp \br{- \frac{\abs{x - y}^{2}}{2 t}} \cdot \exp \br{- \alpha(t) \br{\abs{x}^{2} + \abs{y}^{2}}},
\end{equation}
where $\alpha$ is defined by
$$
\alpha(t) := \frac{\sqrt{1 + t^{2}} - 1}{2 t}
$$
for all $x, \, y$ in $\R^{d}$ and $t > 0$. The operator $T^{*}_{\abs{x}^{2}}$ is
then given by
$$
T^{*}_{\abs{x}^{2}}f(x) := \sup_{t > 0} \int_{\R^{d}} \tilde{k}_{t}(x,y)
\abs{f(y)} dy
$$
for $f \in L^{1}_{loc}(\R^{d})$ and $x \in \R^{d}$. 
\end{lem}

Let $\tilde{A}^{V,m}_{t,c} = \tilde{A}_{t,c}^{\rho_{V},m}$ and
$\tilde{M}_{V,c}^{m} = \tilde{M}_{\rho_{V},c}^{m}$ be as defined in
Definition \ref{def:AdaptedAveragingH}.

\begin{prop} 
 \label{prop:Harmonic} 
 Let $V(x) = \abs{x}^{2}$. There exists $c, \, m > 0$ such that if
 $\norm{T^{*}_{\abs{x}^{2}}}_{L^{p}(w)} < \infty$ for $1 < p < \infty$ then $w \in H_{p,c}^{V,m}$.
\end{prop}

\begin{prof}  
Fix $1 < p < \infty$ and suppose that the weight $w$ satisfies
$\norm{T^{*}_{V}}_{L^{p}(w)} < \infty$. It must be proved that there
exists $c, \, m > 0$  for which $w \in
H_{p,c}^{V,m}$. This will be accomplished by demonstrating
that there exists $c, \, m > 0$ for which the operator
$$
\tilde{M}_{V,c}^{m}f(x) := \sup_{t > 0} \tilde{A}_{t,c}^{V,m} \abs{f}(x) := \sup_{r > 0}
\frac{1}{\Phi^{V}_{m,c}(x,r)\abs{B(x,r)}} \int_{B(x,r)} \abs{f(y)}
\, dy
$$
is bounded on $L^{p}(w)$. Proposition
\ref{prop:HardyLittlewood2} will then allow us to conclude that $w \in
H_{p,c}^{V,2m(k_{0} + 1)}$.  To prove the
boundedness of $\tilde{M}_{V,c}^{m}$ on $L^{p}(w)$ it is sufficient to
prove the pointwise bound $\tilde{A}_{t,c}^{V,m}\abs{f}(x) \lesssim
T^{*}_{V}f(x)$ for any $t > 0$ and $x \in \R^{d}$. 

\vspace*{0.1in}

Fix $x \in \R^{d}$ and $t > 0$. For $y \in B(x,t)$,
\begin{align*}\begin{split}  
 \alpha(t^{2}) \br{\abs{x}^{2} + \abs{y}^{2}} &\lesssim \alpha(t^{2})
 \br{\abs{x}^{2} + \abs{x - y}^{2}} \\
 &\lesssim \alpha(t^{2}) \br{\abs{x}^{2} + t^{2}}.
\end{split}\end{align*}
It then follows from the simple estimates $\alpha(t^{2}) \lesssim
t^{2}$ and $\alpha(t^{2}) \lesssim 1$,
\begin{align*}\begin{split}  
  \alpha(t^{2}) \br{\abs{x}^{2} + \abs{y}^{2}} &\lesssim t^{2} \abs{x}^{2}
  + t^{2} \\
  &\lesssim t^{2} (1 + \abs{x})^{2} \\
  &\simeq \br{\frac{t}{\rho_{\abs{x}^{2}}(x)}}^{2},
\end{split}\end{align*}
where the last line follows from $\rho_{\abs{x}^{2}}(x) \simeq (1 +
\abs{x})^{-1}$. This proves that there exists $c > 0$ for which
$$
\Phi_{c,2}^{\abs{x}^{2}}(x,t)^{-1} = \exp \br{- c \br{1 + \frac{t}{\rho_{\abs{x}^{2}}(x)}}^{2}}
\leq \exp (- \alpha(t^{2})(\abs{x}^{2} + \abs{y}^{2})).
$$
This estimate combined with the trivial estimate $e^{-\frac{\abs{x -
      y}^{2}}{2 t^{2}}} \simeq 1$ for $y \in B(x,t)$ implies that
\begin{align*}\begin{split}  
 \tilde{A}_{t,c}^{V,2} \abs{f}(x) &= \frac{1}{\Phi_{c,2}^{V}(x,t)\abs{B(x,t)}}
 \int_{B(x,t)} \abs{f(y)} \, dy \\
 &\lesssim \int_{B(x,t)} \frac{1}{t^{d}} \exp \br{-\alpha(t^{2})
   (\abs{x}^{2} + \abs{y}^{2})} \exp \br{-\frac{\abs{x - y}^{2}}{2
     t^{2}}} \abs{f(y)} \, dy \\
 &\lesssim T^{*}_{V}f(x).
\end{split}\end{align*}
Therefore $\norm{\tilde{M}_{V,c}^{2}}_{L^{p}(w)} < \infty$ and
Proposition \ref{prop:HardyLittlewood2} allows us to conclude
that $w \in H_{p, c}^{V,4 (k_{0} + 1)}$.
 \end{prof}

\bibliographystyle{siam}
\bibliography{C:/Users/julian/Desktop/TeX/texmf/bibtex/bibmain} 

\vspace*{0.1in}

\small \textsc{Julian Bailey, School of Mathematics, University of Birmingham, Edgbaston,
  Birmingham, B15 2TT, UK,}

\textit{Email address:} \texttt{j.bailey.1@bham.ac.uk}

\end{document}